\documentclass [12pt] {article}
\title {More 1-cocycles for classical knots}
\author{Thomas Fiedler}
\usepackage{amsmath, amsfonts, amssymb, latexsym, graphicx}
\begin{document}
\newtheorem{proposition}{Proposition}
\newtheorem{theorem}{Theorem}
\newtheorem{lemma}{Lemma}
\newtheorem{corollary}{Corollary}
\newtheorem{example}{Example}
\newtheorem{remark}{Remark}
\newtheorem{definition}{Definition}
\newtheorem{question}{Question}
\newtheorem{conjecture}{Conjecture}
\newtheorem{observation}{Observation}
\maketitle

\begin{abstract}

Let $M^{reg}$ be the topological moduli space of long knots up to regular isotopy, and for any natural number $n > 1$ let  $M^{reg}_n$ be the moduli space of all $n$-cables of framed long knots which are twisted by a string link to a knot in the solid torus $V^3$. We upgrade the Vassiliev invariant $v_2$ of a knot to an integer valued combinatorial 1-cocycle for $M^{reg}_n$ by a very simple formula.  This 1-cocycle depends on a natural number $a \in \mathbb{Z}\cong H_1(V^3;\mathbb{Z})$ with $0<a<n$ as a parameter and we obtain a polynomial-valued 1-cocycle by taking the Lagrange interpolation polynomial with respect to the parameter. We show that it induces a non-trivial pairing on $H_0(M^{reg}_n) \times H_0(M^{reg})$ already for $n=2$.
\end{abstract}

\tableofcontents

\section{Introduction}

In the monograph {\em "Polynomial One-cocycles for Knots and Closed Braids"} \cite{F2} we have laid the foundation of the theory of combinatorial 1-cocycles which depend on integer parameters for knots in the solid torus. In particular, these 1-cocycles give invariants for knots in $\mathbb{R}^3$, often called  {\em classical knots}, when they are evaluated on certain {\em canonical} loops  in the topological moduli spaces of knots in the solid torus, i.e. loops which are universally defined in all connected components of the moduli space. (We will often use \cite{F2} as a reference for definitions, notations and conventions.)

The present paper is the sequel of the monograph and we construct a new combinatorial 1-cocycle.
\vspace{0,4cm}

There is a natural projection $pr: \mathbb{C} \times \mathbb{R} \to \mathbb{C}$ of the 3-space into the plan and knots in 3-space can be given by knot diagrams, i.e. a smoothly embedded oriented circle in 3-space together with its projection into the plan. It is well known that oriented knot types are in one-one correspondence with knot types of oriented long knots. We close now the long knot with the 1-braid to a knot in the standard embedded solid torus $V^3 \subset \mathbb{R}^3$. Moreover, a (blackboard-)framed long knot $K'$ can be replaced by its parallel $n$-cable, called $nK'$ (with the same orientation of all strands), which we close to a knot in $V^3$ by a cyclic permutation braid or more generally by any fixed n-component string link $T$, but which induces a cyclic permutation of its end points. We denote the resulting knot in $V^3$ by $K=T \cup nK'$. The projection into the plan becomes now a projection  into the annulus $pr: \mathbb{C}^* \times \mathbb{R} \to \mathbb{C}^*$. We chose as generator of $H_1(V^3)$ the class which is represented by the closure of the oriented 1-braid. Hence the knot in $V^3$, which is obtained from the $n$-cable, represents the homology class $n \in H_1(V^3)$. We consider the infinite dimensional space $M_n$ of all diagrams of knots $K$ in $V^3$, which represent the homology class $n$, and such that there is a compressing disc of $V^3$ which intersects $K$ transversely and in exactly $n$ points. $M_n$ is called the {\em moduli space of knot diagrams without negative loops in the solid torus}. The knot types in 3-space correspond to the connected components of $M_1$. Given a generic knot diagram $K \subset V^3$ we consider the oriented curve $pr(K)$ in the annulus. A {\em loop in $pr(K)$} is a piecewice smoothly oriented immersed circle in $pr(K)$ which respects the orientation of $pr(K)$. In other words, we go along $pr(K)$ following its orientation and at a double point we are allowed to switch perhaps to the other branch, but still following the orientation of $pr(K)$. Naturally, a loop in $pr(K)$ is called {\em negative} (respectively {\em positive}) if it represents a negative (respectively positive) homology class in $H_1(V^3)$. One easily sees that $pr(K)$ contains only positive loops if and only if $K \subset V^3$ is isotopic to a closed braid with respect to the disc fibration of $V^3$, and that knots which arise as cables of long knots by the above construction, contain never negative loops. The space $M_n^{reg}$ is the subspace of $M_n$ which consists of all those knots on which $pr$ induces an immersion. The space $M_n$ is a natural subspace of the space $M$ of {\em all knots in the solid torus} and which represent the homology class $n$ (without loss of generality we can assume that $n \geq 0$), i.e. negative loops in diagrams are now allowed as well. The natural inclusions $M_n^{reg} \subset M_n \subset M$ are all strict.
\vspace{0,4cm}

Victor Vassiliev has started the combinatorial study of the space of long knots in {\em "Combinatorial formulas of cohomology of knot spaces"} \cite{V2} by studying the simplicial resolution of the discriminant of long singular knots in $\mathbb{R}^3$. Our approach is very different, because it is based on the study of another discriminant, namely the discriminant of non-generic diagrams of knots in the solid torus \cite{F2}. There is also a big difference for $H_1$ of the moduli spaces. Let $K_1 \# K_2$ be the long knot which corresponds to the connected sum of two different knots. In the moduli space of long knots we have a loop which consists of pushing the knot $K_1$ from the left to the right through the knot $K_2$ and then pushing $K_2$ through $K_1$ from the left to the right too \cite{Gr}. In $M_n^{reg}$ we can just push a knot $T \cup nK_1$ through $nK_2$ from the left to the right and then sliding without any Reidemeister moves $T \cup nK_1$ again in its initial position along the remaining part of the solid torus. This is already a loop. We show that our 1-cocycles define in this way  non-trivial {\em pairings on $H_0(M^{reg}_n) \times H_0(M^{reg})$} for $n>1$.
\vspace{0,4cm}

{\em Why do we need to construct 1-cocycles? Classical knots can be transformed into knots in the solid torus. The moduli space of knots in the solid torus is essentially the only moduli space of knots in 3-manifolds (and which are not contained in a 3-ball) with infinite $H_1$, compare \cite{F2} (the moduli space of a long hyperbolic knot deformation retracts onto a 2-dimensional torus \cite{H} and hence the study of $H_1$ will be sufficient). We have to use this, because so far the traditional study of $H_0$ (which has given plenty of invariants)  was not enough to distinguish all classical knots! } 

But why do we need to construct 1-cocycles in a combinatorial way? Of course it would be of great importance for better understanding (and perhaps in order to make connections with string theory) to construct  differential 1-forms on the moduli space, which depend on $a \in H_1(V^3; \mathbb{R})$ as parameter and which represent the same cohomology classes as our combinatorial 1-cocycles if the parameter is in $H_1(V^3; \mathbb{Z})$. But they do not exist yet and we have to put up with the difficult combinatorial approach. 

Combinatorial integer 0-cocycles are usually called {\em Gauss diagram formulas}, compare \cite{PV} and also \cite{F1}. They correspond to finite type invariants and are solutions of the 4T- and 1T-relations. Dror Bar-Natan has shown in {\em "On the Vassiliev knot invariants"} \cite{BN} that such solutions can be constructed systematically by using the representation theory of Lie algebras. Arnaud Mortier has constructed for 1-cocycles of finite type for long knots the analog of the Kontsevich integral  in {\em "A Kontsevich integral of order 1"} \cite{Mor3}. The 4T-relations are now replaced by three 16T- and three 28T-relations and 4x4T-relations! There is actually no representation theory which could help to construct such solutions. The reason for this is simple. The well known representation theory related to the tetrahedron equation, see e.g. \cite{KaKo}, is as usual of a local nature. But $H_1$ of the moduli space of (closed) non-satellite knots in $\mathbb{R}^3$ is only torsion (in contrast to $H_0$), see \cite{H}, \cite{Bu} and \cite{Bu2}, and hence all integer 1-cocycles from local solutions of the tetrahedron equation are trivial (in contrast to the local solutions from the Yang-Baxter equation, see e.g. \cite{K})! We construct therefore in a combinatorial way solutions of the {\em global tetrahedron equation}, i.e. the contribution of a R III move depends on the whole knot in the solid torus and not only on the local picture of the move. This is an equation which is much more complicated as the well known Yang-Baxter equation.
\vspace{0,4cm}

\begin{figure}
\centering
\includegraphics{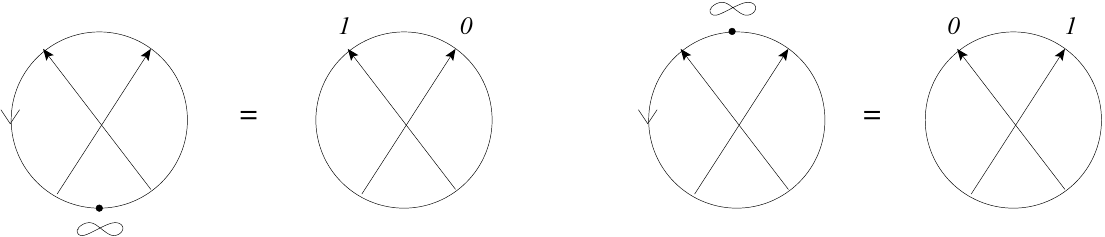}
\caption{\label{PVmark} Polyak-Viro formulas for $v_2(K)$ of long knots}  
\end{figure}

So far, all our 1-cocycles constructed in \cite{F2} have used {\em linear weights} for the contributions of R III moves to the 1-cocycles, i.e. besides the triangle (which corresponds to the move in the Gauss diagram) we consider just the position of individual arrows (which correspond to the crossings) with respect to the triangle. More generally, the construction of  combinatorial 1-cocycles needs as an input {\em Gauss diagram formulas for finite type invariants of long knots}, see e.g.  \cite{PV}, \cite{GPV}, \cite{CP}, \cite{CKR}, \cite{CDM}, \cite{F1}. In the present paper we make use of the Polyak-Viro formula for $v_2(K)$ of long knots shown on the left in 
Fig.~\ref{PVmark}. 

{\em We use the Polyak-Viro formula in order to define an integer combinatorial 1-cocycle in $M_n^{reg}$, called  $R_a^{(2)}$.}

The natural number $a \in \mathbb{Z}\cong H_1(V^3;\mathbb{Z})$ with $0<a<n$ is a parameter for the 1-cocycle. 

The construction of the 1-cocycle is complex but the outcome is a rather beautiful formula. Surprisingly, it becomes now really essential that the projection of $T \cup nK'$ into the annulus contains no loops which represent a negative homology class and that the long knots are framed. (Usually, the invariance of quantum knot invariants under the R I moves is just a matter of normalization. This is no longer the case for 1-cocycles on 1-parameter families of diagrams: the place of the R I move in the diagram becomes very important.)

\begin{definition}
Let $T \cup nK'$ be a knot in $M_n$. We consider the cyclic $n$-fold covering of $V^3$ and we lift $T \cup nK'$ to a knot $\tilde K$ in the covering. $\tilde K$ can be naturally identified with a long knot and we call it the {\em underlying long knot} for $T \cup nK'$.

\end {definition}

It turns out that $R_a^{(2)}$ is a 1-cocycle  even in $M_n$ (and not only in $M_n^{reg}$) if and only if for the underlying long knot $v_2(\tilde K)=0$. Notice, that this can always easily be achieved by putting an appropriate number of trefoils or figure-eight knots in small 3-balls on a branch of $T$. 
\vspace{0,4cm}

To obtain our combinatorial 1-cocycles there are only three types of equations to solve, and in the order given below, because we consider only regular isotopy and only R III moves contribute to the 1-cocycle (compare \cite{F2}):

(a) {\em the commutation relations}: a R III move with simultaneously another Reidemeister move;

(b) the {\em positive global tetrahedron equation}: going around a positive quadruple crossing in the moduli space;

(c) the {\em cube equations}: going around triple crossings where two branches are ordinary tangential.

\vspace{0,4cm}

We explain now briefly our method. We want to construct a {\em weight} $W_2(p)$ for each R III move $p$ which is defined by using couples of crossings (i.e. arrows in the Gauss diagram which always go from the under-cross to the overcross, compare \cite{F2}) outside of the R III move $p$ (i.e. arrows not of the triangle in the Gauss diagram which corresponds to the R III move). But these couples have to be related to the R III move $p$. The commutation relations (a) force the weight $W_2(p)$ to be a knot invariant if we consider all couples, and not only those which are related to the move $p$. R I moves force the weight not to contain isolated arrows with homological markings $0$ or $n$. R II moves force that each crossing of the weight contributes with its sign (writhe). Consequently, it is then sufficient to show that the weight $W_2(p)$ is invariant under a simultaneous R III move with the move $p$.

\begin{figure}
\centering
\includegraphics{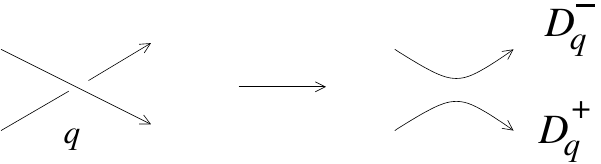}
\caption{\label{splitq} The two ordered knot diagrams associated to a crossing $q$}  
\end{figure}

An invariant of long knots which is given by a Gauss diagram formula, which uses only couples of arrows, is an invariant of degree 2. Hence, it can be defined by the beautiful formulas of Michael Polyak and Oleg Viro for $v_2(K)$ of long knots \cite{PV}, see Fig.~\ref{PVmark}. The point on the circle corresponds to the point at infinity on the knot. For a crossing $q$ we call $D_q^+$ the knot which is obtained by smoothing the crossing $q$ from the under-cross to the over-cross (and the remaining knot is called $D_q^-$), see Fig.~\ref{splitq}.  The homological marking of $q$ is the homology class in $\mathbb{Z}\cong H_1(V^3)$ represented by $D_q^+$ (compare \cite{F2}), i.e. by the knot which corresponds to the arc of the circle from the over-cross to the under-cross. Consequently, the point at infinity is in $D_q^+$ if and only if the homological marking $[q]=1$. The sign of a crossing $q$ is denoted as usual by $w(q)$.

\begin{proposition}
Let $G$ be a Gauss diagram formula for long knots which defines a knot invariant. Let $K \subset V^3$ be a knot which belongs to $M_n$ (i.e. no negative loops). We keep the markings $0$ in $G$ but we replace each marking $1$ by the marking $n$. Then the resulting Gauss diagram formula defines a knot invariant for $K$ in $M_n$.

\end{proposition}

{\em Proof.}

We consider the cyclic $n$-fold covering $V^3$ over $V^3$ and we lift $K$ to the underlying long knot $\tilde K$ in it. Crossings with marking $0$ stay crossings with marking $0$, crossings with marking $n$ become crossings with marking $1$ and all other crossings disappear. The fact that $G$ is an invariant for long knots implies the result.
 $\Box$

\begin{figure}
\centering
\includegraphics{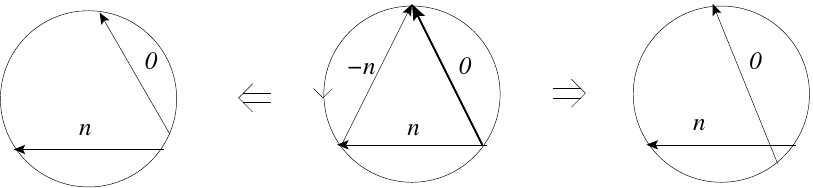}
\caption{\label{PVM} The Polyak-Viro formula for $v_2(K)$ is not true in $M$}  
\end{figure}

Notice that the result is no longer true in $M$ already for the invariant of order two (because more markings appear for $\tilde K$), see Fig.~\ref{PVM}.
\vspace{0,4cm}

Let's come back to the formula for $v_2$, which is important in this paper.

\begin{figure}
\centering
\includegraphics{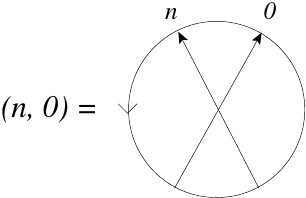}
\caption{\label{weightn} Configuration (or arrow diagram) $(n,0)$}  
\end{figure}

We will use for $W_2(p)$ the following configuration, which corresponds to the Polyak-Viro formula on the left in Fig.~\ref{PVmark} (where as usual we take the product of the signs of the two crossings), and which we denote shortly by $(n,0)$, see Fig.~\ref{weightn}.

To each Reidemeister move of type III corresponds a diagram with a {\em triple 
crossing} $p$: three branches of the knot (the highest, middle and lowest with respect to the projection $pr: \mathbb{C}^* \times \mathbb{R} \to \mathbb{C}^*$) have a common point in the projection into the 
plane. A small perturbation of the triple crossing leads to an ordinary diagram with three crossings near $pr(p)$.

\begin{figure}
\centering
\includegraphics{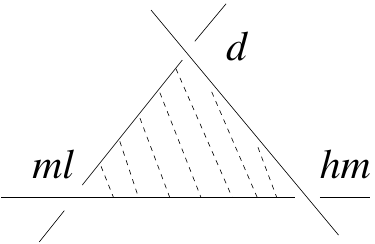}
\caption{\label{names} The names of the crossings in a R III-move}  
\end{figure}

\begin{figure}
\centering
\includegraphics{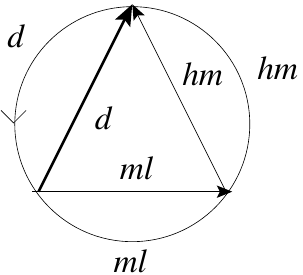}
\caption{\label{namesarc} The names of the arcs in a R III-move}  
\end{figure}

\begin{definition}
 We call the crossing between the 
highest and the lowest branch of the triple crossing $p$ the {\em distinguished crossing} of $p$ and we denote it by $d$ ("$d$" stands for distinguished). The crossing between the highest branch and the middle branch is denoted by $hm$ and that of the middle branch with the lowest is denoted by $ml$,compare Fig.~\ref{names}. For better visualization we draw the crossing $d$ always with a thicker arrow. Moreover, we give the same name to the adjacent arc of the crossing in the circle, compare Fig.~\ref{namesarc}

\end{definition}

Let $p$ be a triple crossing. The triangle in the Gauss diagram cuts the circle into three arcs, which we denote by the name of the adjacent crossing in the triangle. We denote the union of the {\em open} arcs $d$ and $ml$ (to which we add just the common foot of $d$ and $ml$) by $dml$. Of course, we could identify $dml$ with $D_{hm}^+$.

\begin{definition}
Let $p$ be a triple crossing. An ordinary crossing $q$ of the diagram is called an {\em $f$-crossing for $p$} if $[q]=n$ and the foot of $q$ is in the open arc $dml$. An ordinary crossing $q$ is called an $n$-crossing respectively $0$-crossing if its homological marking $[q]=n$ respectively $[q]=0$.

\end{definition}
 
 We make a connection of the weight with the triple crossing $p$ by defining 
 the {\em weight of $p$ of order 2} as $W_2(p)= \sum (n,0)$, where the sum is taken only over all $f$-crossings for $p$ (and no restrictions on the foots of the $0$-crossings here).

\vspace{0,4cm}

It turns out that the weight $W_2(p)$ solves the commutation relations (a), but it does not solve the positive global tetrahedron equation (b). The eight strata of R III moves in the meridian of a positive quadruple crossing in the moduli space come in pairs with different signs: $P_i$ and $\bar P_i$ (compare \cite{F2}). They can differ by $n$-crossings for which the foot has slide over the head or the foot of a crossing in the triangle of $P_i$ or $\bar P_i$ and hence sometimes these  crossings do no longer contribute to $W_2(p)$. The key point is that we consider only R III moves of a very particular global type (see the next section). It turns then out, that a $n$-crossing which changes from $P_i$ to $\bar P_i$ the arc of its foot, is always a crossing $hm$ for another R III move $P_j$ and $\bar P_j$ of the same  particular global type in the meridian of the quadruple crossing. We use this to define a weight of order 1 for $P_j$ and $\bar P_j$, which we multiply by a certain refined linking number and which will contribute to the 1-cocycle too. The weight of order 1 is essentially the same for $P_j$ and $\bar P_j$ but the refined linking numbers are different. This is really complex but leads finally to a solution of the positive global tetrahedron equation (b). The solution of the cube equations (c) is then relatively easy and leads just to some simple linear correction terms in the 1-cocycle.

Let us mention that the construction of  $R_a^{(2)}$ completely breaks down for the parameters $a=0$ and $a=n$ (and hence we can not apply it for $n=1$, and in particular it becomes essential again that knots have to be framed, because the isotopy class of cables depends on the framing).
We will give the precise result in the next section. 

Notice that we have made several choices. The two Polyak-Viro formulas for $v_2$ lead to four different "dual" 1-cocycles: head instead of foot for $n$-crossings, foot or head in $dhm$ instead of $dml$. If we replace a Gauss diagram formula by its mirror image then (1) and (2) are preserved and in (3) foots are replaced by heads. Also the group $\mathbb{Z}/2\mathbb{Z} \times \mathbb{Z}/2\mathbb{Z}$ acts naturally on the moduli spaces, generated by orientation reversing of $K$ together with the hyper elliptic involution of the solid torus, and by taking the mirror image. Some of the corresponding 1-cocycles  correspond simply to the different choices in the construction. 
\vspace{0,4cm}

In this paper we concentrate on the construction of the new 1-cocycle. Calculations by hand of interesting examples of the pairing (and of the values of the 1-cocycle on the important  Fox-Hatcher loops which make relations to geometry, compare \cite{F2}) are too complicated because we need that $n \ge 2$. We calculate just the easiest examples in order to show that it is not always trivial already for $n=2$.

{\em  In common work with Roland van der Veen and Jorge Becerra we hope to create a computer program and to calculate lots of examples and to use it to answer many open questions in knot theory in 3-space.}

But our 1-cocycles are not only useful in order to define knot invariants. A well known question in 3-dimensional knot theory is: given a complicated knot diagram, is this the unknot, or more generally, do two given knot diagrams represent the same knot type?

{\em An analogue of this question on the next level is: given two loops in the moduli space $M_n^{reg}$ by complicated movies of diagrams, are these loops homotopic? (There isn't known any algorithm in order to simplify the loops.)  However, our 1-cocycles allow one sometimes to give a negative answer to this question, even without recognizing the knot type.}

\section{The integer 1-cocycle $R_a^{(2)}$}

Let $\gamma$ be a generic oriented loop in $M$. 

\begin{figure}
\centering
\includegraphics{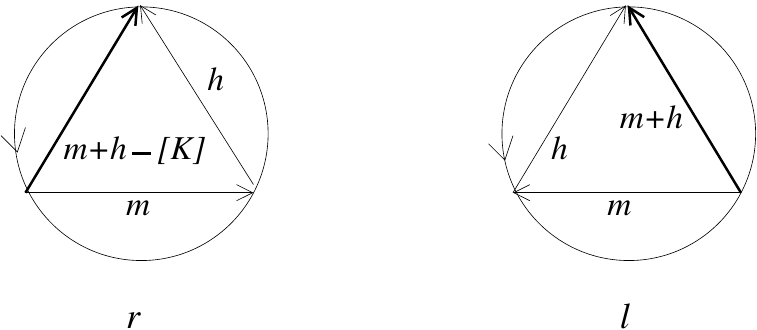}
\caption{\label{RVIII} The global types of R III moves for knots in the solid torus}  
\end{figure}

A Reidemeister III move $p$ in $\gamma$ corresponds to a triangle in the Gauss diagram. The {\em global type of a Reidemeister III move} is now shown in Fig.~\ref{RVIII}, where $m+h$, $m+h-[K]$, $m$ and $h$ are the homological markings of the corresponding arrows. Here, $[K]=n$ is the homology class represented by the knot $K$. Moreover, we indicate whether the arrow $ml$ in the triangle goes to the left, denoted by $l$, or it goes to the right, denoted by $r$. Notice that each two of the markings determine always the third marking. 

{\em It is convenient to encode the global type of a R III move in the following way: $r([d],[hm],[ml])$ and respectively $l([d],[hm],[ml])$.}

We will consider Gauss diagram formulas for R III moves. The arrows in the configurations which are not arrows of the triangle are called the {\em weights} of the formulas.

\begin{figure}
\centering
\includegraphics{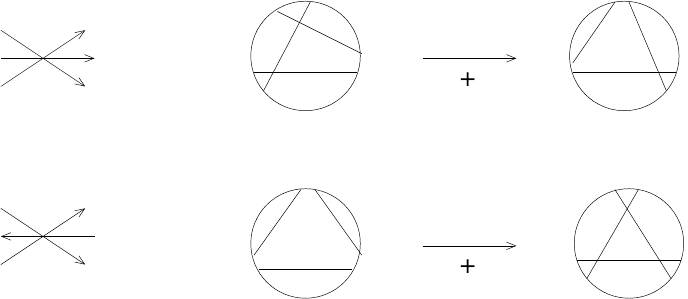}
\caption{\label{coorient} The coorientation for Reidemeister III-moves}  
\end{figure}

\begin{definition}
The {\em coorientation} for a Reidemeister III move  is the direction from two intersection points of the corresponding three arrows to one intersection point and of no intersection point of the three arrows to three intersection points, compare 
Fig.~\ref{coorient}. (In \cite{F2} we have studied the discriminant $\Sigma$ of non generic diagrams in $M$. In particular, we have shown in the cube equations for $\Sigma^{(2)}_{trans-self}$, that the two coorientations for triple crossings fit together for the strata of $\Sigma^{(1)}_{tri}$ which come together in $\Sigma^{(2)}_{trans-self}$.) Evidently, our coorientation is completely determined by the corresponding planar curves and therefore we can draw just chords instead of arrows in Fig.~\ref{coorient}.  We call the side of the complement of the codimension 1 strata of R III moves $\Sigma_{tri}^{(1)}$ in $M$ into which points the coorientation, the {\em positive side} of  $\Sigma_{tri}^{(1)}$.

\end{definition}

Each transverse intersection point $p$ of an oriented generic arc in $M$ with $\Sigma^{(1)}_{tri}$ has now an intersection index $+1$ or $-1$, called $sign(p)$, by comparing the orientation of the arc with the coorientation of $\Sigma_{tri}^{(1)}$.

\begin{figure}
\centering
\includegraphics{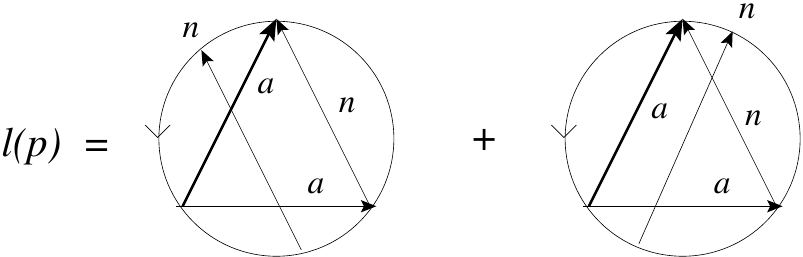}
\caption{\label{lofp} The refined linking number $l(p)$ for the type $r(a,n,a)$}  
\end{figure}

\begin{figure}
\centering
\includegraphics{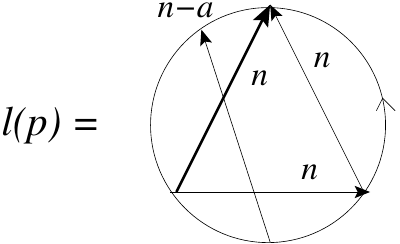}
\caption{\label{lnnn} The refined linking number $l(p)$ for the type $r(n,n,n)$}  
\end{figure}

Let $0<a<n$ be fixed. We consider in $R_a^{(2)}$ {\em only} R III moves $p$ of the global types $r(a,n,a)$ and $r(n,n,n)$, compare e.g. Fig.~\ref{lofp} and Fig.~\ref{lnnn}.

\begin{definition}

Let $p$ be of type $r(a,n,a)$. Then the {\em weight $W_2(p)$} is defined by applying $(n,0)$ to the underlying long knot but only by summing up over those configurations were the $n$-crossing is an $f$-crossing for $p$ (i.e. an $n$-crossing with the foot in the open arc $dml$). 

The {\em weight $W_2(hm)$} is defined by applying $(n,0)$ to the underlying long knot but only by summing up over those configurations were the $n$-crossing is the crossing $hm$ in $p$.

The {\em refined linking number $l(p)$} is defined as the sum of the signs of all $n$-crossings which cut the crossing $ml$ of $p$ (i.e. the corresponding arrows in the Gauss diagram intersect) and which have their foot in the open arc $ml$, compare Fig.~\ref{lofp}.

\vspace{0,4cm}

Let $p$ be of type $r(n,n,n)$. There is no weight $W_2(p)$ defined but only the weight $W_2(hm)$, and which is defined exactly in the same way as for $p$ of type $r(a,n,a)$.

The {\em refined linking number $l(p)$} is defined now as the sum of the signs of all $(n-a)$-crossings (i.e. $[q]=n-a$) which cut the crossing $ml$ of $p$ (i.e. the corresponding arrows in the Gauss diagram intersect) and which have their foot in the open arc $ml$, compare Fig.~\ref{lnnn} (there is only one possibility to cut the triangle for the $(n-a)$-crossings here, because there are no negative loops in the diagrams).
\vspace{0,4cm}

For a given crossing $c$ different from $hm$ we denote by $W_2(c)$ the sum of all $(n,0)$ where $c$ is the $n$-crossing. 

\end{definition}

The following is our central definition.

\begin{definition}

Let $\gamma$ be a generic oriented loop in $M_n$.

The integer-valued 1-cochain $R_a^{(2)}$ is defined by 

\begin{align*}
& R_a^{(2)}(\gamma)= \sum_{p=r(a,n,a) \in \gamma}sign(p)(W_2(p)+(l(p)+w(hm)-1)W_2(hm)w(hm)) \\
&- \sum_{p=r(n,n,n) \in \gamma}sign(p)l(p)W_2(hm)w(hm).\\
\end{align*}

\end{definition}

We want to apply the 1-cochains to canonical loops in $M_n^{reg}$. In fact, it turns out that we can slightly weaken the condition of regular isotopy to {\em semi-regular isotopy}: in a R I move a new crossing appears or disappears which has marking $0$ or $n$. In a semi-regular isotopy we allow R I moves with the marking $0$, but we do not allow R I moves with the marking $n$ and we call the corresponding moduli space $M_n^{semi-reg}$. Notice that for semi-regular isotopy of a long knot $K$ we have a finite type invariant $w_1(K)$ of order 1, namely the algebraic sum of all crossings in $K$ with marking $1$.

\begin{definition}
Let $K'$ be an oriented  framed long knot and let $T$ be a string link which induces a cyclical permutation of its end points. We denote by $K=T \cup nK'$ the knot in $M_n^{reg}$ which is obtained by closing the parallel $n$-cable of $K'$ (with respect to the framing and with the induced orientation) by the string link $T$ to a knot in the solid torus. The loop $push(T, nK')$ is defined by pushing once $T$ through the parallel n-cable of $K'$ in the solid torus in counter-clockwise direction \cite{F2}.

\end{definition}

\begin{figure}
\centering
\includegraphics{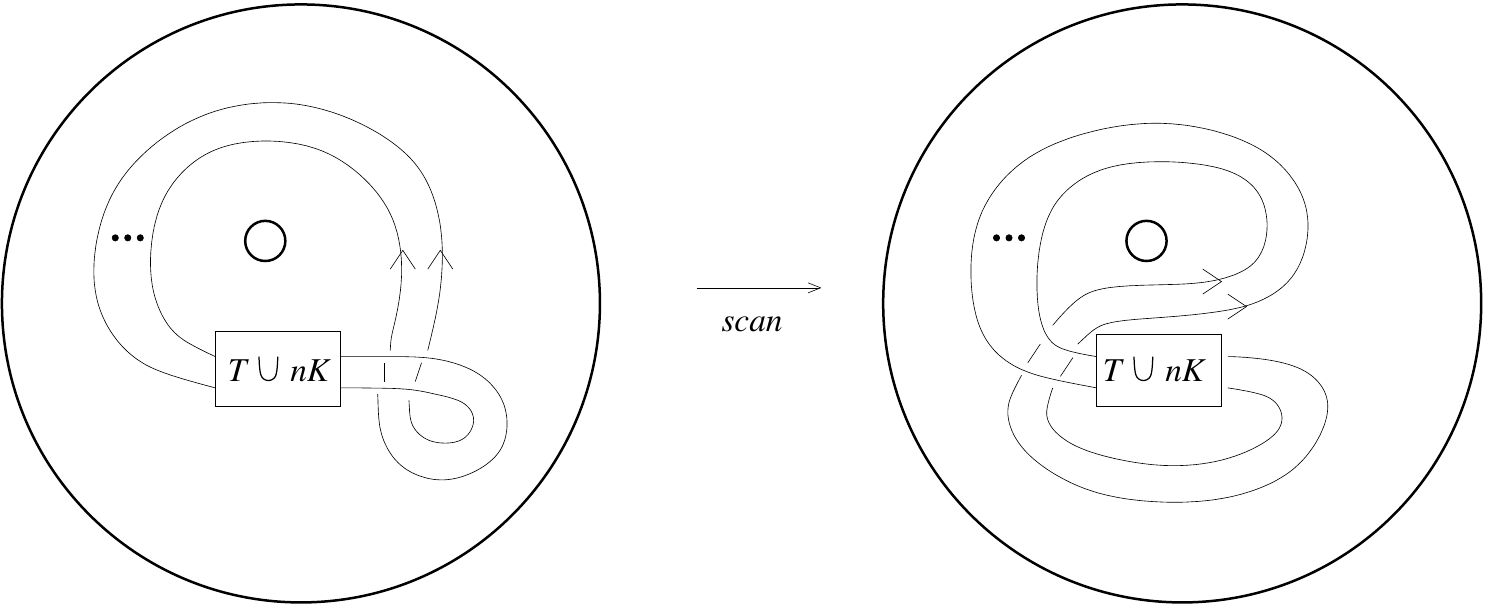}
\caption{\label{scannK} The scan of $T \cup nK'$}  
\end{figure}

Adding to $T \cup nK'$ a full-twist in form of a positive $n$-curl, compare Fig.~\ref{scannK}, we could push then $T \cup  nK'$ once through the curl. This is an example for our pairing, where the long knot is just a positive curl with $w_1=0$. On the other hand it is a nice representative in $M_n^{reg}$ of Gramain's loop $rot(K=T \cup nK')$, see \cite{Gr}, which is induced by the full rotation of $V^3$ around its core. We call the first half of the loop $rot$ the {\em scan-arc} and denote it by $scan(T \cup nK')$.

A picture of $scan(T \cup nK')$ is given in Fig.~\ref{scannK}. Notice, that it is only an arc in $M_n^{reg}$ and not a loop.

(Instead of classical knots our invariants are of course also defined for arbitrary coherently oriented $n$-component string links $T$ in $\mathbb{R}^3$, by evaluating the 1-cocycles e.g. on $scan(T)$ or on $rot(T)$.)

One easily sees that the $n$-curl is semi-regular isotopic to the full twist $\Delta^2$ (the generator of the center of the braid group) in the braid group $B_n$. Hence, if $T$ is a braid then we could push also the $n$-curl  through $T \cup nK'$. Geometrically, this corresponds to the loop which is generated by the full rotation of the companion solid torus $V_1^3 \subset V^3$ (from the cabling construction) around its core.
\vspace{0,4cm}

The following theorem is our main result.

\begin{theorem}
$R_a^{(2)}$ is an integer 1-cocycle in $M_n^{semi-reg}$ for each $n>1$, each $0<a<n$ and which represents a non-trivial cohomology class. In particular, the corresponding pairings $H_0(M^{semi-reg}_n) \times H_0(M^{reg})$ are non-trivial already for $n=2$ and when $T$ in $T \cup 2K'$ is just the standard generator $\sigma_1$ of the braid group $B_2$.

Moreover, already $R_a^{(2)}(scan(T \cup nK'))$ is an invariant of $T \cup nK'$ up to semi-regular isotopy (and hence a knot invariant of $K'$).

$R_a^{(2)}$ is a 1-cocycle  even in $M_n$ if and only if for the underlying long knot $v_2(\tilde K)=0$.

$R_a^{(2)}$ depends in general non-trivially on the integer parameter $0<a<n$ and hence it can be made a polynomial-valued 1-cocycle by taking the Lagrange interpolation polynomial with respect to the parameter.

\end{theorem}

Notice that all the knot invariants from Theorem 1 can be calculated with polynomial complexity with respect to the number of crossings of $T \cup nK'$.

Let $T \in M^{semi-reg}_n$ be the parallel $n$-cable of a long framed knot $K'$ with the strands twisted by the permutation braid $\sigma_1 \sigma_2 ... \sigma_{n-1}$ and let $K'' \in M^{reg}$.

We have a pairing on $\oplus_{n >1} ( H_0(M^{semi-reg}_n) \times H_0(M^{reg}))$ induced by 

$\oplus_n \oplus_a  (R_a^{(2)}(push(T, nK'')))$. Is it symmetric with respect to $K'$ and $K''$?
\vspace{0,4cm}

It is worth noting that the Polyak-Viro formula for $v_2$ is just the first case of the Chmutov-Khoury-Rossi Gauss diagram formulas for the coefficients of the Conway polynomial. For the convenience of the reader we repeat these formulas here, compare \cite{CKR}.

Let $k \in \mathbb{N}$ be fixed and let $A_{2k}$ be an arrow diagram (i.e. an abstract Gauss diagram without signs on the arrows, we call it often a configuration) with one (oriented) circle, $2k$ arrows and a base point.  $A_{2k}$ is called {\em ascending one-component} if by going along the oriented circle starting from the base point and each time jumping along the arrow, we meet each arrow first at its foot and the traveling meets the whole circle.

\begin{figure}
\centering
\includegraphics{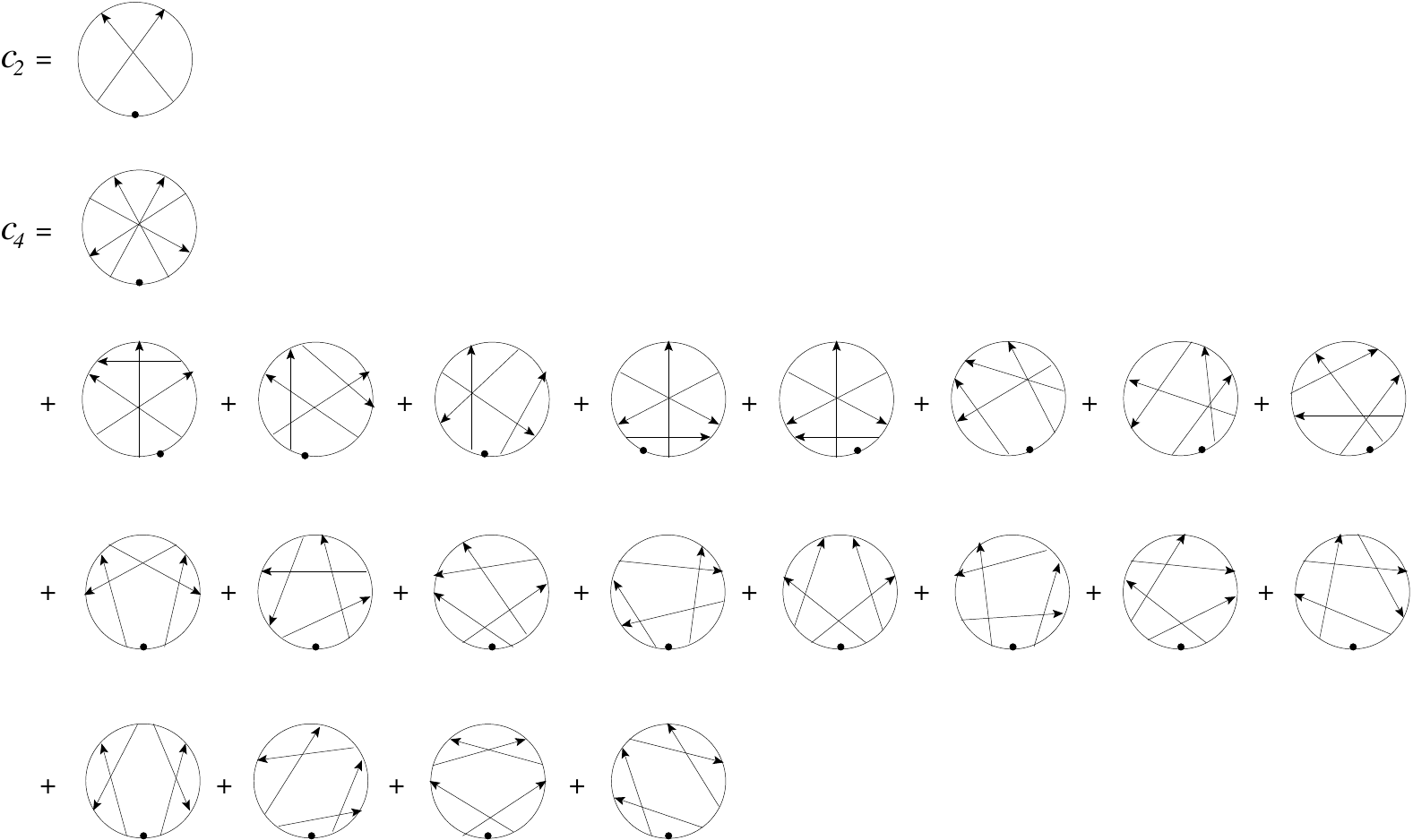}
\caption{\label{ckr} The Chmutov-Khoury-Rossi formulas for $k=1$ and $k=2$}  
\end{figure}

The Gauss diagram formula $C_{2k}$ is simply the sum of {\em all} ascending one-component arrow diagrams $A_{2k}$. We show as examples $C_2$ and $C_4$ in Fig.~\ref{ckr}.

Chmutov, Khoury and Rossi have proven by using skein relations that $C_{2k}$ evaluated on a long knot diagram $K$ is the coefficient of $z^{2k}$ in the Conway polynomial of $K$.  Brandenbursky \cite{Br} has interpreted $C_{2k}$ evaluated on a long knot diagram $K$ as the algebraic number of certain surfaces with one boundary component related to $K$ and of Euler characterisics $1-2k$ (in the spirit of Gromov-Witten invariants), and he has given a direct proof that it is invariant under the Reidemeister moves (and in particular, there is a natural $1-1$ correspondence of contributing configurations before and after the R III move). Hence the formula $C_{2k}$ gives even an invariant for virtual long knots (this is very important for us in the proof, even if our 1-cocycles are only well-defined for real knots, because each loop for virtual knots would contain forbidden moves).
\vspace{0,4cm}

\begin{figure}
\centering
\includegraphics{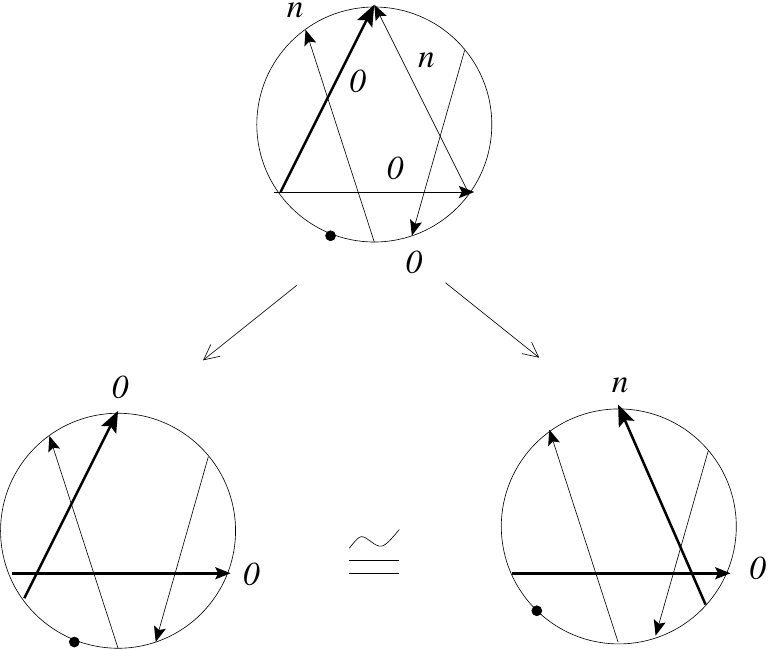}
\caption{\label{vanishn}  A couple $(0,0)$ is replaced by a couple $(n,0)$}  
\end{figure}

A natural first attempt to generalize our result is to define the weights by demanding that all $n$-crossings in all configurations of $C_{2k}$ have to be $f$-crossings. However, this does not work already for the commutation relations.  Indeed, let $p$ be a move of type $r(a,n,a)$ and let $p'$ be another independent R III move. If $p'$ is of type $r(0,n,0)$ and of local type $1$ (i.e. all three crossings positive) then the $n$-crossing can appear or disappear in the corresponding configurations before and after the move $p'$ in the $1-1$ correspondence. If the $n$-crossing would be not an $f$-crossing then it could be that the configuration without the $n$-crossing contributes to $W_{2k}(p)$ but the corresponding configuration with the $n$-crossing does not. We show an example for $R_a^{(4)}$ for the Conway polynomial in Fig.~\ref{vanishn} (the two configurations correspond to the second and the fourth configurations in the third line of Fig.~\ref{ckr}).

\vspace{0,4cm}
The Chmutov-Khoury-Rossi formulas have two remarkable {\em properties}: 

$1)$  if we travel on the oriented circle starting from the point at infinity then we meet in each configuration always first the foot of an $n$-crossing and we meet last the foot of an $0$-crossing (this is an immediate consequence of their definition). We call this correspondingly the {\em first $n$-crossing} and the {\em last $0$-crossing} in each configuration. 

$2)$ there is a 1-1 correspondence of contributing configurations before and after a R III move. If all crossings in the R III move are positive and a couple of  intersecting crossings in a configuration (i.e. the arrows intersect) contributes before the move (no matter the direction of the move) then it is replaced by exactly one different couple but also with intersecting crossings in a configuration after the move (this is a consequence of Brandenbursky's proof \cite{Br}). Consequently, the corresponding analogue is true as well for couples of non-intersecting crossings.
\vspace{0,4cm}

It turns out that the following is the right definition: {\em we fix a point on the knot in the disc at infinity.  $W_{2k}(p)$  is defined by demanding simply that the first $n$-crossing in each configuration has to be an $f$-crossing (i.e. it has its foot in the open arc $dml$) and $W_{2k}(hm)$ is defined by demanding that $hm$ is the first $n$-crossing in each configuration. The definitions of the refined linking numbers $l(p)$ stay the same. Definition 6 together with these changings gives now a lift of all the Chmutov-Khoury-Rossi-Brandenbursky formulas to non-trivial combinatorial 1-cocycles! But the proof is much more complicated.  
Very recently we have succeeded to construct yet another much more complicated 1-cocycle with quadratic weights by demanding in addition, that the last $0$-crossing in each configuration should has its foot in the open arc $dml$ too. We will come back to these constructions in another paper or perhaps in another book.}

\section{Examples}
 
Let $K$ be the positive long trefoil with $w_1(K)=1$, let $n=2$ and let $T= \sigma_1$. We denote by $2K$ the blackboard framed 2-cable of the long knot $K$ and we consider $R_1^{(2)}(push(\sigma_1,2K))$.

\begin{figure}
\centering
\includegraphics{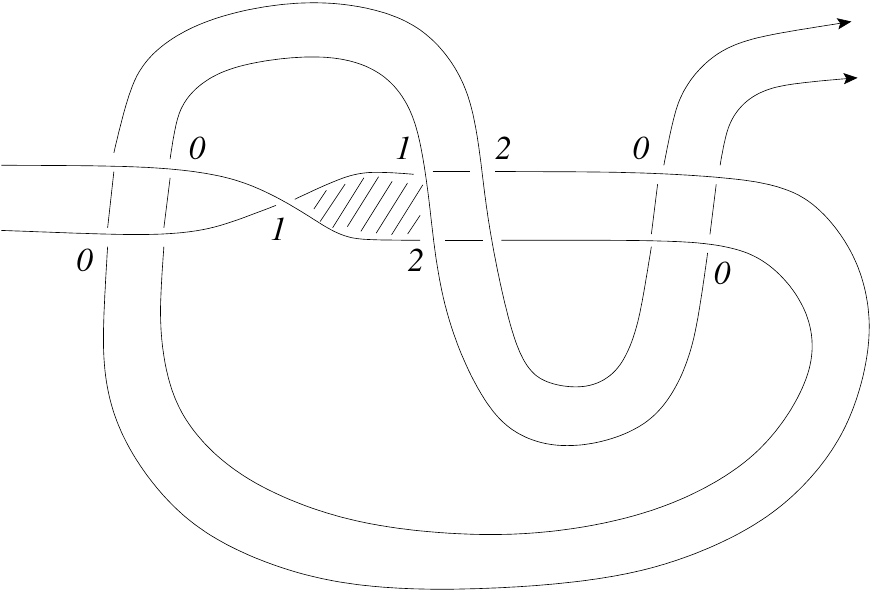}
\caption{\label{trefoil} The single R III move which contributes to $R_1^{(2)}(push(\sigma_1,2K))$ for the right trefoil}
\end{figure}

\begin{figure}
\centering
\includegraphics{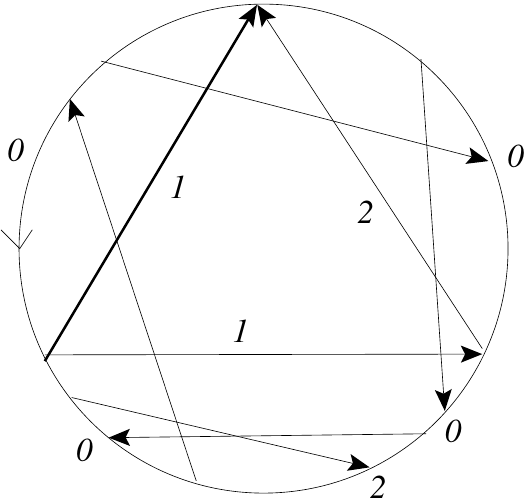}
\caption{\label{gausstrefoil} The Gauss diagram of the single R III move which contributes to $R_1^{(2)}(push(\sigma_1,2K))$ for the right trefoil}
\end{figure}

If $T$ is a braid, then it has no $n$-crossings and consequently there are no R III moves of type $r(n,n,n)$ at all in the loop $push(T, nK)$. One easily sees that there is exactly one R III move of global type $r(a,n,a)$ in the loop and it has a positive sign, compare
Fig.~\ref{trefoil}, where we show the diagram just before the move. The Gauss diagram with markings of the move (without the crossings of marking $a=1$ which are not in the triangle) is shown  in Fig.~\ref{gausstrefoil}. Notice that the homological markings $0$ and $n$ always determine the point at infinity (up to symmetry of the configuration) and hence we can calculate the contributions of the weights without going to the cyclic $n$-fold covering of the solid torus.

\begin{figure}
\centering
\includegraphics{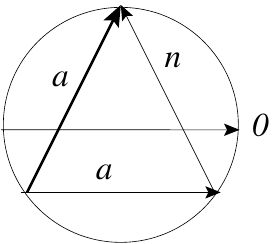}
\caption{\label{restcross} The only possible $0$-crossings in the weight $W_2(hm)$}  
\end{figure}

\begin{figure}
\centering
\includegraphics{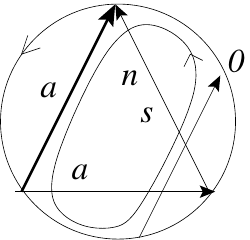}
\caption{\label{negloop} The $0$-crossing would imply the existence of a loop $s$  with $[s]=0-a<0$}  
\end{figure}

Actually, the possible positions of the $0$- and $n$-crossings in the weights $W_2$ are often determined by the fact that there are no negative loops in diagrams for knots in $M_n$. We show examples in Fig.~\ref{restcross} and Fig.~\ref{negloop}. 

The single 2-crossing in our example is an $f$-crossing (because its foot is in $dml$). It follows that $W_2(p)=1$, $l(p)=0$, $w(hm)-1=0$ and $W_2(hm)=1$. Consequently, $R_1^{(2)}(push(\sigma_1,2K))=1$. Notice, that this is of course a knot invariant of $K$. Indeed, for another knot $K'$ we can add small curls such that $w_1(K)=w_1(K')$. The corresponding compact knots $K$ and $K'$ are isotopic in $S^3$ if and only if the long knots $K$ and $K'$ are semi-regular isotopic as long knots, compare e.g.  \cite{F1}.
\vspace{0,4cm}

\begin{figure}
\centering
\includegraphics{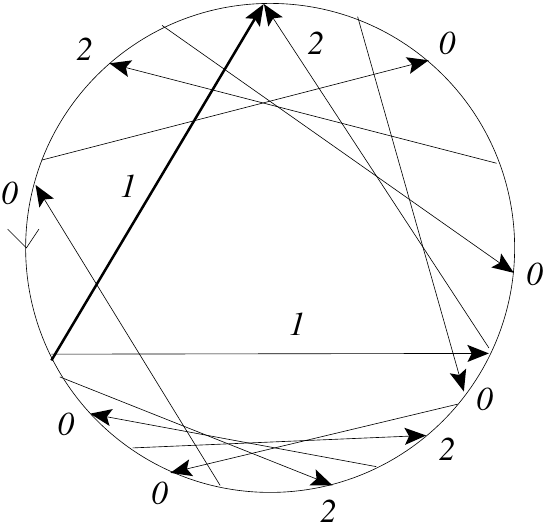}
\caption{\label{gaussfive} The Gauss diagram of the first R III move which contributes to $R_1^{(2)}(push(\sigma_1,2K))$ for the torus knot $(2,5)$}
\end{figure}

\begin{figure}
\centering
\includegraphics{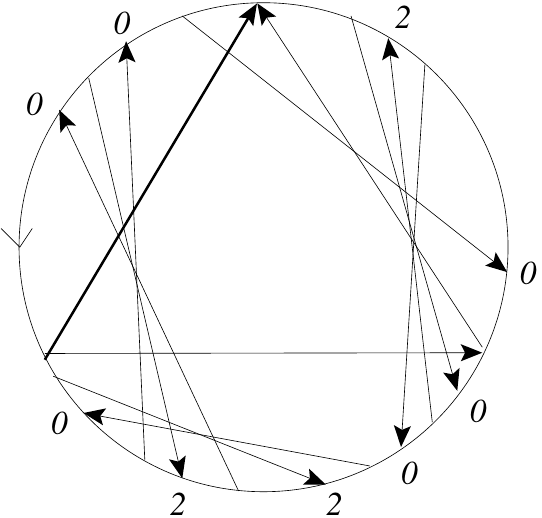}
\caption{\label{gaussfive2} The Gauss diagram of the second R III move which contributes to $R_1^{(2)}(push(\sigma_1,2K))$ for the torus knot $(2,5)$}
\end{figure}

Let $K$ be the positive torus knot of type $(2,5)$ with $w_1(K)=2$ and again let $n=2$ and let $T= \sigma_1$. In analogy to the positive trefoil there are now exactly two moves of type $r(a,n,a)$ in the loop and they have both a positive sign. We show their Gauss diagrams in  Fig.~\ref{gaussfive} and  Fig.~\ref{gaussfive2}. In the first move $l(p)=0$, $W_2(p)= 2+1$. In the second move $l(p)=1$, $W_2(p)=2+2+1$, $W_2(hm)=1$. Consequently, $R_1^{(2)}(push(\sigma_1,2K))=9$. 
\vspace{0,4cm}

Let us  take the same examples but for $n=4$ and $T=\sigma_1 \sigma_2 \sigma_3$. One easily sees that if we change $\sigma_2$ in $T$ to $\sigma_2^{-1}$ then  $R_1^{(2)}(push(T,4K))$ and  $R_3^{(2)}(push(T,4K))$ stay the same as before  but  $R_2^{(2)}(push(T,4K))$ (which was positive) changes the sign, because the sign of the corresponding move (which is now a different one) has changed, but the homological markings of the crossings $\sigma_2$ and $\sigma_2^{-1}$ stay the same, namely $2$. This simple example shows already that  $R_a^{(2)}(push(T,4K))$ depends non-trivially on the parameter $a$.
\vspace{0,4cm}

The calculation of $R_1^{(2)}(rot(\sigma_1 \cup 2K))$ (which can be seen as a particular case of our pairing) is already more complicated. We have calculated that if $K$ is the positive trefoil  with $w_1(K)=1$ then $R_1^{(2)}(rot(\sigma_1 \cup 2K))=R_1^{(2)}(scan(\sigma_1 \cup 2K))=-2$, and if $K$ is the figure eight knot with $w_1(K)=-1$ then $R_1^{(2)}(rot(\sigma_1 \cup 2K))=R_1^{(2)}(scan(\sigma_1 \cup 2K))=-6$.
\vspace{0,4cm}

The curl is semi-regular isotopic to the full-twist $\sigma_1^2$. Consequently, $rot(\sigma_1 \cup 2K)$ followed by $push(\sigma_1^2,\sigma_1 \cup2K)$ is a loop in the moduli space of long 2-cables up to semi-regular isotopy (instead of knots in the solid torus). Our calculations above show that $R_1^{(2)}$ vanishes on this loop in the case that $K$ is the right trefoil. Hence in order to get non-trivial 1-cocycles it was very important to replace long knots by knots in the solid torus (where $rot$ and $push$ are {\em each} already a loop) and also to {\em break the symmetry}: the global type $r(a,n,a)$ of R III moves contributes, but the global type $r(a,a,n)$ does not contribute to the same 1-cocycle. These two global types are interchanged by the symmetry which is generated by orientation reversing of the knot together with the hyper elliptic involution of the solid torus (hence the knot represents still the homology class $n$). In particular, it could be a priory possible that our 1-cocycles detect the non-invertibility of a knot because they use only exactly one of these two global types of R III moves. (Again, a computer program would be very helpful.)

\section{Proof}

\subsection{Generalities and preparations}

We use the technology which was developed in \cite{F2}. For the convenience of the reader we recall here the main lines of our approach.

We study the discriminant $\Sigma$ of non-generic diagrams in $M$ or $M_n^{reg}$, together with its natural stratification. Our strategy is the following: for an oriented  generic loop in $M$ or $M_n^{reg}$ we associate an integer to the intersection with each stratum in $\Sigma^{(1)}_{tri}$, i.e. to each Reidemeister III move, and we sum up over all moves in the loop. 

\begin{figure}
\centering
\includegraphics{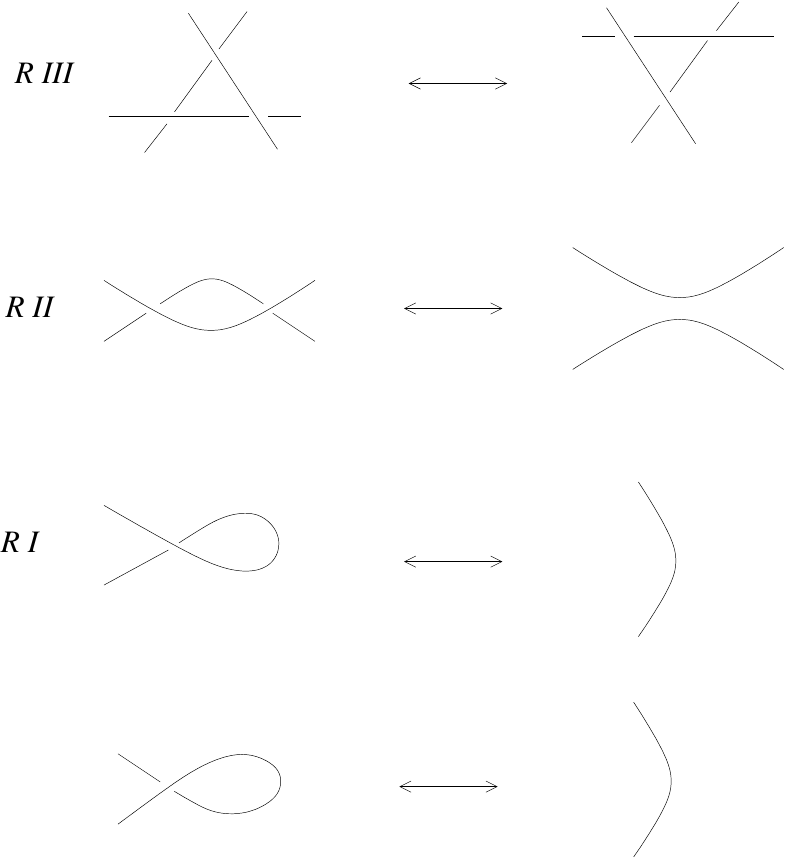}
\caption{\label{moveR} The Reidemeister moves for unoriented knots}  
\end{figure}

The local types of Reidemeister moves for unoriented knots are shown in Fig.~\ref{moveR}.

\begin{figure}
\centering
\includegraphics{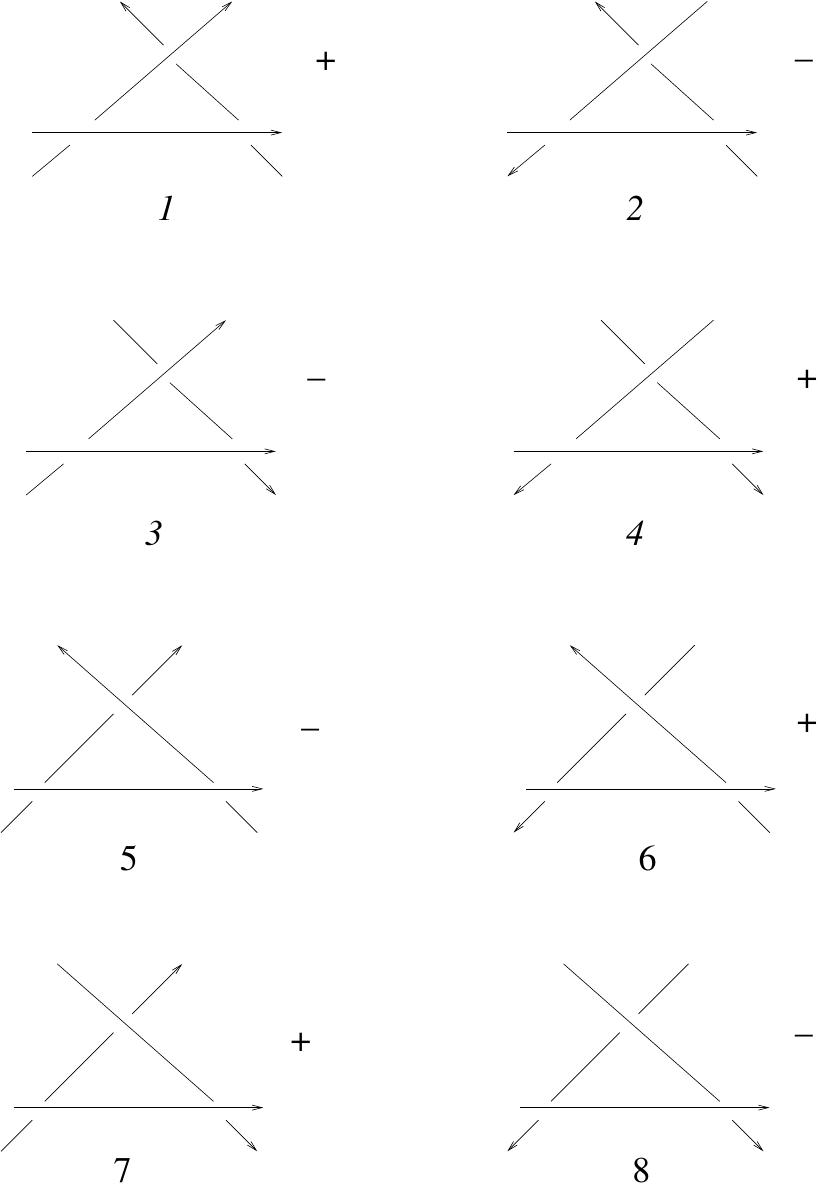}
\caption{\label{loctricross} Local types of a triple crossing}  
\end{figure}

For oriented knots there are exactly eight local types of R III moves, see Fig.~\ref{loctricross}. The sign corresponds to the side of the complement of the discriminant. It coincides with our coorientation for the global type $r$ and it is the opposite for the global type $l$.

\vspace{0,4cm}
In order to show that our 1-cochains are 1-cocycles, we have to prove that the sum is 0 for each meridian of strata of codimension 2, i.e. in $\Sigma^{(2)}$. This is very complex but we have used strata from $\Sigma^{(3)}$ in order to  reduce the proof to a few strata in $\Sigma^{(2)}$. It follows that our sum is invariant under generic homotopies of loops in $M$ or $M_n^{reg}$. But it takes its values in an abelian ring and hence it is a 1-cocycle. Showing that the 1-cocycle is $0$ on the meridians of the quadruple crossings $\Sigma^{(2)}_{quad}$ is by far the hardest part. This corresponds to finding a new solution of the {\em tetrahedron equation}.

Consider four oriented straight lines which form a braid and such that the intersection of their projection into $\mathbb{C}$ consists of a single point. We call this an {\em ordinary quadruple crossing}. After a generic perturbation of the four lines we will see now exactly six ordinary crossings. We assume that all six crossings are positive and we call the corresponding quadruple crossing a {\em positive quadruple crossing}. Quadruple crossings form smooth strata of codimension 2 in  the topological moduli space of lines in 3-space which is equipped with a fixed projection $pr$. Each generic point in such a stratum is adjacent to exactly eight smooth strata of codimension 1. Each of them corresponds to configurations of lines which have exactly one ordinary triple crossing besides the remaining ordinary crossings. We number the lines from 1 to 4 from the lowest to the highest (with respect to the projection $pr$).

\begin{figure}
\centering
\includegraphics{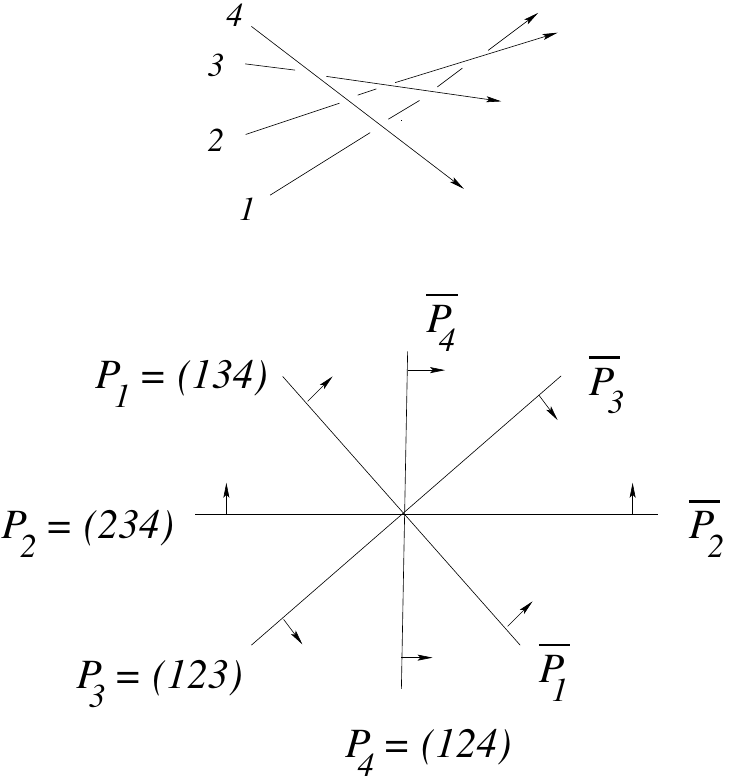}
\caption{\label{quadcros} The intersection of a normal 2-disc of a positive quadruple crossing with the strata of triple crossings}  
\end{figure}

\begin{figure}
\centering
\includegraphics{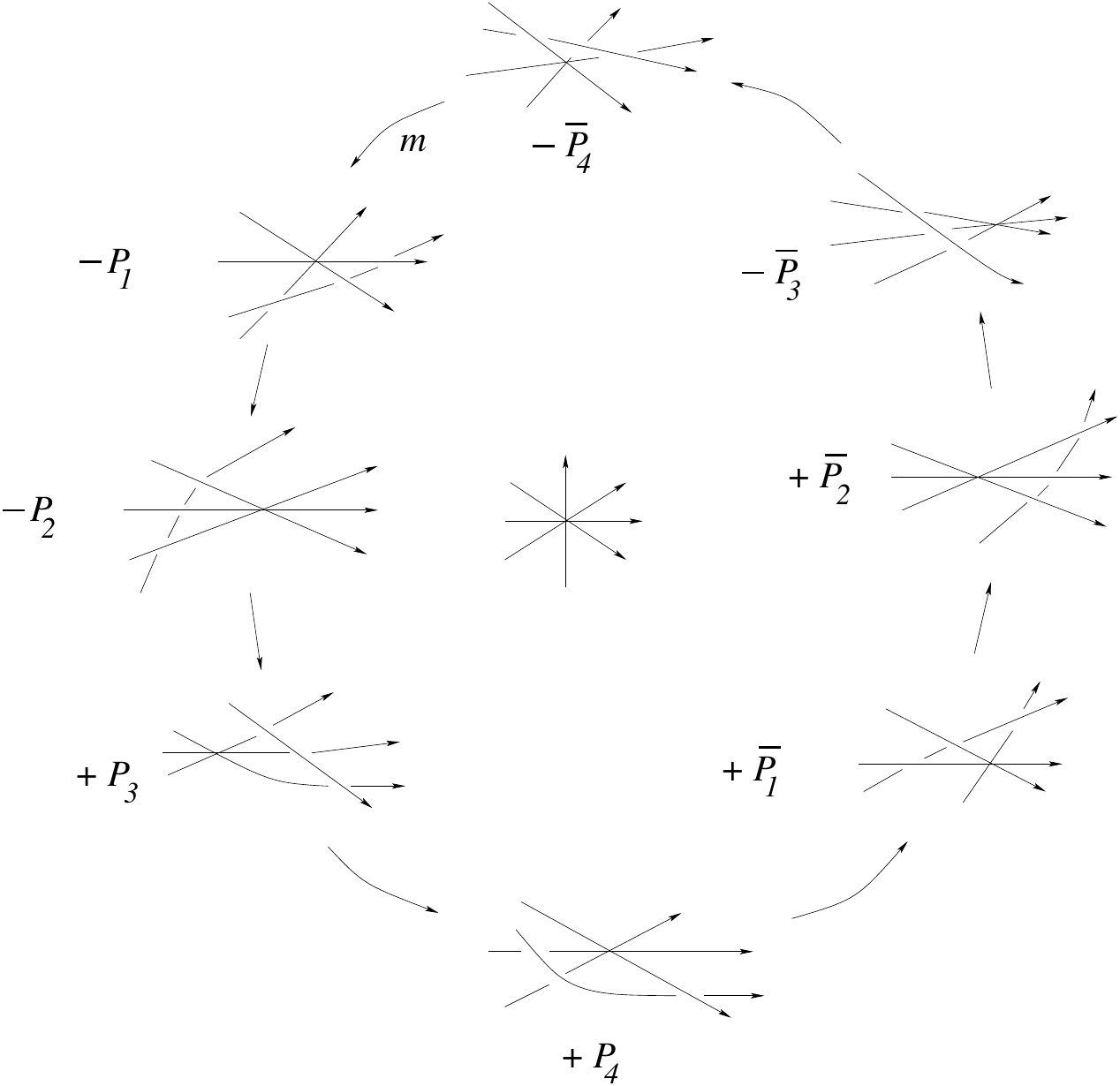}
\caption{\label{unfoldquad} Unfolding of a positive quadruple crossing}  
\end{figure}

The eight strata of triple crossings glue pairwise together to form four smooth strata which intersect pairwise transversely in the stratum of the quadruple crossing, see Fig.~\ref{unfoldquad}, compare \cite{F2} and \cite{FK}. The strata of triple crossings are determined by the names of the three lines which give the triple crossing. For shorter writing we give them names from $P_1$ to $P_4$ and $\bar P_1$ to $\bar P_4$ for the corresponding stratum on the other side of the quadruple crossing.  We show the intersection of a normal 2-disc of the stratum of codimension 2 of a positive quadruple crossing with the strata of codimension 1 in Fig.~\ref{quadcros}. The strata of codimension 1 have a natural coorientation, compare Definition 2. We could interpret the six ordinary crossings as the edges of a tetrahedron and the four triple crossings likewise as the vertices or the 2-faces of the tetrahedron. For the classical tetrahedron equation one associates to each stratum $P_i$, i.e. to each vertex or equivalently to each 2-face of the tetrahedron, some operator (or some R-matrix)\index{R-matrix} which depends  only on the names of the three lines and to each stratum $\bar P_i$ the inverse operator. The tetrahedron equation says now that if we go along the meridian then the composition of these operators is equal to the identity. Notice, that in the literature, see e.g. \cite{KaKo}, one considers planar configurations of lines. But this is of course equivalent to our situation because all crossings are positive and hence the lift of the lines into 3-space is determined by the planar picture. Moreover, each move of the lines in the plane which preserves the transversality lifts to an isotopy of the lines in 3-space. The tetrahedron equation has many solutions, the first one was found by Zamolodchikov, see e.g. \cite{KaKo}.

\index{normal 2-disc of a positive quadruple crossing}

\begin{figure}
\centering
\includegraphics{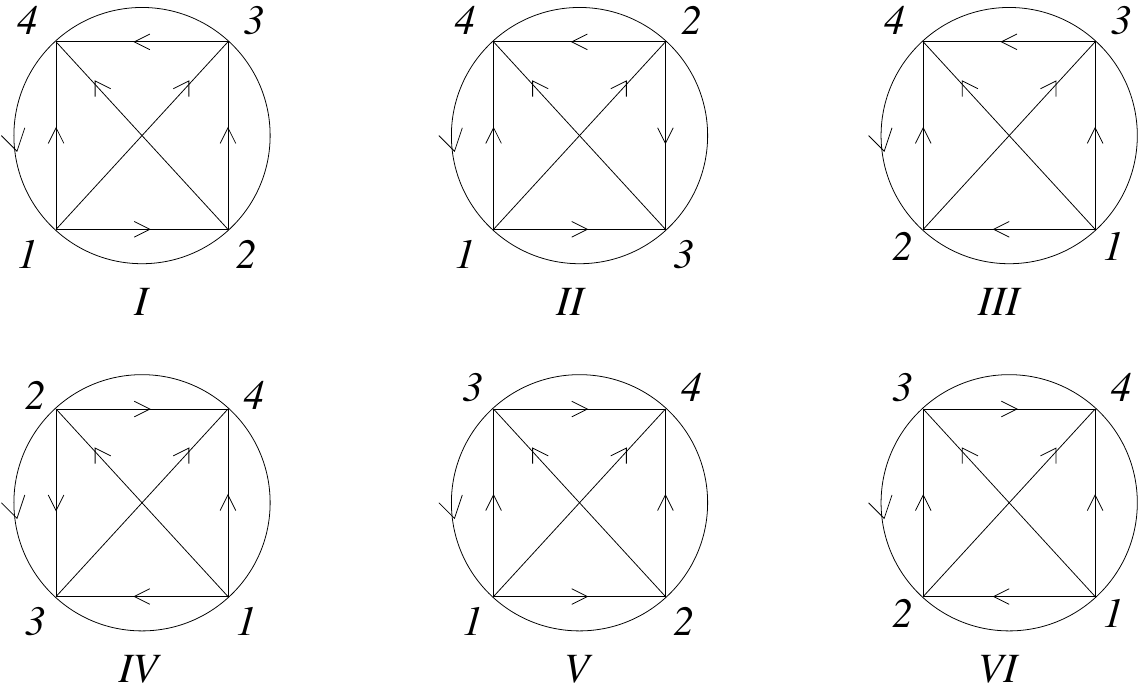}
\caption{\label{globquad} The global types of quadruple crossings}  
\end{figure}

However, the solutions of the classical tetrahedron equation are not well adapted in order to construct 1-cocycles for moduli spaces of knots. A local solution of the tetrahedron equation is of no use for us, because as already pointed out there are no integer valued 1-cocycles for all knots in the 3-sphere. We have to replace them by long knots or more generally by points in $M_n$. For knots in the solid torus $V^3$ we can now associate to each crossing in the diagram a winding number (i.e. a homology class in $H_1(V^3)$) in a canonical way. Therefore we have to consider six different positive tetrahedron equations, corresponding to the six different abstract closures of the four lines to a circle and in each of the six cases we have to consider all possible winding numbers of the six crossings. We call this the {\em positive global tetrahedron equations}. 
There are exactly six global types of positive quadruple crossings without the homological markings. We show them in Fig.~\ref{globquad}

One easily sees that there are exactly {\em forty eight} local types of quadruple crossings (analog to the eight local types of triple crossings).
\vspace{0,4cm}

\begin{figure}
\centering
\includegraphics{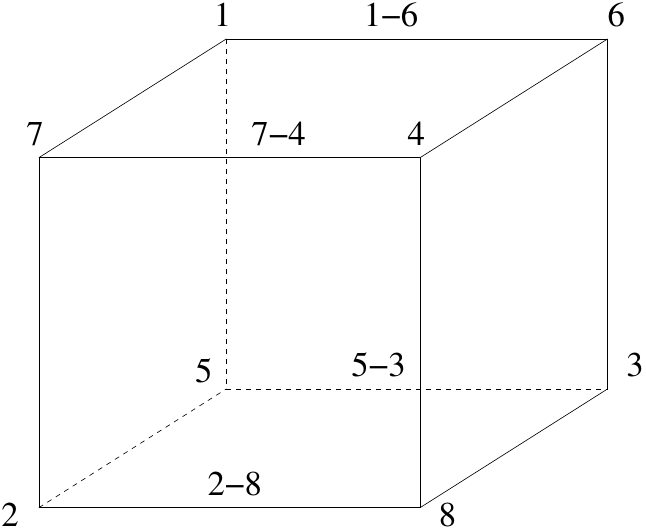}
\caption{\label{gamma} The graph $\Gamma$}  
\end{figure}

We study the relations of the local types of triple crossings in what we call the {\em cube equations}.
Triple crossings come together in points of $\Sigma^{(2)}_{trans-self}$, i.e. an auto-tangency with in addition a transverse branch \cite{F2}. But one easily sees that the global type of the triple crossings (i.e. its Gauss diagram with the homological markings but  without the writhe) is always preserved. We make now a graph $\Gamma$ for each global type of a triple crossing in the following way: the vertices correspond to the different local types of triple crossings. We connect two vertices by an edge if and only if the corresponding strata of triple crossings are adjacent to a stratum of $\Sigma^{(2)}_{trans-self}$. We have shown that the resulting graph is the 1-skeleton of the 3-dimensional cube $I^3$, see Fig.~\ref{gamma}. 
In particular, it is connected. The edges of the graph $\Gamma = skl_1(I^3)$ correspond to the types of strata in $\Sigma^{(2)}_{trans-self}$. The solution of the positive tetrahedron equation tells us what is the contribution to the 1-cocycle of a positive triple crossing (i.e. all three involved crossings are positive). The meridians of the strata from $\Sigma^{(2)}_{trans-self}$ give equations which allow us to determine the contributions of all other types of triple crossings. However, a global phenomenon occurs: each loop in $\Gamma$ could give an additional equation. Evidently, it suffices to consider the loops which are the boundaries of the 2-faces from $skl_2(I^3)$. We call all the equations which come from the meridians of $\Sigma^{(2)}_{trans-self}$ and from the loops in $\Gamma = skl_1(I^3)$ the {\em cube equations}. (Notice that a loop in $\Gamma$ is more general than a loop in $M$ or $M_n^{reg}$. For a loop in $\Gamma$ we come back to the same local type of a triple crossing but not necessarily to the same whole diagram of the knot.)
\vspace{0,4cm}

Our strategy is the following: after solving the commutation relations we find a solution of the positive global tetrahedron equation. We solve then the cube equations by adding {\em correction terms} for different local types of triple crossings, but which vanish for positive triple crossings. We have shown in \cite{F2} that the resulting 1-cochain is then a 1-cocycle in $M^{reg}$ or $M_n^{reg}$.
\vspace{0,4cm}

\begin{figure}
\centering
\includegraphics{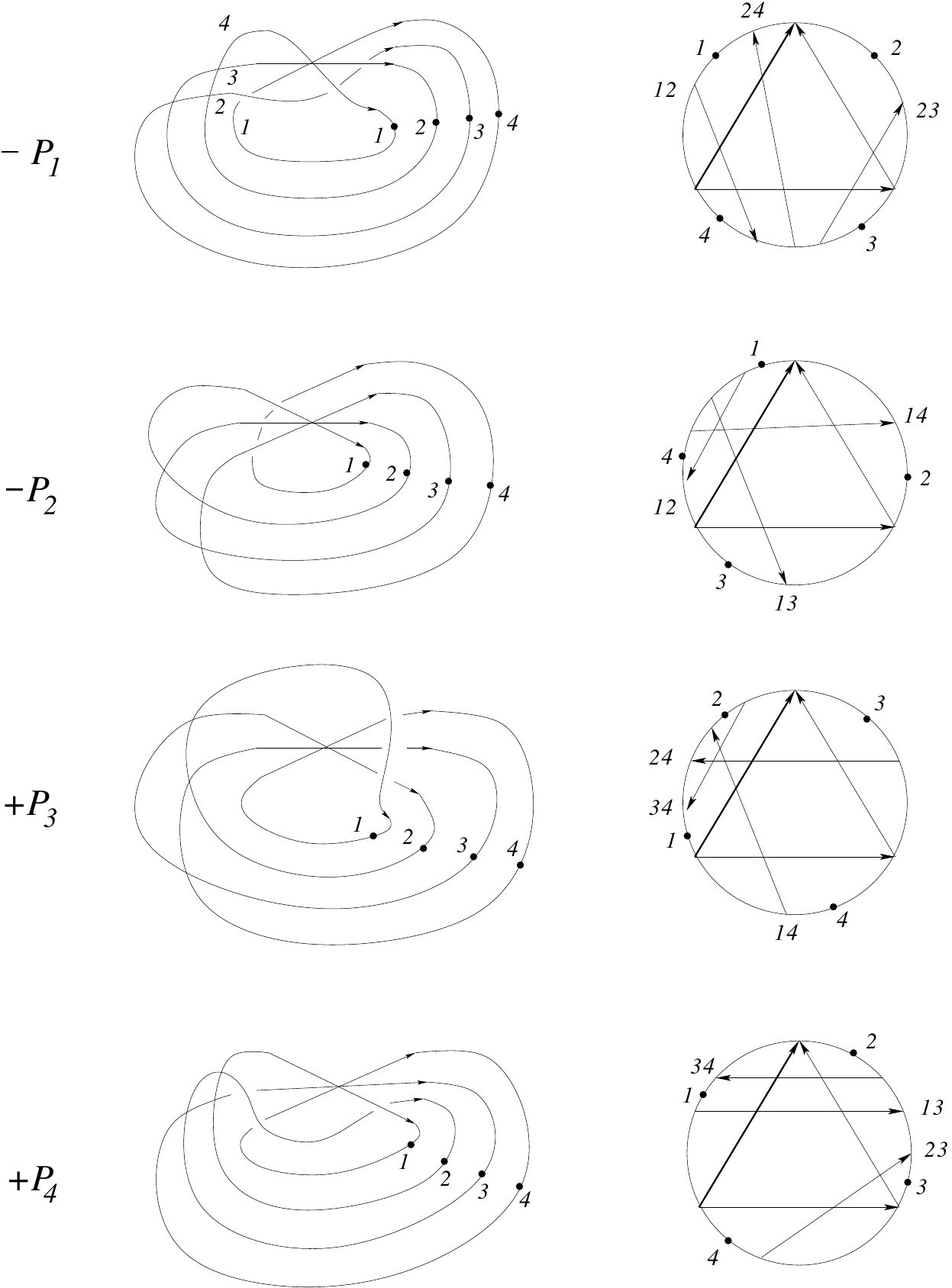}
\caption{\label{Iglob} The first half of the meridian for global type I}  
\end{figure}

\begin{figure}
\centering
\includegraphics{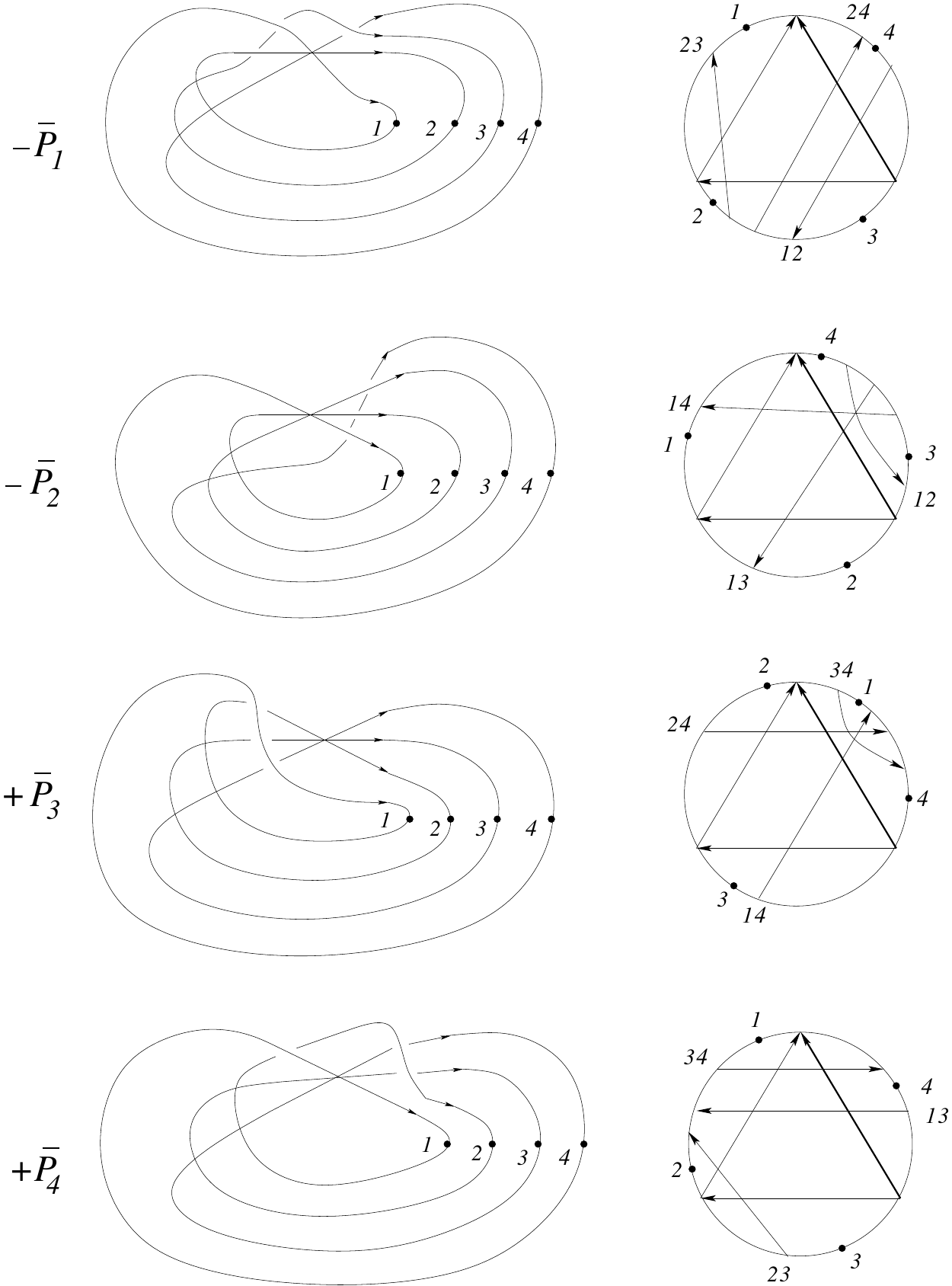}
\caption{\label{VI2glob} The second half of the meridian for global type VI}  
\end{figure}

For the convenience of the reader we give here again the figures which are our main tool \cite{F2}. Let consider the global positive quadruple crossings. We naturally identify crossings in an isotopy outside Reidemeister moves of type I and II. The Gauss diagrams of the unfoldings of the quadruple crossings are given in Fig.~\ref{Iglob} up to Fig.~\ref{VI2glob}. For the convenience of the reader (and for further research) we indicate also the different possibilities for the point at infinity in the case of long knots. The (positive) crossing between the local branch $i$ and the local branch $j$ is always denoted by $ij$. 

\begin{figure}
\centering
\includegraphics{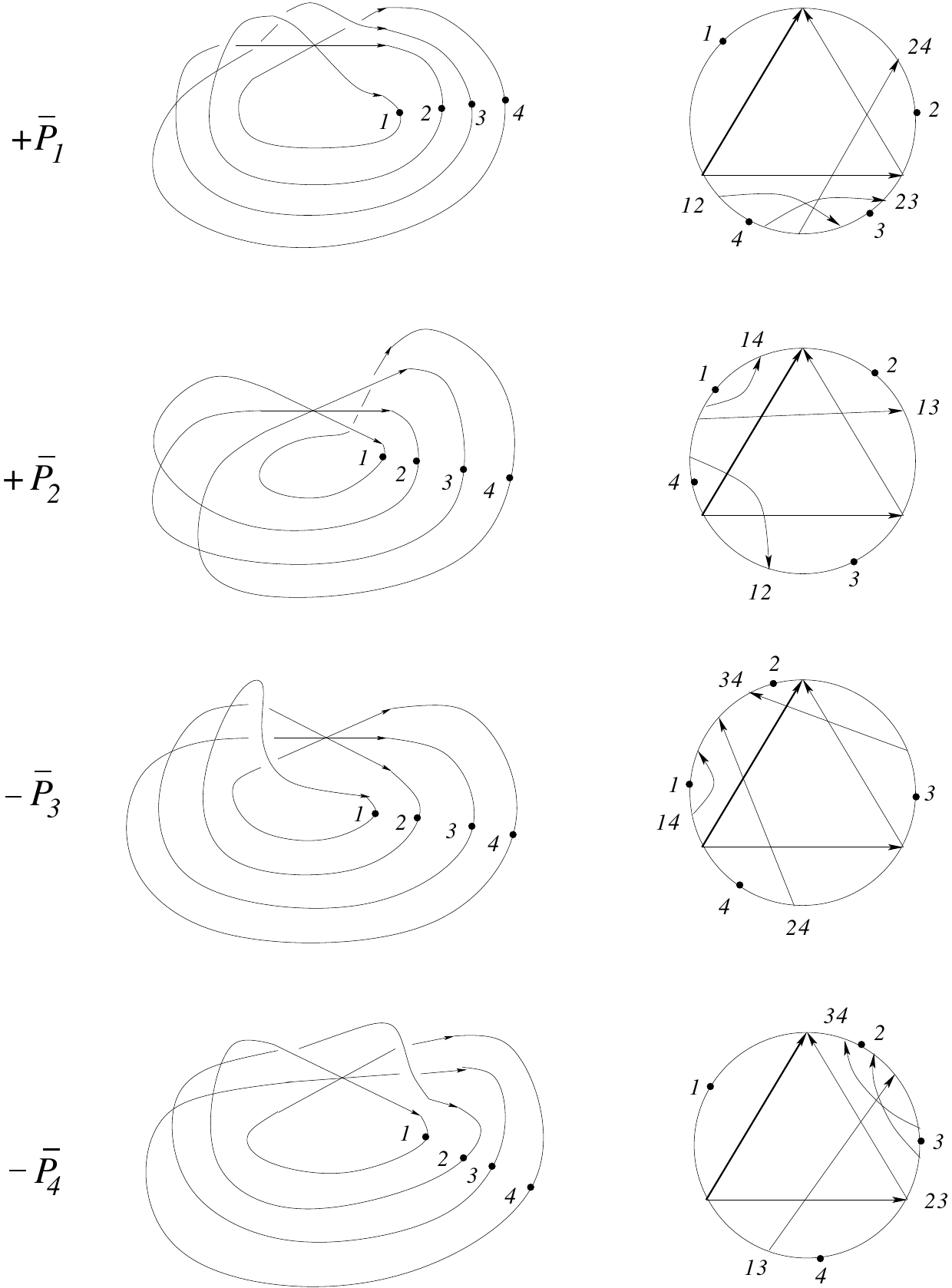}
\caption{\label{I2glob} The second half of the meridian for global type I}  
\end{figure}

\begin{figure}
\centering
\includegraphics{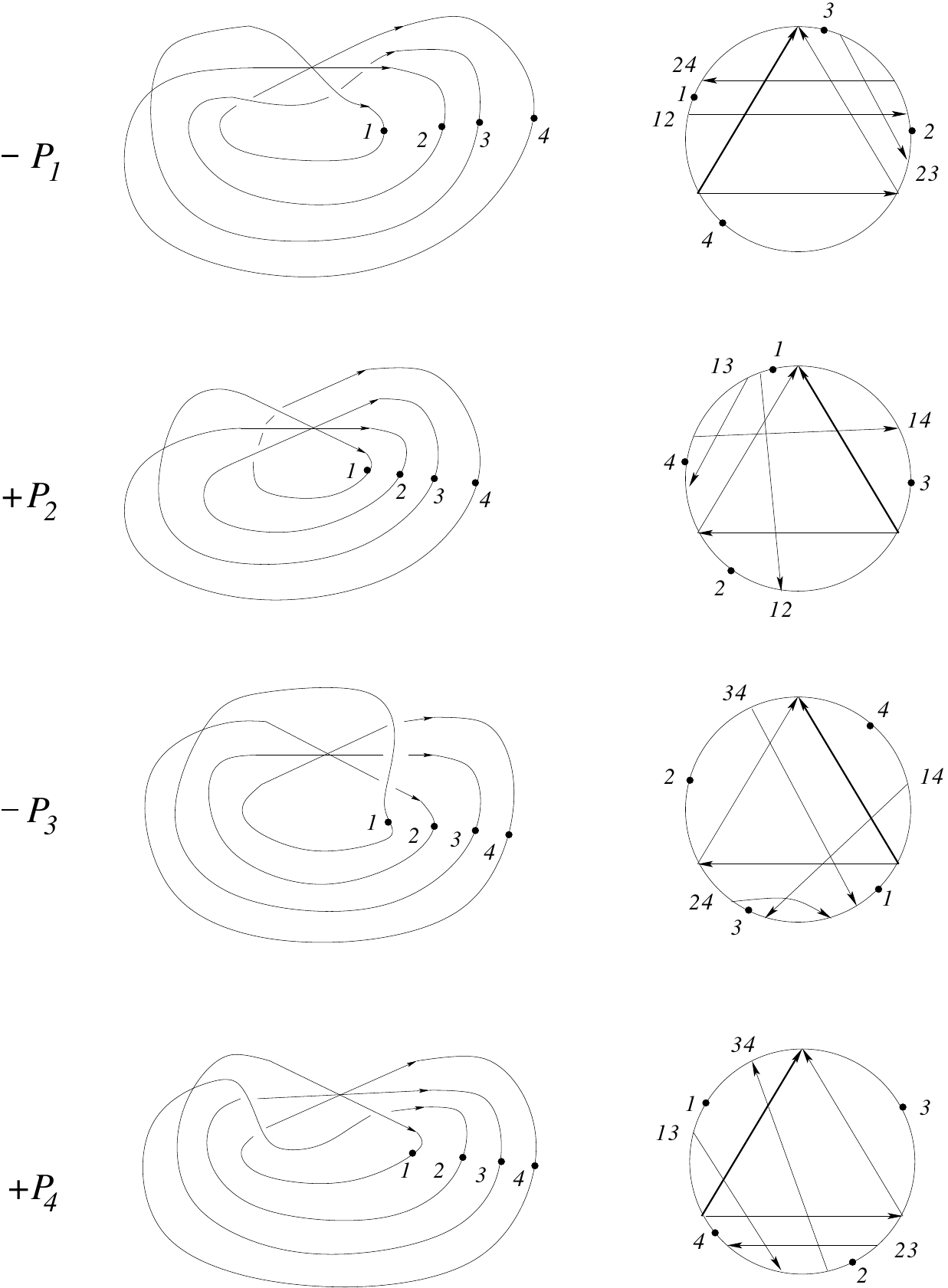}
\caption{\label{IIglob} The first half of the meridian for global type II}  
\end{figure}

\begin{figure}
\centering
\includegraphics{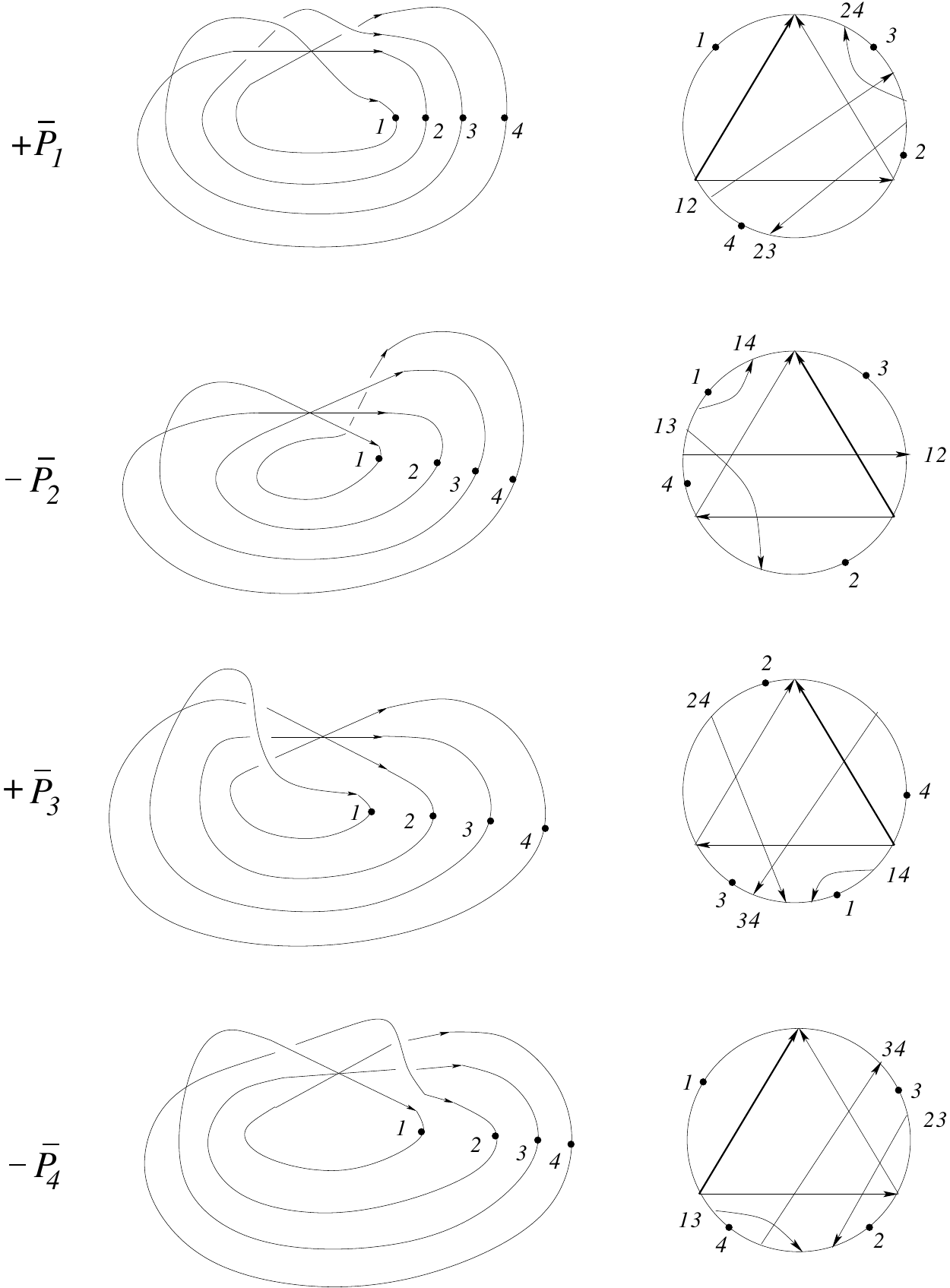}
\caption{\label{II2glob} The second half of the meridian for global type II}  
\end{figure}

\begin{figure}
\centering
\includegraphics{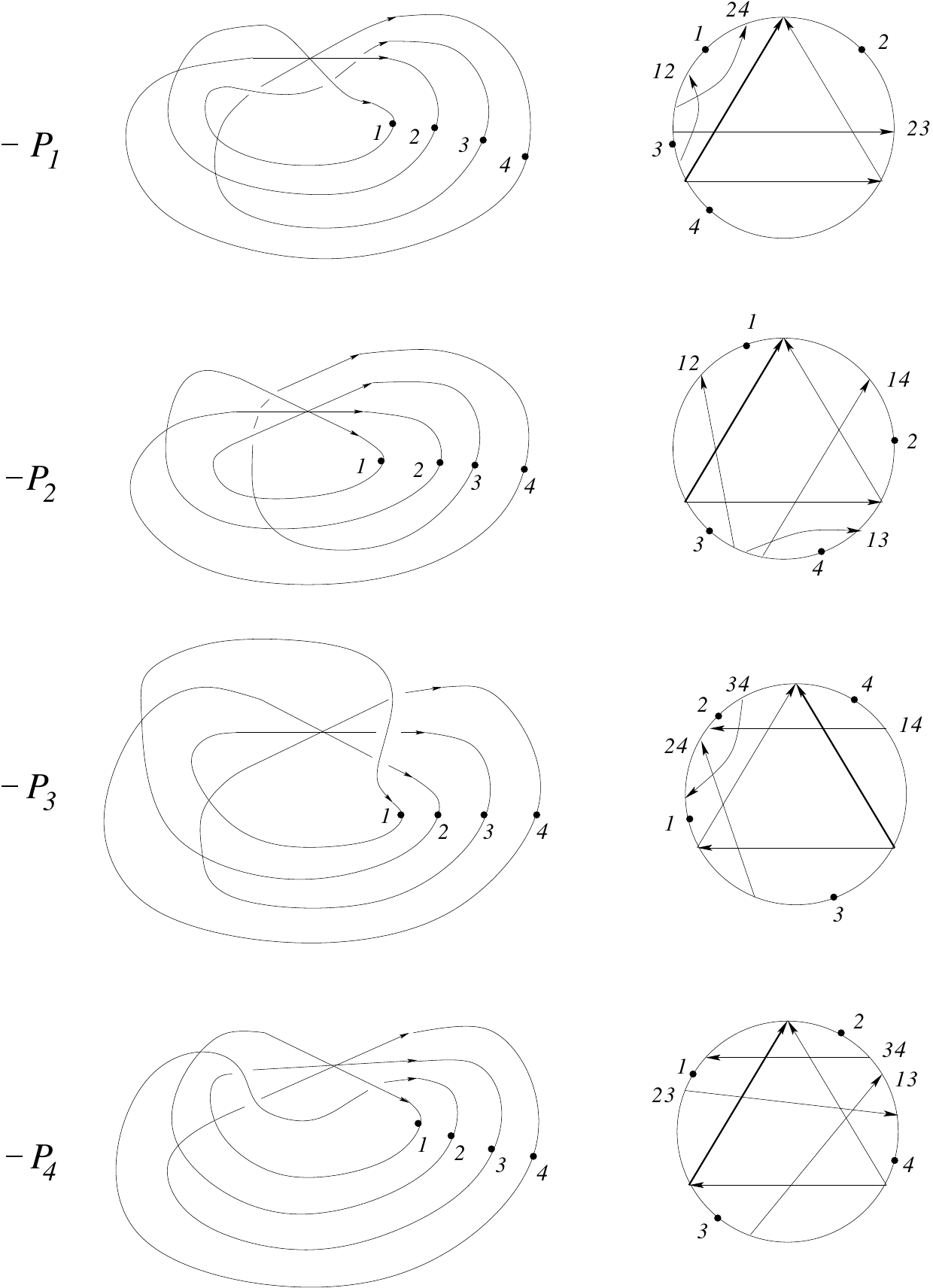}
\caption{\label{IIIglob} The first half of the meridian for global type III}  
\end{figure}

\begin{figure}
\centering
\includegraphics{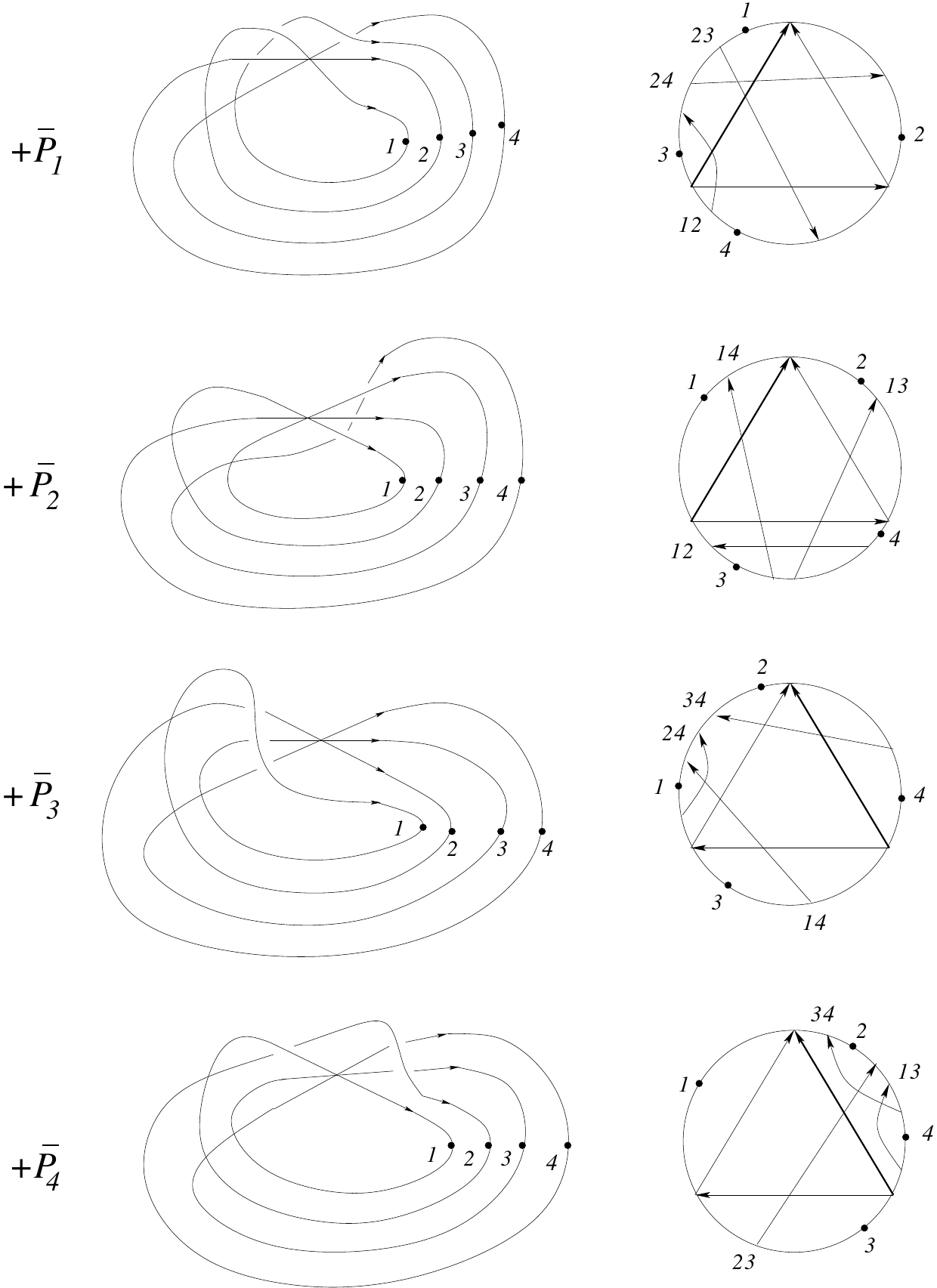}
\caption{\label{III2glob} The second half of the meridian for global type III}  
\end{figure}

\begin{figure}
\centering
\includegraphics{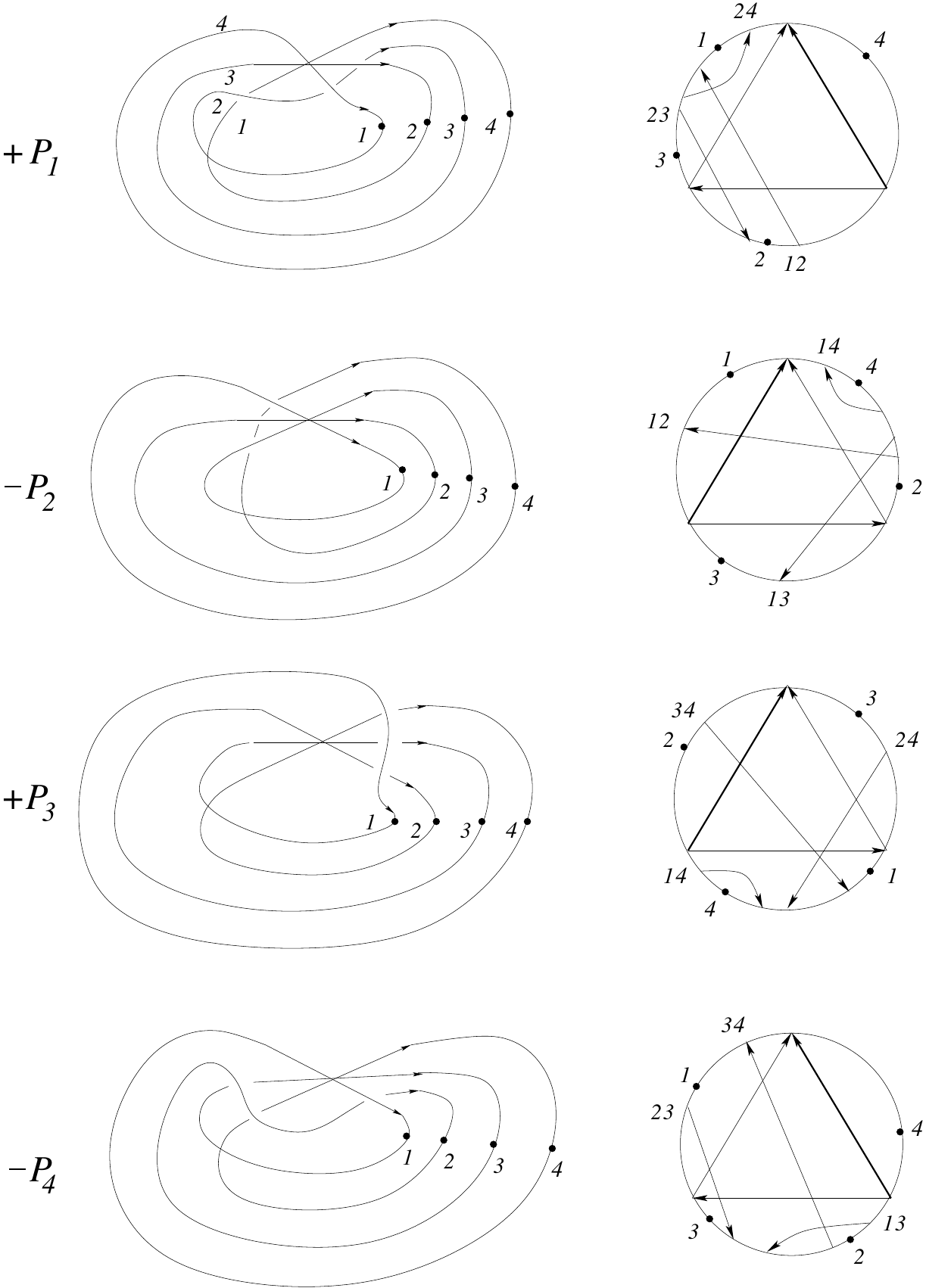}
\caption{\label{IVglob} The first half of the meridian for global type IV}  
\end{figure}

\begin{figure}
\centering
\includegraphics{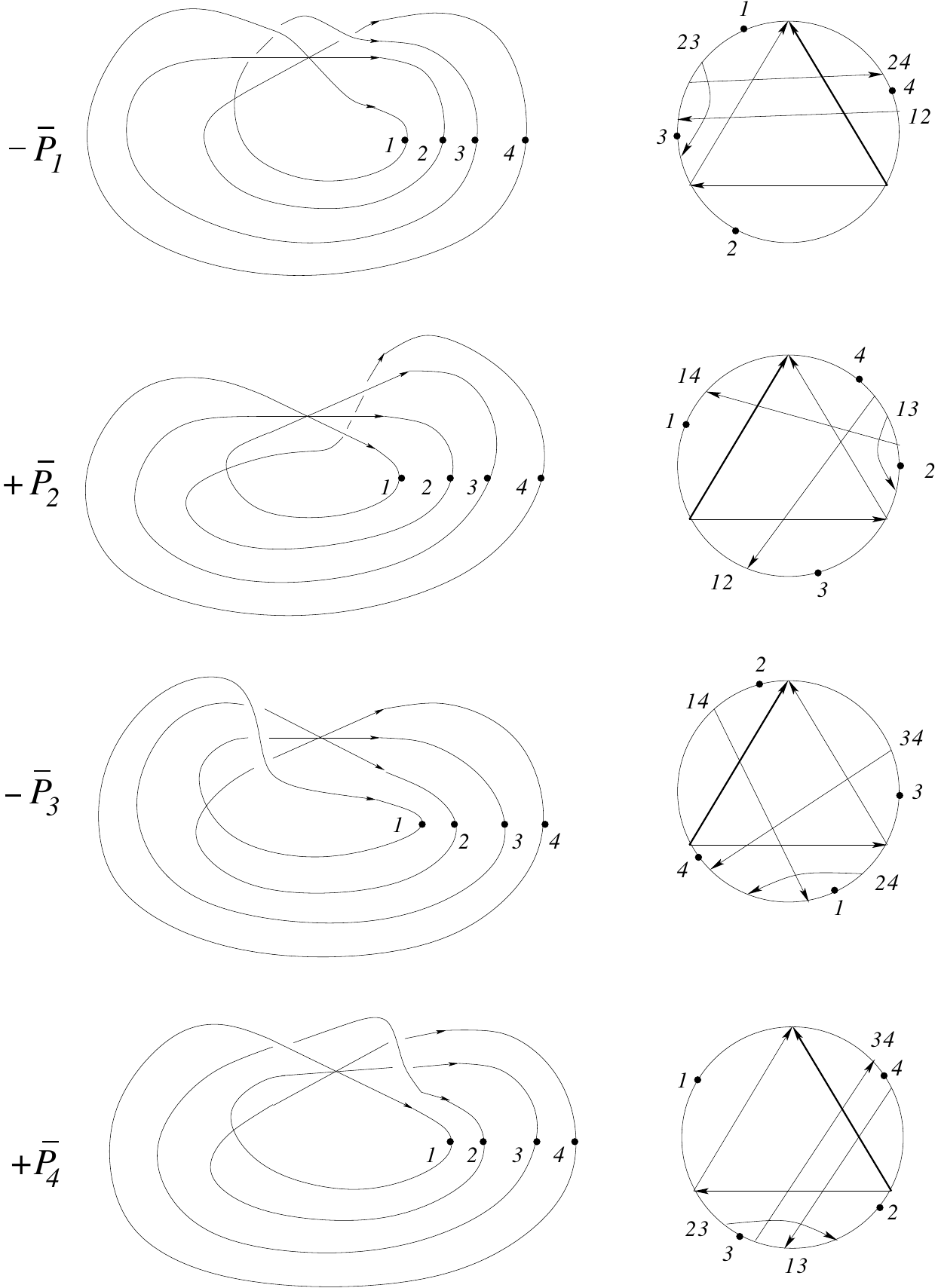}
\caption{\label{IV2glob} The second half of the meridian for global type IV}  
\end{figure}

\begin{figure}
\centering
\includegraphics{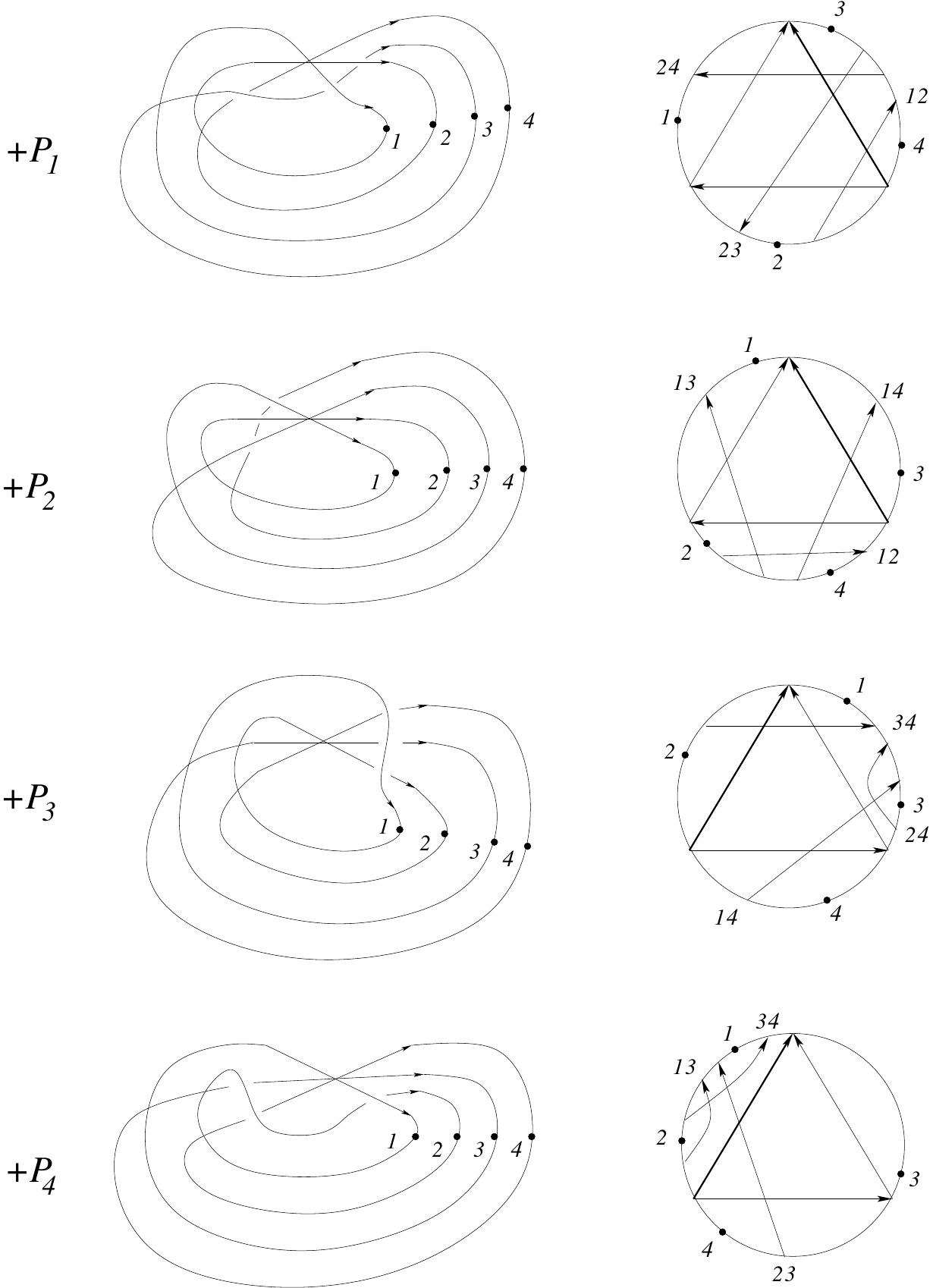}
\caption{\label{Vglob} The first half of the meridian for global type V}  
\end{figure}

\begin{figure}
\centering
\includegraphics{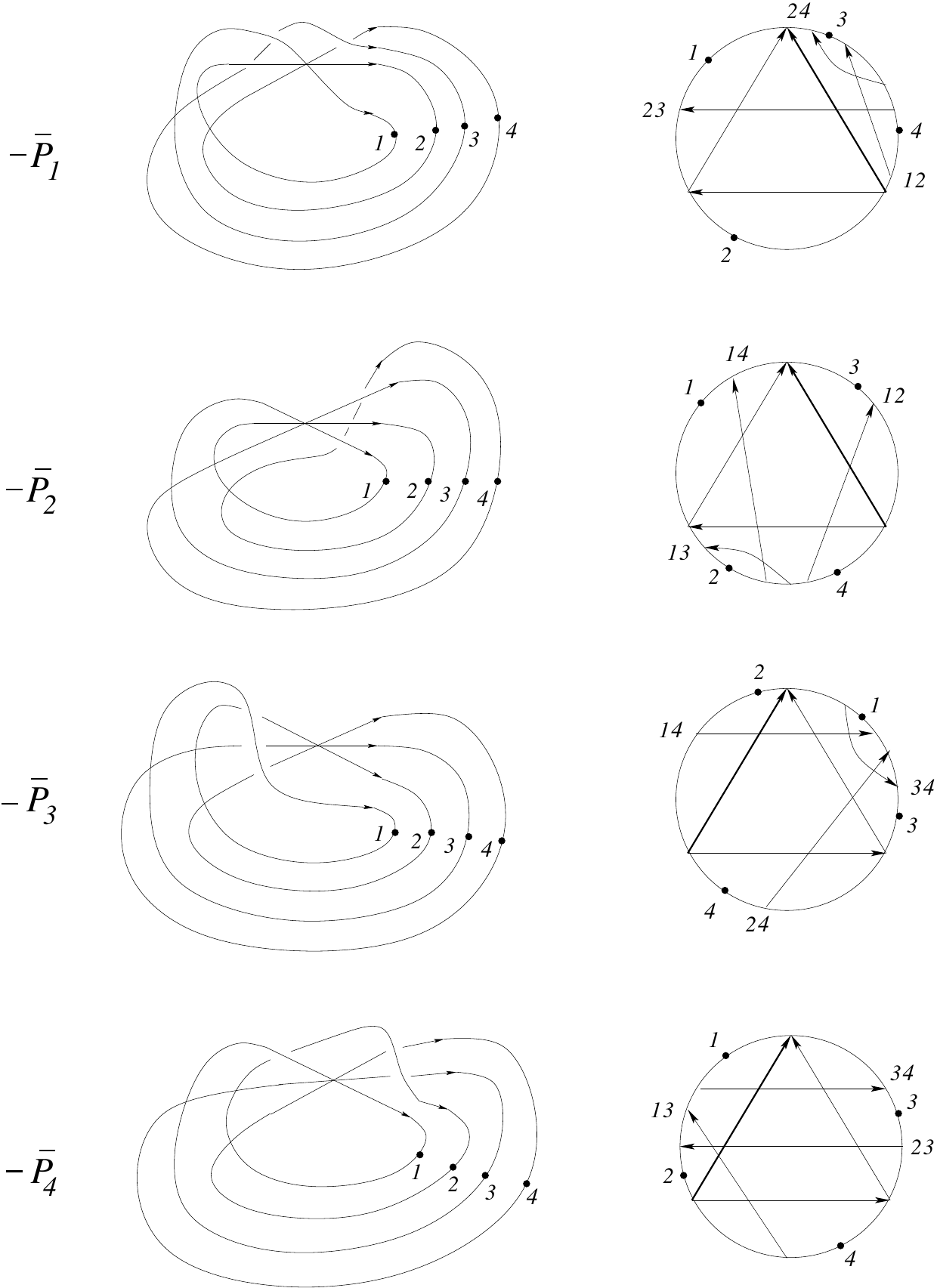}
\caption{\label{V2glob} The second half of the meridian for global type V}  
\end{figure}

\begin{figure}
\centering
\includegraphics{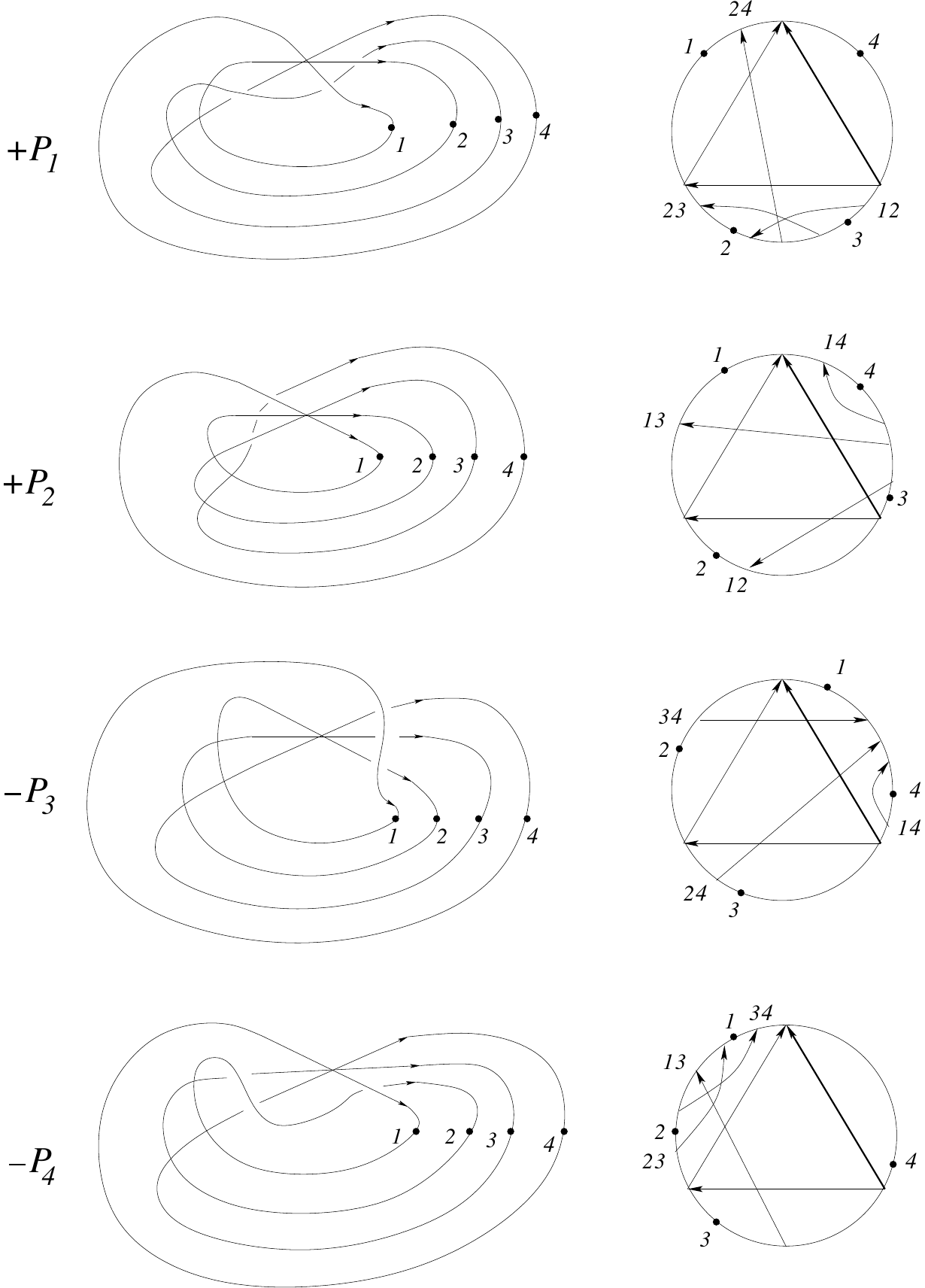}
\caption{\label{VIglob} The first half of the meridian for global type VI}  
\end{figure}

\subsection{$R_a^{(2)}$ satisfies the commutation relations (a)}

As already pointed out in the Introduction, we have only to check that the contribution of a triple crossing $p$ to $R_a^{(2)}$ does not change if we pass through another independent triple crossing $p'$. Evidently, we have only to study triple crossings $p'$ which contain at least one n-crossing and at least one 0-crossing.

\begin{figure}
\centering
\includegraphics{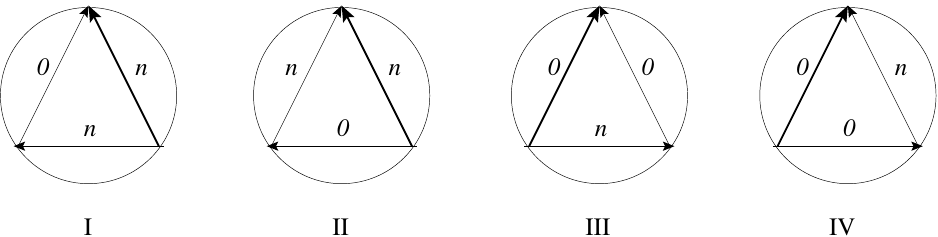}
\caption{\label{globRIIIn} The global types of R III moves with markings $0$ and $n$}  
\end{figure}

In this case there are only four global types of R III moves to study, see Fig.~\ref{globRIIIn}.

\begin{figure}
\centering
\includegraphics{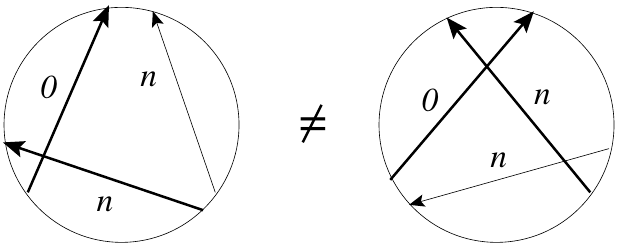}
\caption{\label{indweight} The individual contribution of an $n$-crossing has changed}  
\end{figure}

The individual contribution of each $n$-crossing (i.e. crossing with homological marking $n$) is already invariant for the above global types II, III and IV. However it changes for the global type I, see Fig.~\ref{indweight}.

But we notice that the two $n$-crossings in the case I have their foot in the same arc for $p$ on the circle! Hence, they are $f$-crossings for $p$ only simultaneously and $W_2(p)$ stays invariant.

Notice  that each weight of degree 1 satisfies automatically the commutation relations (a) for regular isotopies (i.e. R I moves are not allowed). In particular, the linking numbers $l(p)$ as well as $W_2(hm)$ are defined by using only single crossings and they stay evidently invariant.
$\Box$
\vspace{0,2 cm}

\subsection{$R_a^{(2)}$ satisfies the positive global tetrahedron equation (b)}

\begin{figure}
\centering
\includegraphics{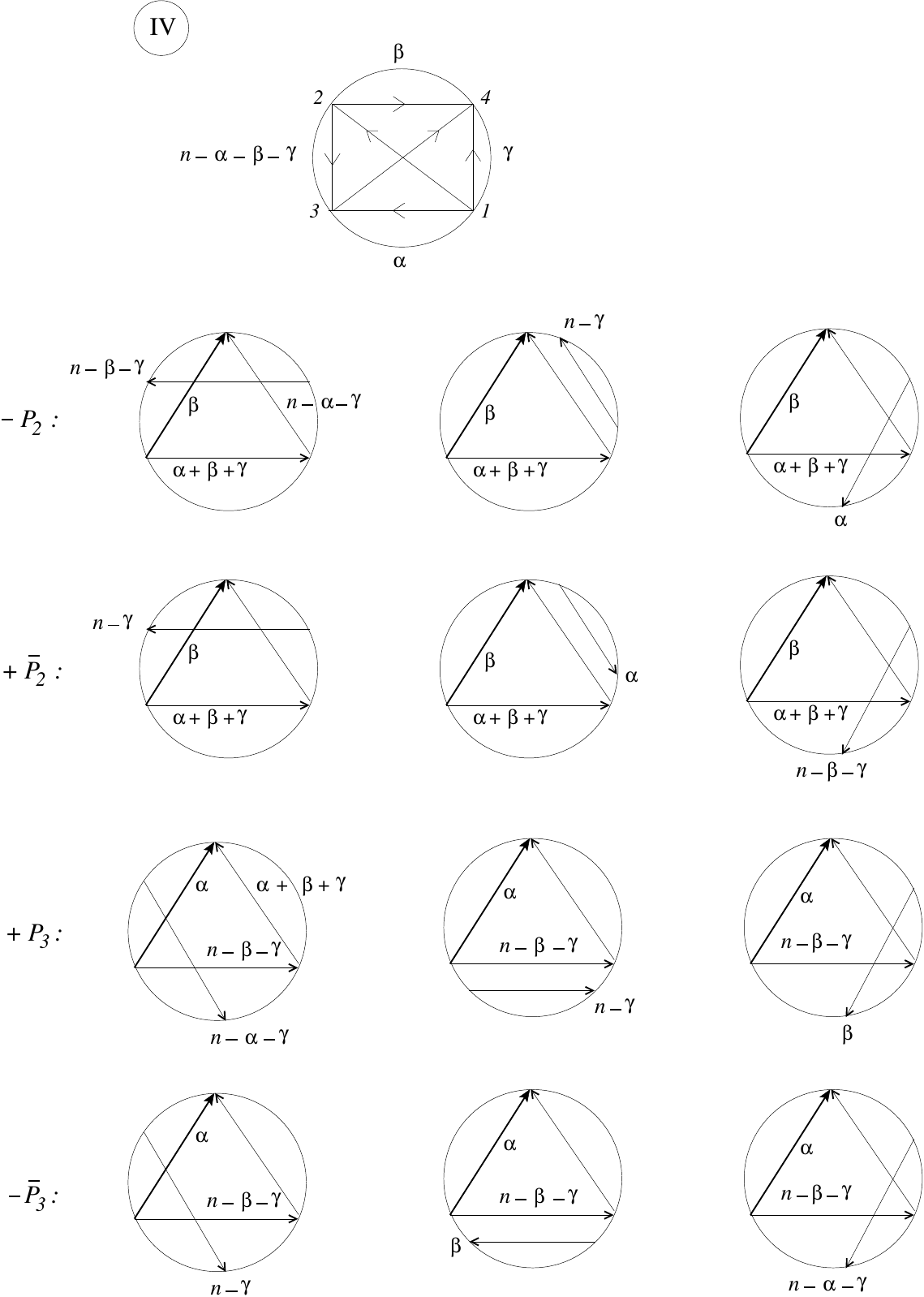}
\caption{\label{ser1} The moves of type $r$ for the global type IV}  
\end{figure}

\begin{figure}
\centering
\includegraphics{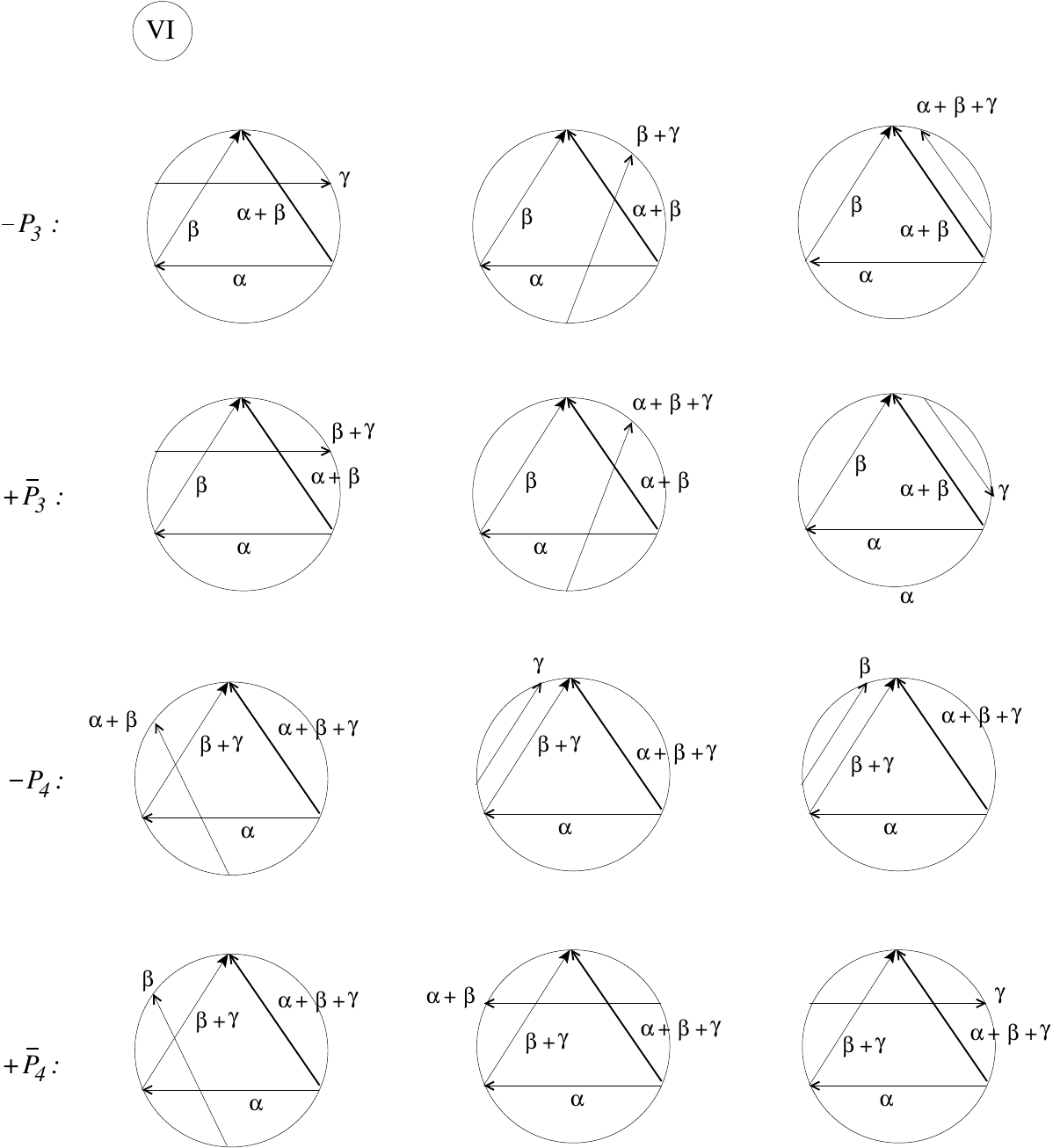}
\caption{\label{ser12}  The remaining moves of type $l$ for the global type VI}  
\end{figure}

We give now the homological markings, which depend on three parameters $\alpha$, $\beta$ and $\gamma$ on the circle.  
Using the previous figures we show in Fig.~\ref{ser1} up to Fig.~\ref{ser12} the homological markings of the arrows. Here $\alpha$, $\beta$ and $\gamma$ are the homology classes represented by the corresponding arcs in the circle (remember that the circle in the plan is always oriented counter-clockwise). 

{\em In $M_n$ we have the unknown parameters  $\alpha$, $\beta$ and $\gamma$ in $\{0, 1, ...,n\}$ and $\alpha+\beta+\gamma \leq n$, because there are no negative loops in the diagrams. }
\vspace{0,2 cm}

\begin{figure}
\centering
\includegraphics{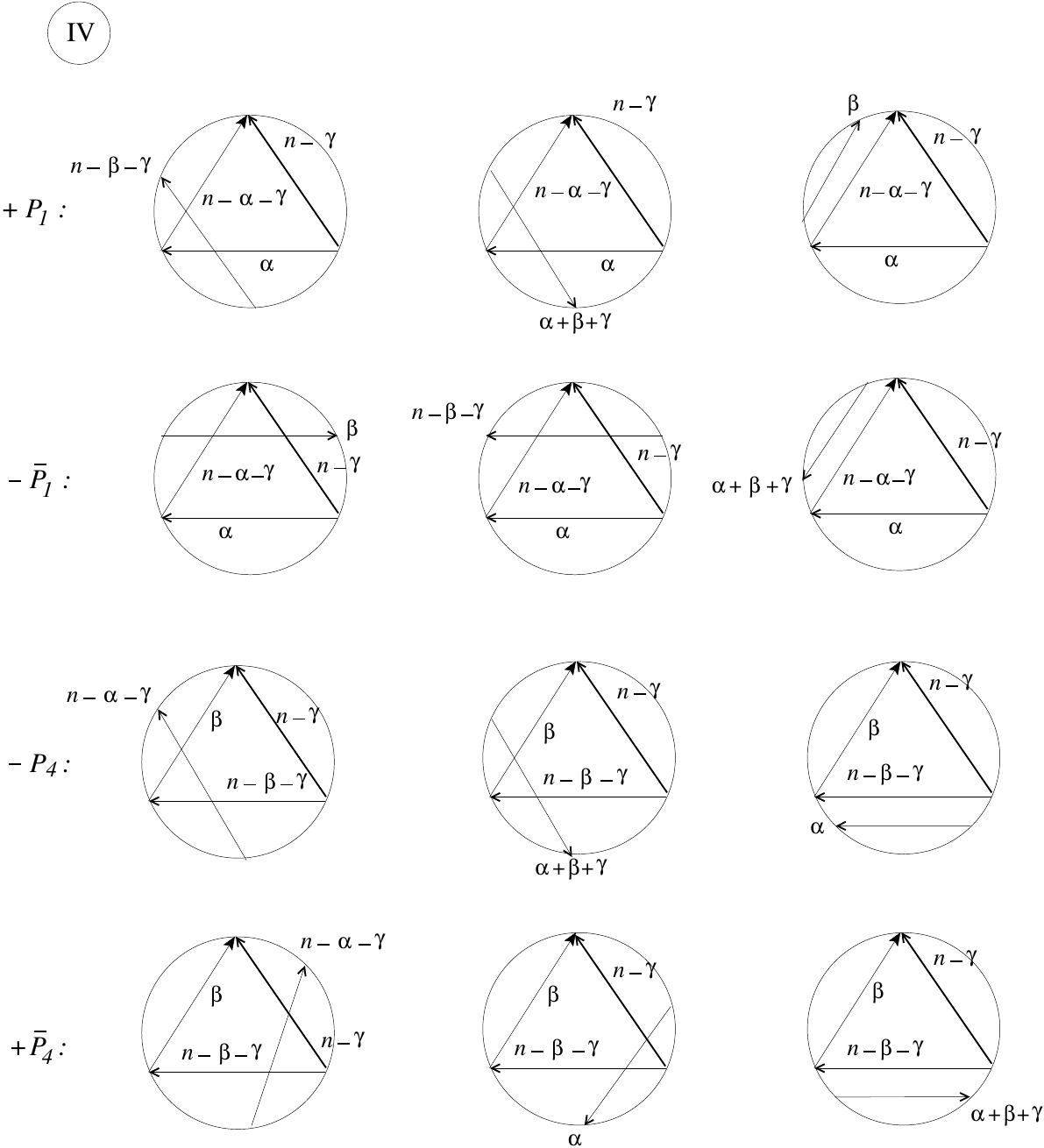}
\caption{\label{ser2} The moves of type $l$ for the global type IV}  
\end{figure}

\begin{figure}
\centering
\includegraphics{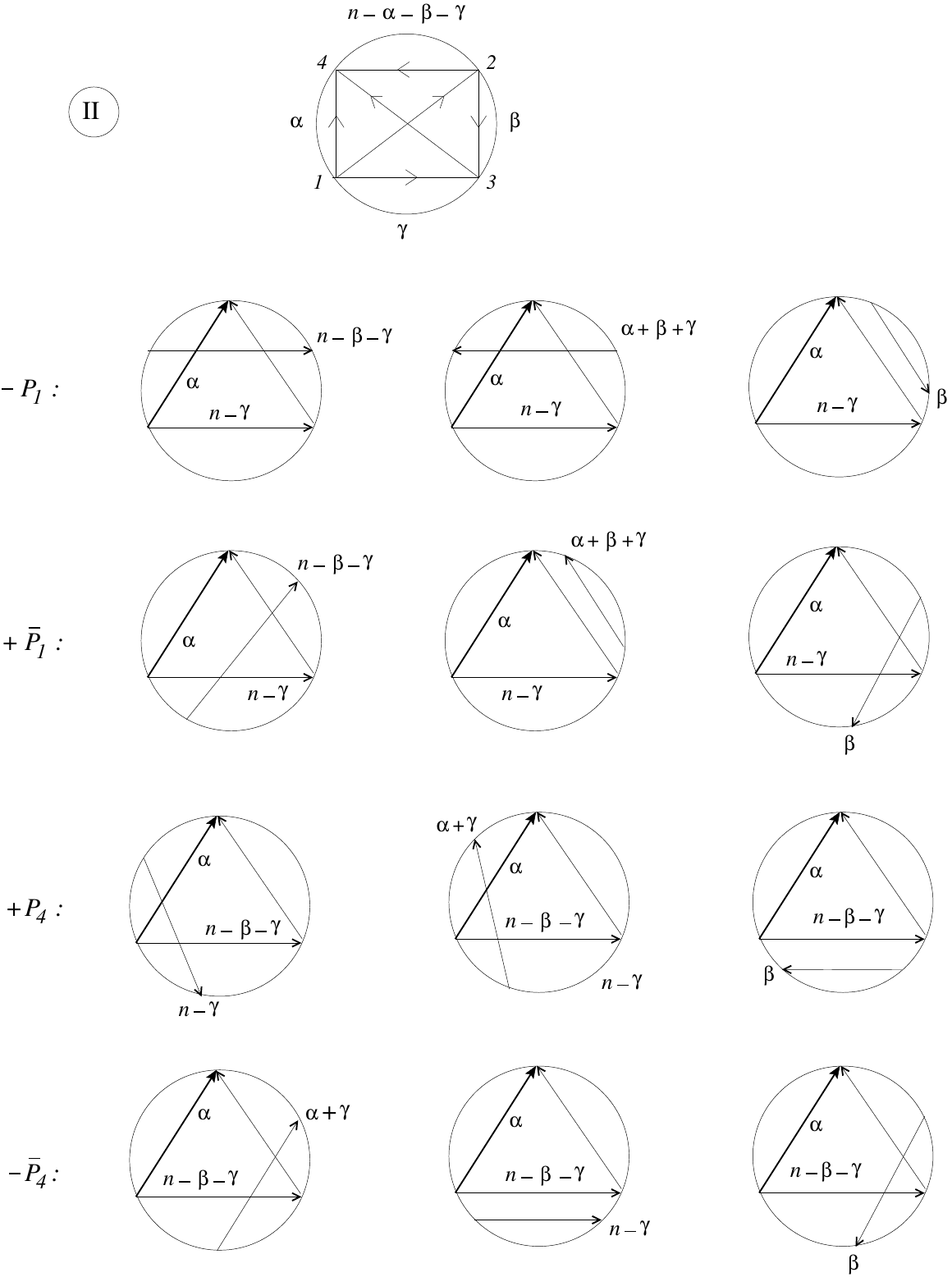}
\caption{\label{ser3} The moves of type $r$ for the global type II}  
\end{figure}

\begin{figure}
\centering
\includegraphics{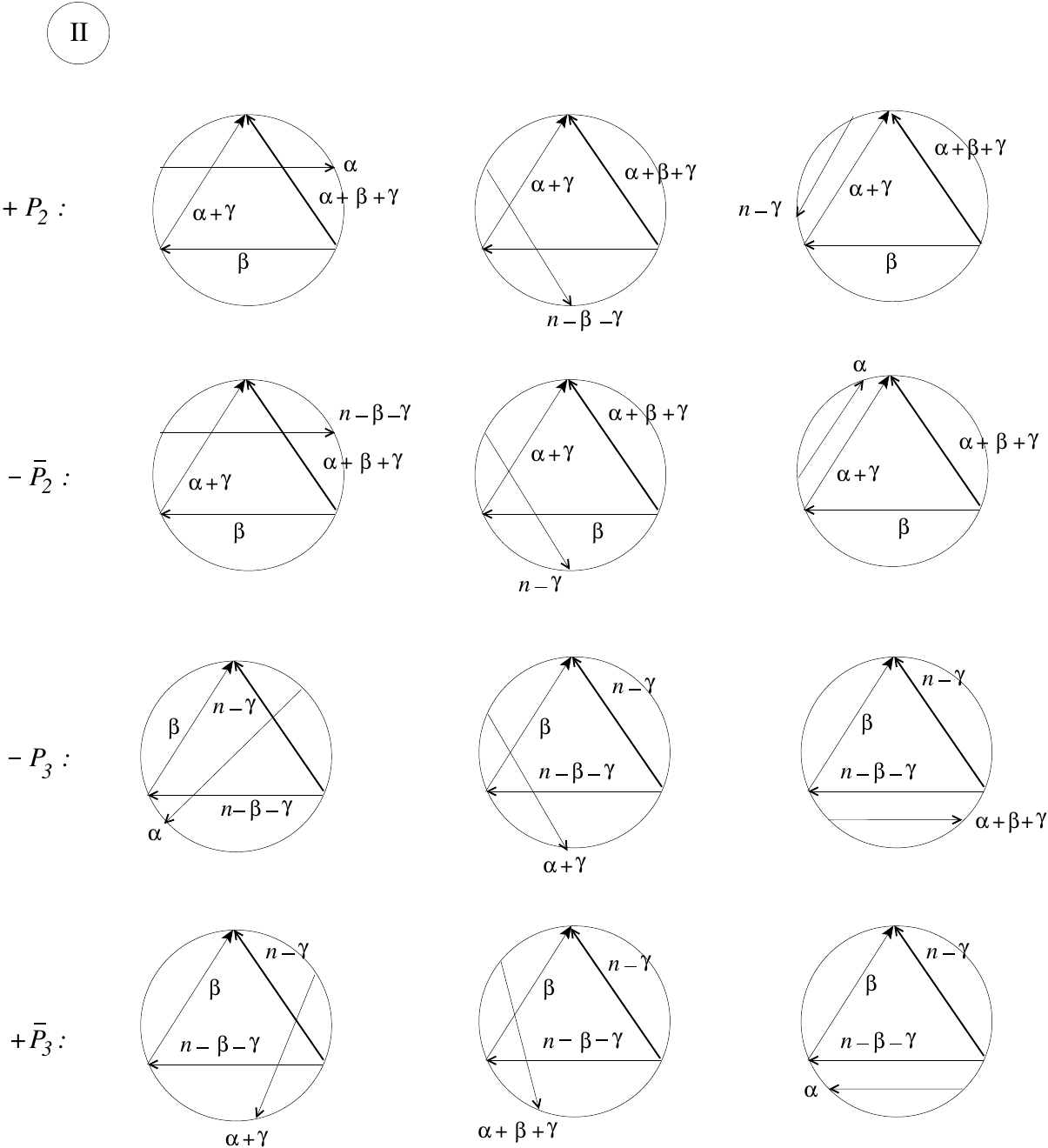}
\caption{\label{ser4} The moves of type $l$ for the global type II}  
\end{figure}

\begin{figure}
\centering
\includegraphics{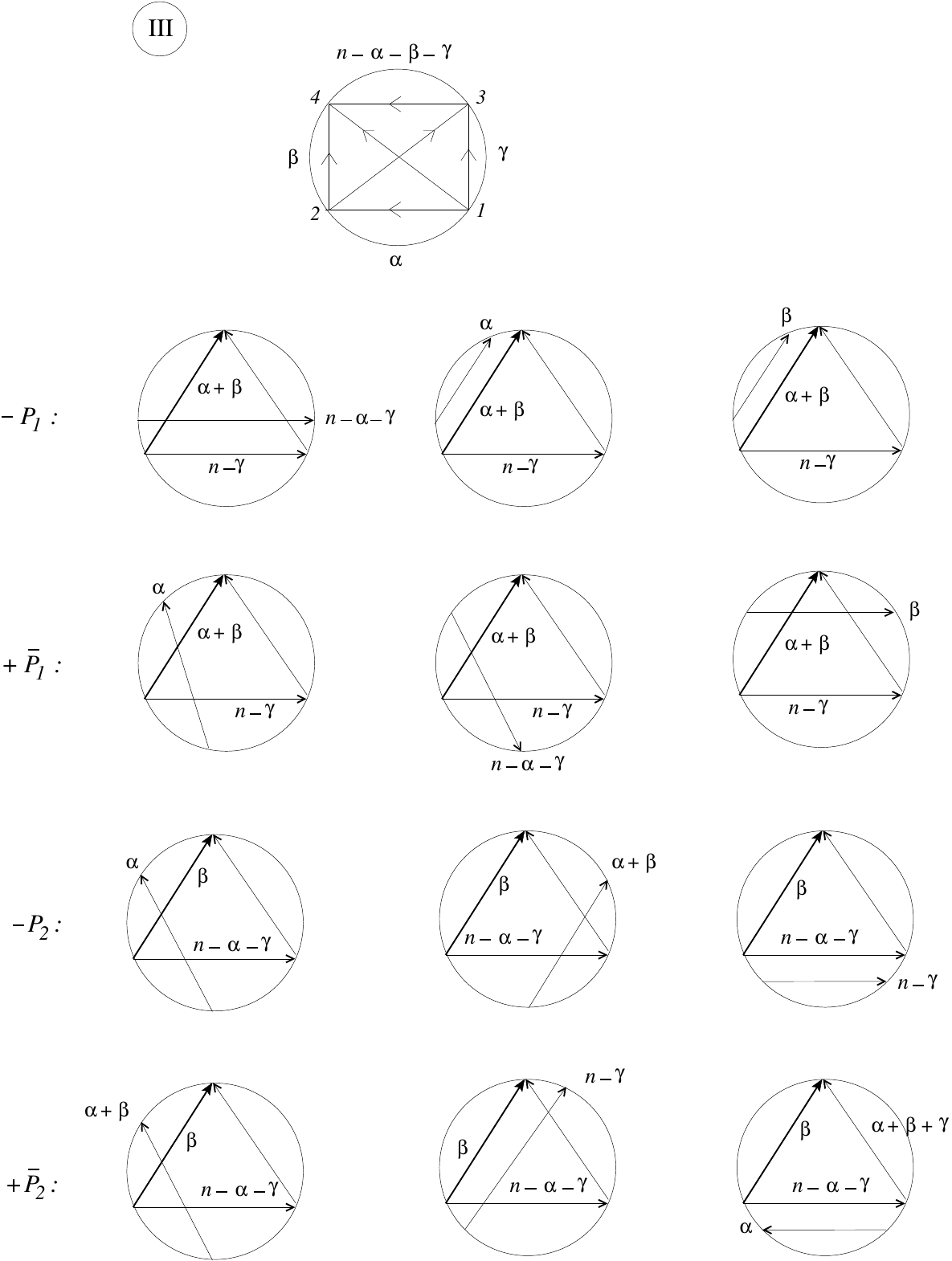}
\caption{\label{ser5} The moves of type $r$ for the global type III}  
\end{figure}

\begin{figure}
\centering
\includegraphics{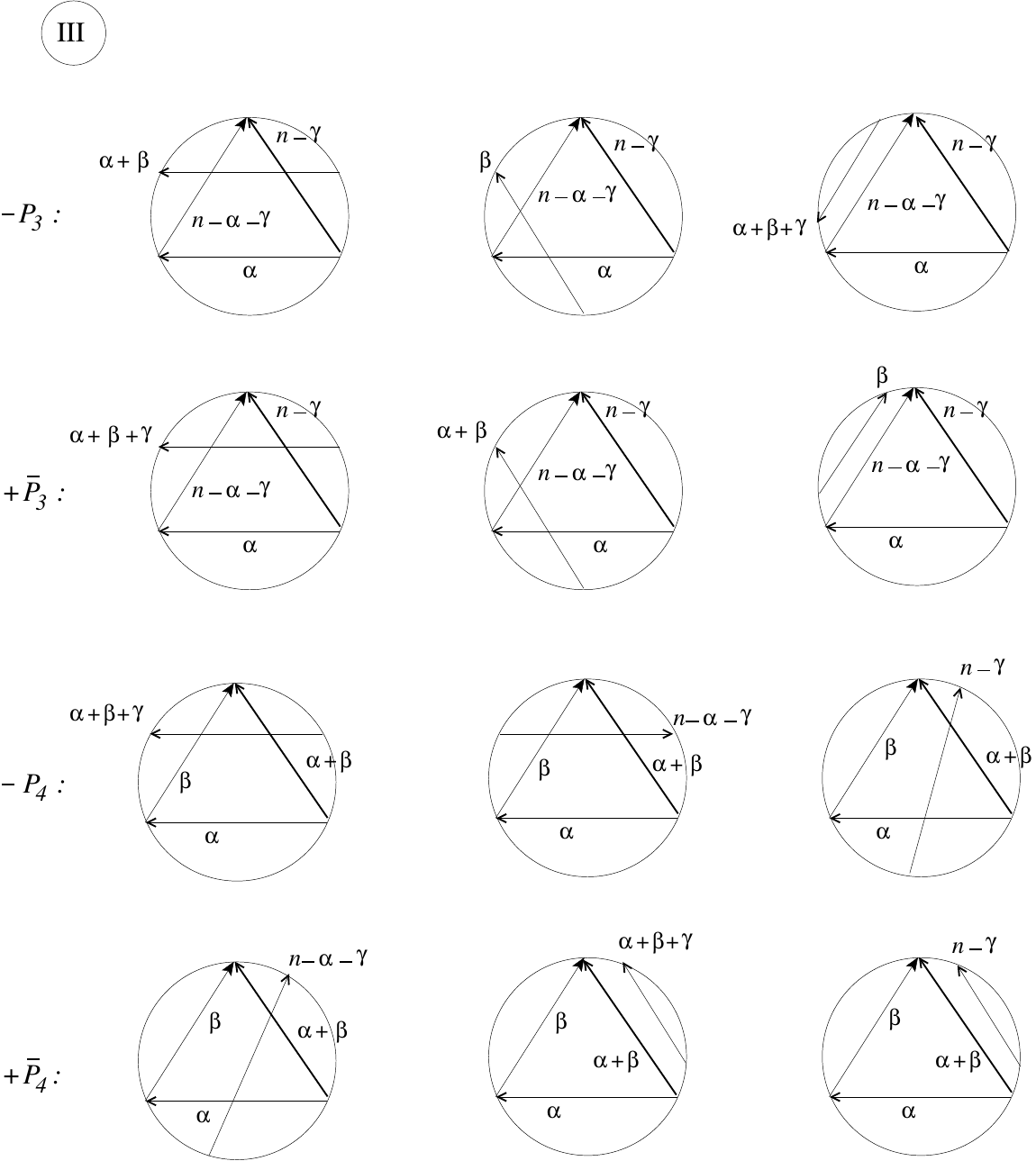}
\caption{\label{ser6} The moves of type $l$ for the global type III}  
\end{figure}

\begin{figure}
\centering
\includegraphics{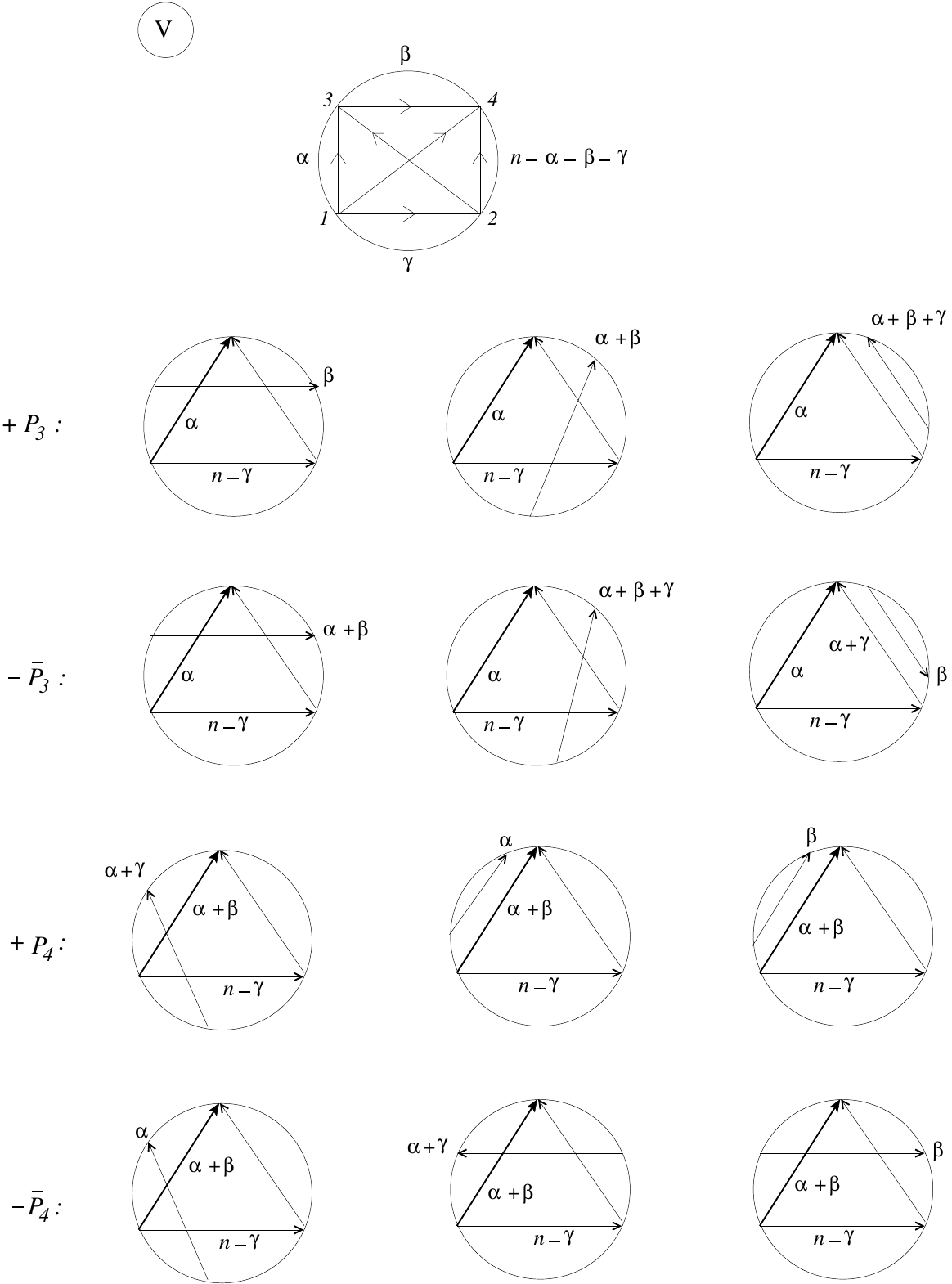}
\caption{\label{ser7} The moves of type $r$ for the global type V}  
\end{figure}

\begin{figure}
\centering
\includegraphics{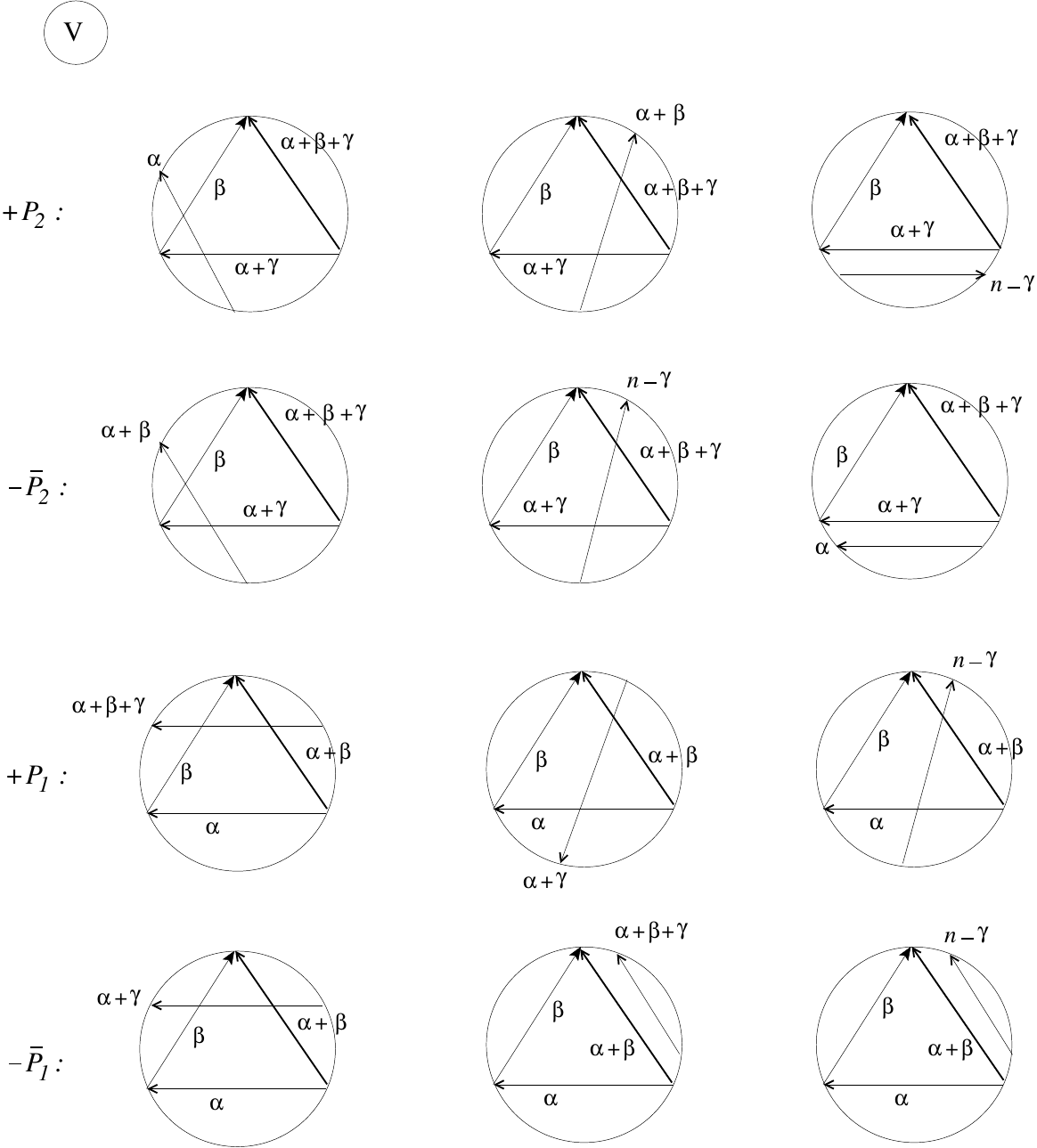}
\caption{\label{ser8} The moves of type $l$ for the global type V}  
\end{figure}

\begin{figure}
\centering
\includegraphics{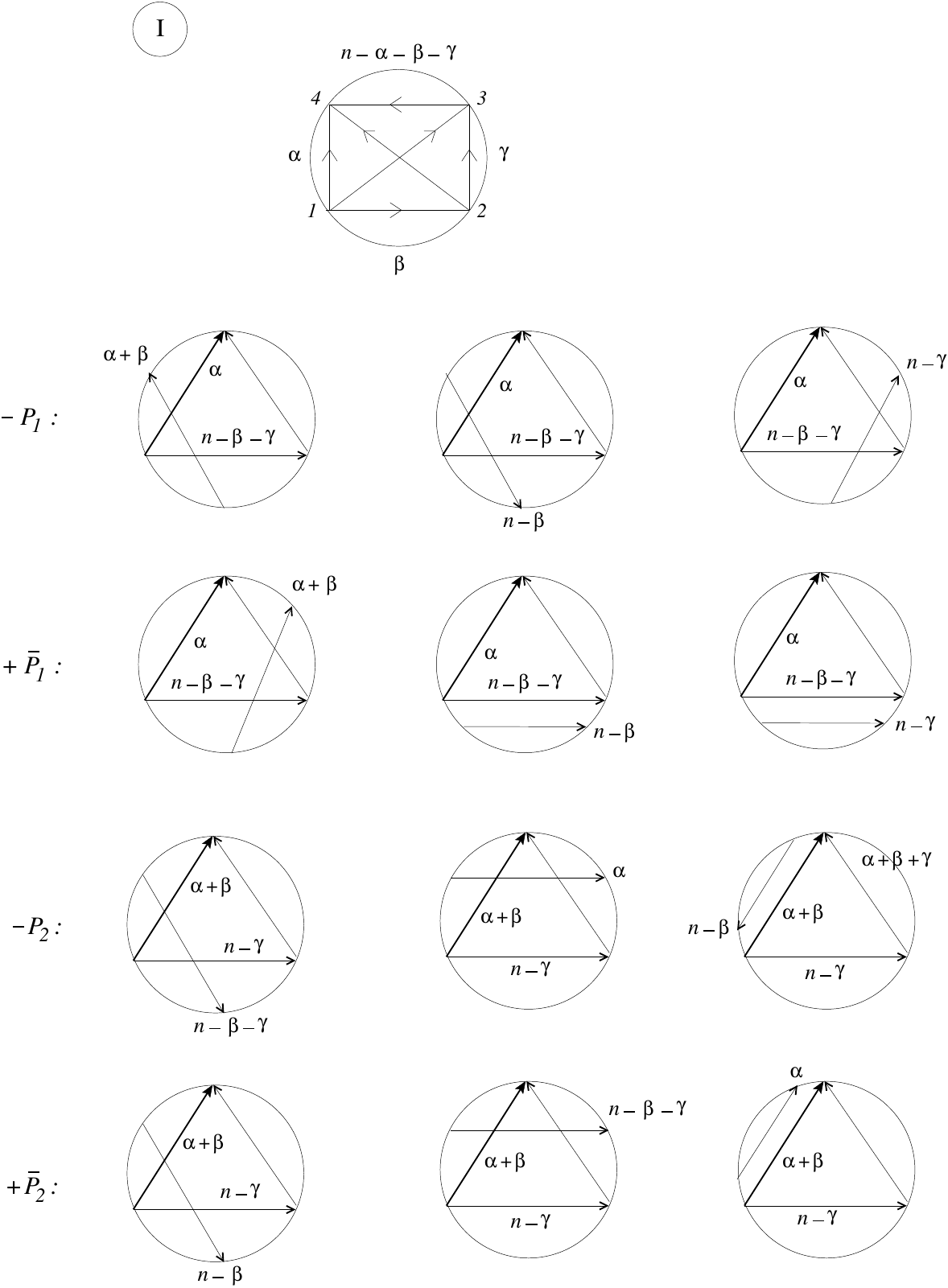}
\caption{\label{ser9} The moves of type $r$ for the global type I}  
\end{figure}

\begin{figure}
\centering
\includegraphics{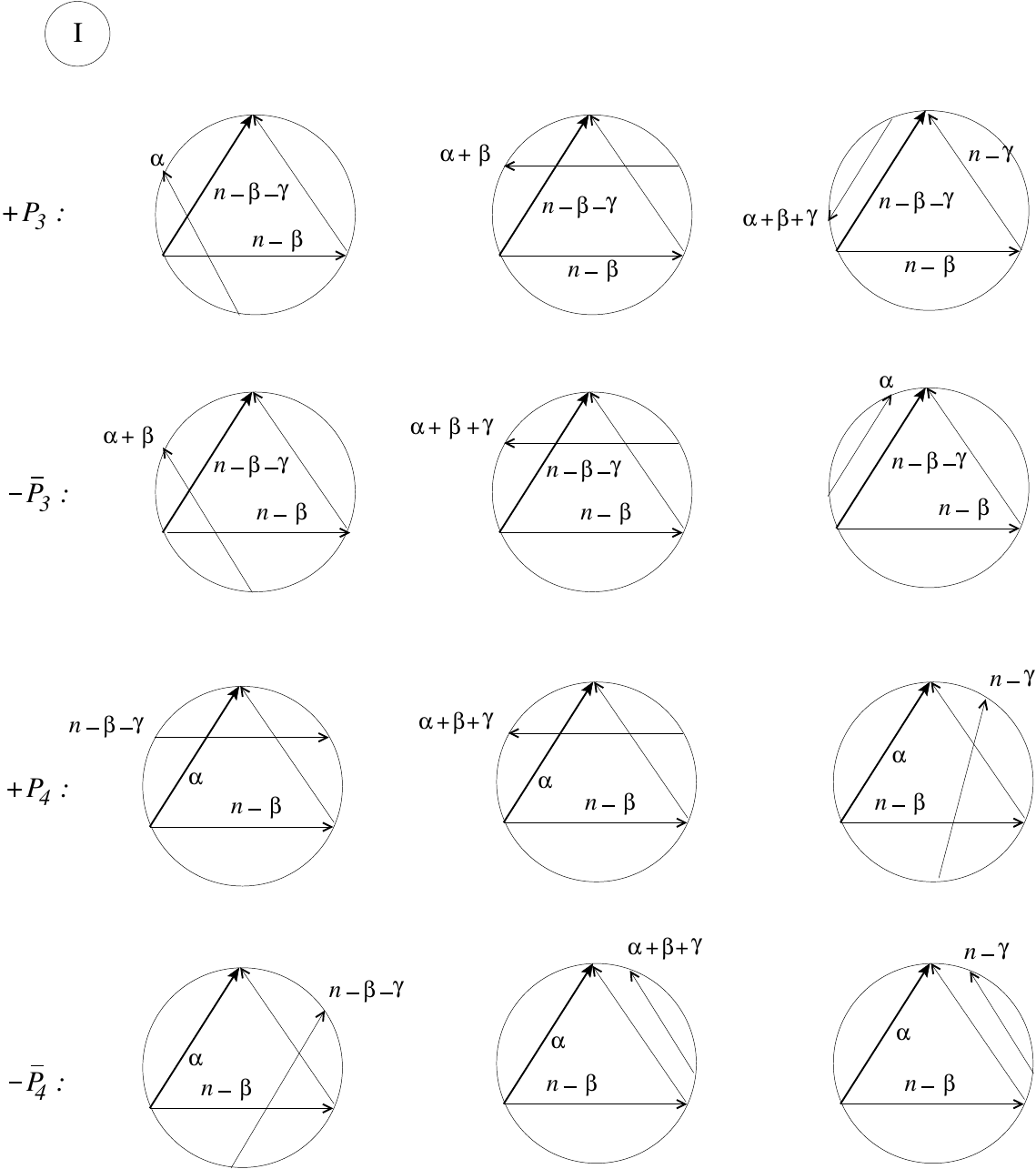}
\caption{\label{ser10} The remaining moves of type $r$ for the global type I}  
\end{figure}

\begin{figure}
\centering
\includegraphics{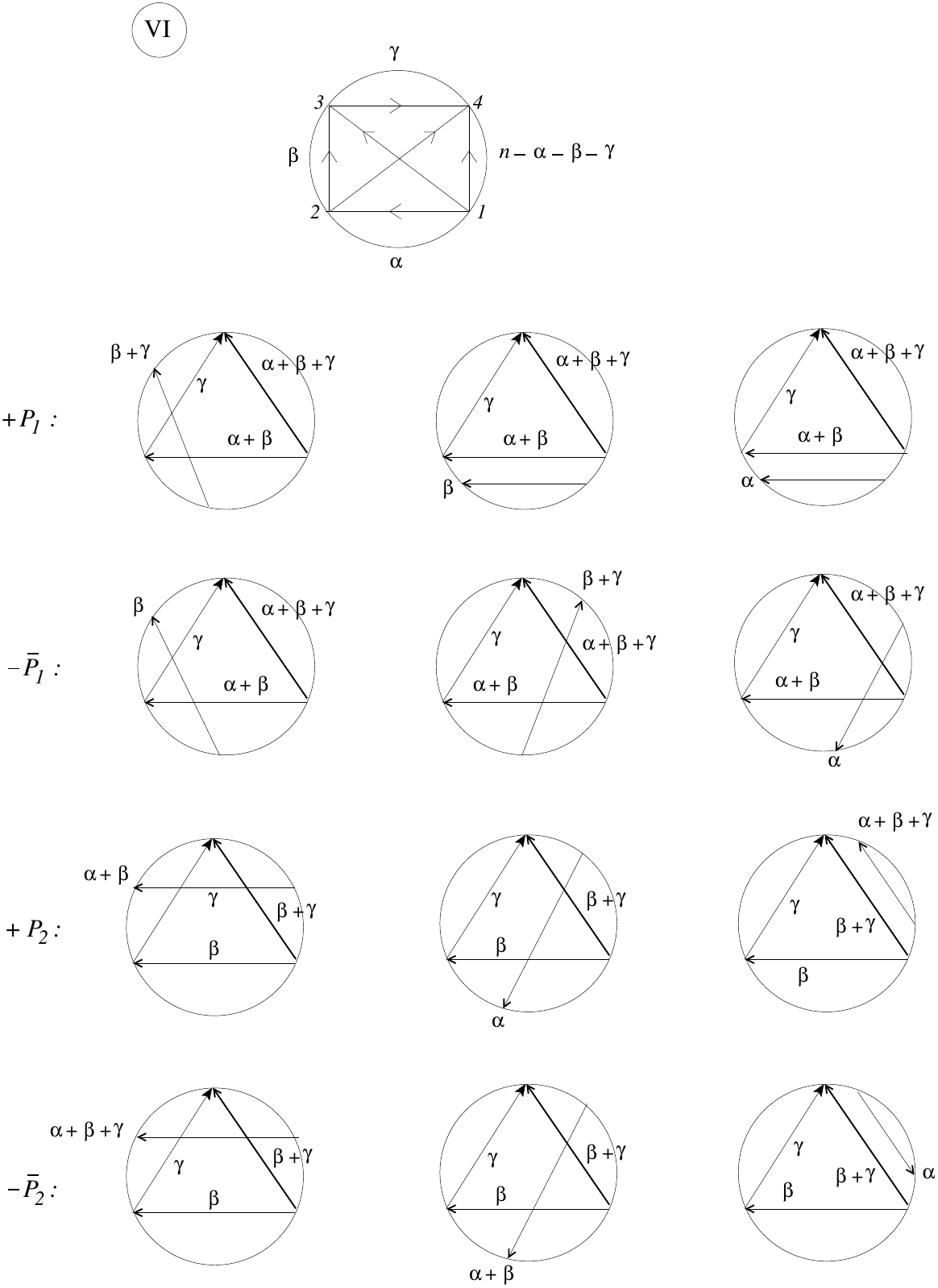}
\caption{\label{ser11} The moves of type $l$ for the global type VI}  
\end{figure}

We have to show now (b), that $R_a^{(2)}$ without the correction term 

$(w(hm)-1)W_2(hm)w(hm)$ (which vanishes automatically for positive triple crossings) vanishes on the meridian of each global positive quadruple crossing, i.e. it satisfies the positive global tetrahedron equation.
\vspace{0,2cm}

The eight strata in the meridian come in pairs with different signs: $P_i$ and $\bar P_i$. Notice, that in the case of positive triple crossings $p$ (i.e. of local type $1$) and of global type $r$ (i.e. the arrow $ml$ goes to the right) the sign of $p$ coincides with the local sign given by "$\sigma_1 \sigma_2 \sigma_1 \to \sigma_2 \sigma_1 \sigma_2$ is the positive direction", compare \cite{F2} (for the global type $l$ it is the opposite). The crossings $d$, $hm$ and $ml$ are of course the same for $P_i$ and $\bar P_i$. These two strata differ with respect to the triangle by exactly the three remaining crossings from the quadruple crossing. But only $n$-crossings with the foot in $dml$ contribute to $W_2(p)$. In particular, we have therefore to study when the foot of an $n$-crossing has slide over the head of the crossing $d$ or when the foot of an $n$-crossing has slide over the head of the crossing $ml$ from $P_i$ to $\bar P_i$.
Inspecting Fig.~\ref{Iglob} up to Fig.~\ref{VI2glob}, we see that we have to study only the following crossings:
\vspace{0,2cm}

{\em $P_3$: foot of $34$ in $dml$, $\bar P_3$: foot of $34$ in $hm$
 
$P_3$: foot of $24$ in $hm$, $\bar P_3$: foot of $24$ in $dml$

$P_4$: foot of $23$ in $dml$, $\bar P_4$: foot of $23$ in $hm$}
\vspace{0,2cm}

Notice, that this never happens for $P_1$ and $\bar P_1$ and for $P_2$ and $\bar P_2$. Indeed, for $P_2$ and $\bar P_2$ the foots stay all in the same region and  for $P_1$ and $\bar P_1$ a foot can only slide over a foot of $ml$, which coincides with the foot of $d$, and hence it stays in the region $dml$. 
\vspace{0,2cm}

In each of the six global types of positive quadruple crossings we have now to consider how $W_2(p)$ changes and how this changing is compensated by a configuration with $W_2(hm)$ for another stratum. But we have also to study the possible changing of $W_2(p)$, $W_2(hm)$ and $l(p)$ from $P_i$ to $\bar P_i$ because of different positions of $0$-crossings and $n$-crossings in the weights. We do all this by inspecting simultaneously a figure from Fig.~\ref{Iglob} up to Fig.~\ref{VI2glob} and the corresponding figure from Fig.~\ref{ser1} up to Fig.~\ref{ser12}. We are conscious that this part of the proof will give a hard time to the interested reader. We apologize for that.
\vspace{0,2cm}

Using  Fig.~\ref{ser1} up to Fig.~\ref{ser12} we go through the list of all cases where at least one of the strata $P_i$ is of the global type $r(a,n,a)$ or $r(n,n,n)$. 

\begin{remark}
 The edges of the tetrahedron correspond to the crossings and the 2-faces correspond to the triple crossings. If there is one 2-face corresponding to $r(a,n,a)$ or $r(n,n,n)$ then the remaining three edges have a common vertex. The homology class of one of these edges (i.e. the homological marking of the corresponding crossing) determines now the homology classes of the remaining two edges. Of course, the homology classes of three edges with a common vertex determine the  homology classes of the remaining three edges.

\end{remark}
 
We call a tetrahedron {\em generic} if it contains one stratum of triple crossings with all markings in $\{0,n\}$ and there is at least one crossing of marking $a$ (there are then automatically two other crossings with markings in $\{a,n-a\}$ and there are no other triple crossings with all markings in $\{0,n\}$).

We call the tetrahedron {\em degenerate} if all six crossings have their markings in $\{0,n\}$. 

{\em Our strategy is the following: we solve the generic tetrahedron equations by using weights of order  $2$. We take then this solution and possibly complete it to a solution of the degenerate tetrahedron equations by using weights of lower order.}

It turns out that for the 1-cocycle $R_a^{(2)}$ actually we do not have to complete the solution in the degenerate case, but just to check that it is still a solution without adding new strata or new weights.
\vspace{0,2cm}

{\em Let us study now the generic tetrahedrons.} We have only to consider the cases when there is at least one stratum of type $r(a,n,a)$ or $r(n,n,n)$ and moreover, there is at least one other $n$-crossing or $(n-a)$-crossing (which could change $l(p)$ or which could change the fact to be an $f$-crossing). This gives us conditions on the parameters $\alpha$, $\beta$ and $\gamma$. The weight $W_2(p)$ can only change if there is at least one other stratum of type $l(n,0,n)$ as was shown in Fig.~\ref{globRIIIn} and  Fig.~\ref{indweight}. If $P_i$ has not the right global type in order to contribute, or if $l(P_i)=l(\bar P_i)$, or if $W_2(P_i)=W_2(\bar P_i)$, or if $W_2(hm)(P_i)=W_2(hm)(\bar P_i)$, then we will simply skip them in the proof (because always $sign(P_i)=-sign(\bar P_i)$ and they cancel out together).
\vspace{0,2cm}

{\em Global type I}
\vspace{0,2cm}

All global types of triple crossings in the meridian are of type $r$ in this case. Consequently, the individual contribution of an $n$-crossing to $(n,0)$ cannot change in the meridian, because it changes only by passing a triple crossing of type $l(n,0,n)$, compare Section 4.2.
\vspace{0,2cm}

$I_1$: $\alpha = a$, $\beta = n-a$ and hence $\gamma = 0$. Then $P_1$, $P_3$ and $P_4$ are of type $r(a,n,a)$ and $P_2$ is of type $r(n,n,n)$. In this case $23$, $24$ and $34$ are all $n$-crossings. $l(P_1)=l(\bar P_1)+1$. There aren't any $(n-a)$-crossings and hence $-l(P_2)+l(\bar P_2)=0$. There aren't any $0$-crossings and hence all weights $W_2$ of crossings stay the same.

$P_3-\bar P_3=W_2(34)-W_2(24)$. The crossing $34$ is the crossing $hm$ in $P_1$. $-l(P_1)+l(\bar P_1)=-1$ and hence $W_2(34)$ cancels out in $R_a^{(2)}$. The crossing $24$ is $hm$ in $P_4$. $l(P_4)-l(\bar P_4)=1$ and hence $W_2(24)$ cancels out in $R_a^{(2)}$ too.

$P_4-\bar P_4=W_2(23)$. The crossing $23$ is $hm$ in $P_3$. $l(P_3)-l(\bar P_3)=-1$ and hence $W_2(23)$ cancels out in $R_a^{(2)}$ too.
\vspace{0,2cm}

$I_2$: $\alpha = a$, $\beta = 0$ and $\gamma = n-a$. Then $P_1$ and $P_2$  are of type $r(a,n,a)$. There aren't any $f$-crossings which could change, because they can appear only for $P_3$ and $P_4$. There aren't any  $0$-crossings and all weights stay the same. $-l(P_1)+l(\bar P_1)=0$ and $-l(P_2)+l(\bar P_2)=0$. Hence, $R_a^{(2)}=0$ on the meridian.
\vspace{0,2cm}

$I_3$: $\alpha = 0$, $\beta = a$ and $\gamma = n-a$. Then $P_2$  is of type $r(a,n,a)$. $-l(P_2)+l(\bar P_2)=0$ and $R_a^{(2)}=0$ on the meridian.
\vspace{0,2cm}

$I_4$: $\alpha = n-a$, $\beta = a$ and hence $\gamma = 0$. Then $P_2$  is of type $r(n,n,n)$. $-l(P_2)+l(\bar P_2)=0$ and $R_a^{(2)}=0$ on the meridian.
\vspace{0,2cm} 

$I_5$: $\alpha = n$ and hence $\beta = 0$ and $\gamma = 0$. Then all four $P_i$ are of type $r(n,n,n)$. There aren't any $(n-a)$-crossings and hence $l(P_i)-l(\bar P_i)=0$ for all four $i$ and $R_a^{(2)}=0$ on the meridian.
\vspace{0,2cm} 

$I_6$: $\alpha = 0$, $\beta = n$ and hence $\gamma = 0$. Then $P_2$  is of type $r(n,n,n)$. $-l(P_2)+l(\bar P_2)=0$ and $R_a^{(2)}=0$ on the meridian.
\vspace{0,2cm} 

$I_7$: $\alpha =$ arbitrary, $\beta = 0$ and $\gamma = 0$. Then $P_3$  is of type $r(n,n,n)$. If $\alpha = n-a$ then $[24]=[14]$ and they contribute simultaneously. Consequently,  $l(P_3)-l(\bar P_3)=0$ and $R_a^{(2)}=0$ on the meridian.
\vspace{0,2cm}

{\em Global type II}
\vspace{0,2cm}

Here only the strata $P_1$ and $P_4$ could contribute.
\vspace{0,2cm}

$II_1$: $\alpha = a$, $\beta = 0$ and $\gamma = n-a$. Then $P_1$ and $P_4$ are of type $r(a,n,a)$. But $P_2$ is of type $l(n,n,0)$ and the weights cannot change (they can change only for $l(n,0,n)$). $-l(P_1)+l(\bar P_1)=0$ because $[12]=a$ and the contributions of $-P_1+\bar P_1$ cancel out. $[23]=0$ and hence $f$-crossings do not change. It follows that $R_a^{(2)}=0$ on the meridian.
\vspace{0,2cm} 

$II_2$: $\alpha = a$, $\beta = n-a$ and hence $\gamma = 0$. Then $P_4$ is of type $r(a,n,a)$. There is no stratum $l(n,0,n)$ and the weights cannot change. $[23]=n-a$ and $f$-crossings do not change. $l(P_4)-l(\bar P_4)=0$ and hence $R_a^{(2)}=0$ on the meridian.
\vspace{0,2cm} 

$II_3$: $\alpha = n$, $\beta = 0$ and $\gamma = 0$. Then $P_1$ and $P_4$ are of type $r(n,n,n)$. $P_2$ is of type $l(n,0,n)$. But there is no $0$-crossing (from the edges of the tetrahedron) which cuts $hm=34$ in $P_1$ or $\bar P_1$. There is also no $0$-crossing which cuts $hm=24$ in $P_4$ or $\bar P_4$. Hence the weights do not change. There are no $(n-a)$-crossings at all and $R_a^{(2)}=0$ on the meridian.
\vspace{0,2cm} 

{\em Global type III}
\vspace{0,2cm}

Here only the strata $P_1$ and $P_2$ could contribute. Consequently, $f$-crossings cannot change.
\vspace{0,2cm}

$III_1$: $\alpha = a$, $\beta = 0$ and $\gamma = n-a$. Then $P_1$ is of type $r(a,n,a)$. There is no stratum $l(n,0,n)$ and the weights cannot change. $[12]=a$ and hence $-l(P_1)+l(\bar P_1)=0$ and $R_a^{(2)}=0$ on the meridian.
\vspace{0,2cm} 

$III_2$: $\alpha = 0$, $\beta = a$ and $\gamma = n-a$. Then $P_1$ and $P_2$ are of type $r(a,n,a)$. There is no stratum $l(n,0,n)$ and the weights cannot change. $[12]=0$ and hence $-l(P_1)+l(\bar P_1)=0$. $[12]=0$, $[13]=a$ and  hence $-l(P_2)+l(\bar P_2)=0$ and $R_a^{(2)}=0$ on the meridian.
\vspace{0,2cm} 

$III_3$: $\alpha = n-a$, $\beta = a$ and $\gamma = 0$. Then $P_1$ is of type $r(n,n,n)$ and $P_2$ is of type $r(a,n,a)$. 
There is no stratum $l(n,0,n)$ and the weights cannot change. But here an unexpected phenomenon arrives: $[13]=n$ and it follows that $-l(P_2)+l(\bar P_2)=1$ (this is the {\em only} case where this happens). Consequently, $-P_2+\bar P_2$ contributes $W_2(hm)=W_2(34)$. But $34$ is also the crossing $hm$ in $P_1$, which is of type $r(n,n,n)$. $[12]=n-a$ and hence $-l(P_1)+l(\bar P_1)=1$. However, the strata $r(n,n,n)$ enter into the formula with a negative sign and it follows  that 
$R_a^{(2)}=0$ on the meridian.

{\em This is the only case where the contributions of strata $r(a,n,a)$ are compensated by the contributions of strata $r(n,n,n)$!}
\vspace{0,2cm}

$III_4$: $\alpha = a$, $\beta = n-a$ and $\gamma = 0$. Then $P_1$ is of type $r(n,n,n)$. There is no stratum $l(n,0,n)$ and the weights cannot change. $[12]=a$ and we have to distinguish two cases. If $a \not = n-a$ then $-l(P_1)+l(\bar P_1)=0$ and they cancel out together. If $a  = n-a$ then $-l(P_1)+l(\bar P_1)=1$ and they contribute $-W_2(34)$  to $R_a^{(2)}$. But in this case $P_2$ is of type $r(n-a=a,n,n-a=a)$. $[13]=n$ and $-l(P_2)+l(\bar P_2)=1$. Hence, $-P_2+\bar P_2$ contribute $W_2(34)$ and $R_a^{(2)}=0$ on the meridian.
\vspace{0,2cm} 

{\em Global type IV}
\vspace{0,2cm}

Here only the strata $P_2$ and $P_3$ could contribute. 
\vspace{0,2cm}

$IV_1$: $\alpha = 0$, $\beta = a$ and $\gamma = 0$. Then $P_2$ is of type $r(a,n,a)$. There is no stratum $l(n,0,n)$ and the weights cannot change. $-l(P_2)+l(\bar P_2)=0$ and $R_a^{(2)}=0$ on the meridian.
\vspace{0,2cm} 

$IV_2$: $\alpha = a$, $\beta = n- a$ and $\gamma = 0$. Then $P_3$ is of type $r(a,n,a)$. There is no stratum $l(n,0,n)$ and the weights cannot change. $[34]=[24]=n-a$ and $f$-crossings cannot change.  $l(P_3)-l(\bar P_3)=0$ and $R_a^{(2)}=0$ on the meridian.
\vspace{0,2cm} 

$IV_3$: $\alpha = a$, $\beta = 0$ and $\gamma = n-a$. Then $P_3$ is of type $r(a,n,a)$. There is no stratum $l(n,0,n)$ and the weights cannot change. $[34]=[24]=0$ and $f$-crossings cannot change. $l(P_3)-l(\bar P_3)=0$ and $R_a^{(2)}=0$ on the meridian.
\vspace{0,2cm} 

$IV_4$: $\alpha = 0$, $\beta = n$ and $\gamma = 0$. Then $P_2$ is of type $r(n,n,n)$. There is no stratum $l(n,0,n)$ and the weights cannot change. $-l(P_2)+l(\bar P_2)=0$ and $R_a^{(2)}=0$ on the meridian.
\vspace{0,2cm} 

{\em Global type V}
\vspace{0,2cm}

Here only the strata $P_3$ and $P_4$ could contribute. 
\vspace{0,2cm}

$V_1$: $\alpha = a$, $\beta = 0$ and $\gamma = n-a$. Then $P_3$ and $P_4$ are of type $r(a,n,a)$. $P_2$ is of type $l(n,0,n)$. $[23]=[24]=n$ and hence both $f$-crossings change. We write the contribution of each of the four strata.

$P_3$: $W_2(P_3) + l(P_3)W_2(hm=23)$

$\bar P_3$: $-W_2(P_3)-W_2(24)-(l(P_3)+1)(W_2(23)-1)$. The correction $-1$ comes from the fact that $[34]=0$ and $34$ cuts $hm$ in $P_3$ but not in $\bar P_3$.

$P_4$: $W_2(\bar P_4)+W_2(23)+(l(\bar P_4)+1)((W_2(hm=24)-1)$. The correction $-1$ comes from the fact that $34$ does not cut $24$. Notice that $l(P_3)=l(\bar P_4)$ because they share the same crossing $ml=12$.

$\bar P_4$: $-W_2(\bar P_4)-l(\bar P_4)(W_2(hm=24)$

It follows that the total contribution of the four strata is:
\vspace{0,2cm}

$W_2(P_3)+l(P_3)W_2(23)-W_2(P_3)-W_2(24)-l(P_3)W_2(23)+l(P_3)-W_2(23)+1+W_2(\bar P_4)+W_2(23)+l(P_3)W_2(24)-l(P_3)+W_2(24)-1-W_2(\bar P_4)-l(P_3)W_2(24)=0$ (We are really lucky here!)
\vspace{0,2cm}

$V_2$: $\alpha = 0$, $\beta = a$ and $\gamma = n-a$. Then $P_4$ is of type $r(a,n,a)$. There is no stratum $l(n,0,n)$ and the weights cannot change. $[23]=n-a$ and hence $f$-crossings do not change. $[23]=n-a$, $[13]=0$ and hence $l(P_4)-l(\bar P_4)=0$ and $R_a^{(2)}=0$ on the meridian.
\vspace{0,2cm}  

$V_3$: $\alpha = n$, $\beta = 0$ and $\gamma = 0$. Then $P_3$ and $P_4$ are of type $r(n,n,n)$. $P_1$ and $P_2$ are both of type $l(n,0,n)$. There are no $(n-a)$-crossings at all and $ml=12$ for both $P_3$ and $P_4$. It follows that $l(P_3)=l(\bar P_3)=l(P_4)=l(\bar P_4)$. Hence we have only to check that the sum over the four strata of $W_2(hm)$ is $0$.

$W_2(hm)(P_3)-W_2(hm)(\bar P_3)=1$ but $W_2(hm)(P_4)-W_2(hm)(\bar P_4)=-1$ (in both cases the difference comes from the $0$-crossing $34$) and they cancel out together.
\vspace{0,2cm}

$V_4$: $\alpha +\beta = n$ and hence $\gamma = 0$. Then $P_4$ is of type $r(n,n,n)$. We can assume that $\beta \not= 0$ because this was the case $V_3$. Consequently, there is no stratum $l(n,0,n)$ and the weights cannot change. We have $[23]=[13]$ because $\gamma =0$. It follows that $l(P_4)-l(\bar P_4)=0$ and hence $R_a^{(2)}=0$ on the meridian.
\vspace{0,2cm}

{\em Global type VI}
\vspace{0,2cm}

All strata are of type $l$ and there is no contribution at all.
\vspace{0,2cm}

\vspace{0,2cm}
We have proven that our 1-cochain $R_a^{(2)}$ satisfies (b): the positive global tetrahedron equation in the generic case. This was the hardest part!
\vspace{0,2cm}

{\em Let us study now the degenerate tetrahedrons.} First of all we observe that if a stratum $P_i$ is of type $r(n,n,n)$ then $l(P_i)=l(\bar P_i)$, because there are no $n-a$ -crossings at all in the tetrahedron. Moreover, $l(p)$ is completely determined by the position with respect to the crossing $ml$ alone, compare Fig.~\ref{lnnn}. Consequently, if two triple crossings $r(n,n,n)$ share the same crossing $ml$ then their refined linking numbers $l(p)$ are the same.

Let $hm$ be a crossing in $r(n,n,n)$. Its weight $W_2(hm)$ can {\em only} change if $hm$ is also a crossing in a stratum of type $l(n,0,n)$, compare Fig.~\ref{indweight}. The homological markings of the edges of the tetrahedron determine in a unique way the point at infinity and vice versa in the case of a degenerate tetrahedron (this is of course not the case for the non-degenerate tetrahedron). We examine now Fig.~\ref{globquad} and we observe that the crossing $hm$ in a stratum $r(n,n,n)$ is also a crossing in a stratum $l(n,0,n)$ exactly once, namely for the global type $V$ and when the point at infinity is on the arc from branch 3 to branch 1. This corresponds to $2=\infty$ in the figures Fig.~\ref{Vglob} and Fig.~\ref{V2glob}. There are exactly two strata of type $r(n,n,n)$, namely $P_3$ and $P_4$ and exactly two strata of type $l(n,0,n)$, namely $P_1$ and $P_2$. The strata $P_3$ and $P_4$ share the crossing $12$ as crossing $ml$ and hence $P_3$, $\bar P_3$, $P_4$ and $\bar P_4$ have all the same refined linking number $l(p)$. Inspecting Fig.~\ref{Vglob} and Fig.~\ref{V2glob} we see that $P_3-\bar P_3$ contributes $+1$ to $W_2(hm)$ and 
that $P_4-\bar P_4$ contributes $-1$ to $W_2(hm)$. Consequently $R_a^{(2)}$ vanishes on the meridian. 

We have proven that $R_a^{(2)}$ satisfies also the degenerate tetrahedron equations.
$\Box$

\subsection{$R_a^{(2)}$ satisfies the cube equations (c)}

\begin{figure}
\centering
\includegraphics{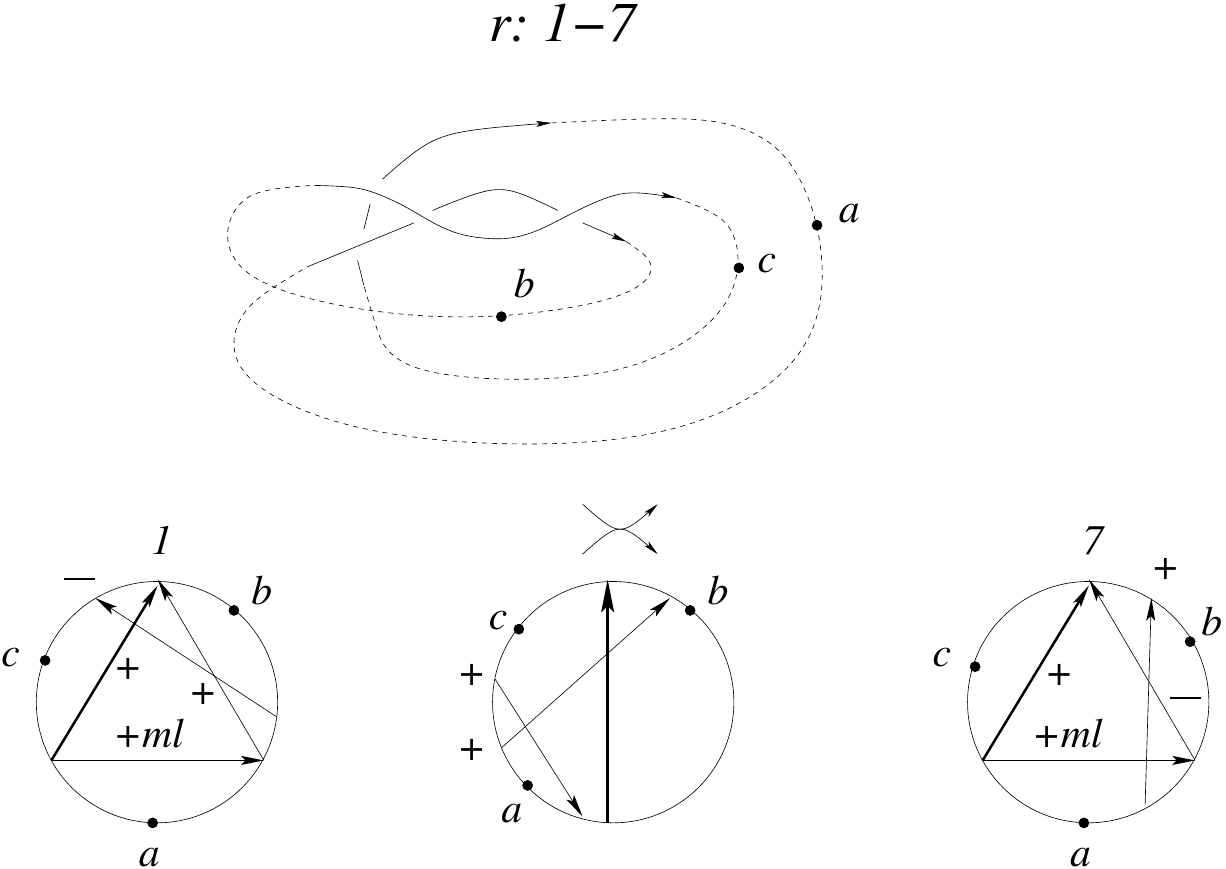}
\caption{\label{r1-7} $r1-7$ }  
\end{figure}

\begin{figure}
\centering
\includegraphics{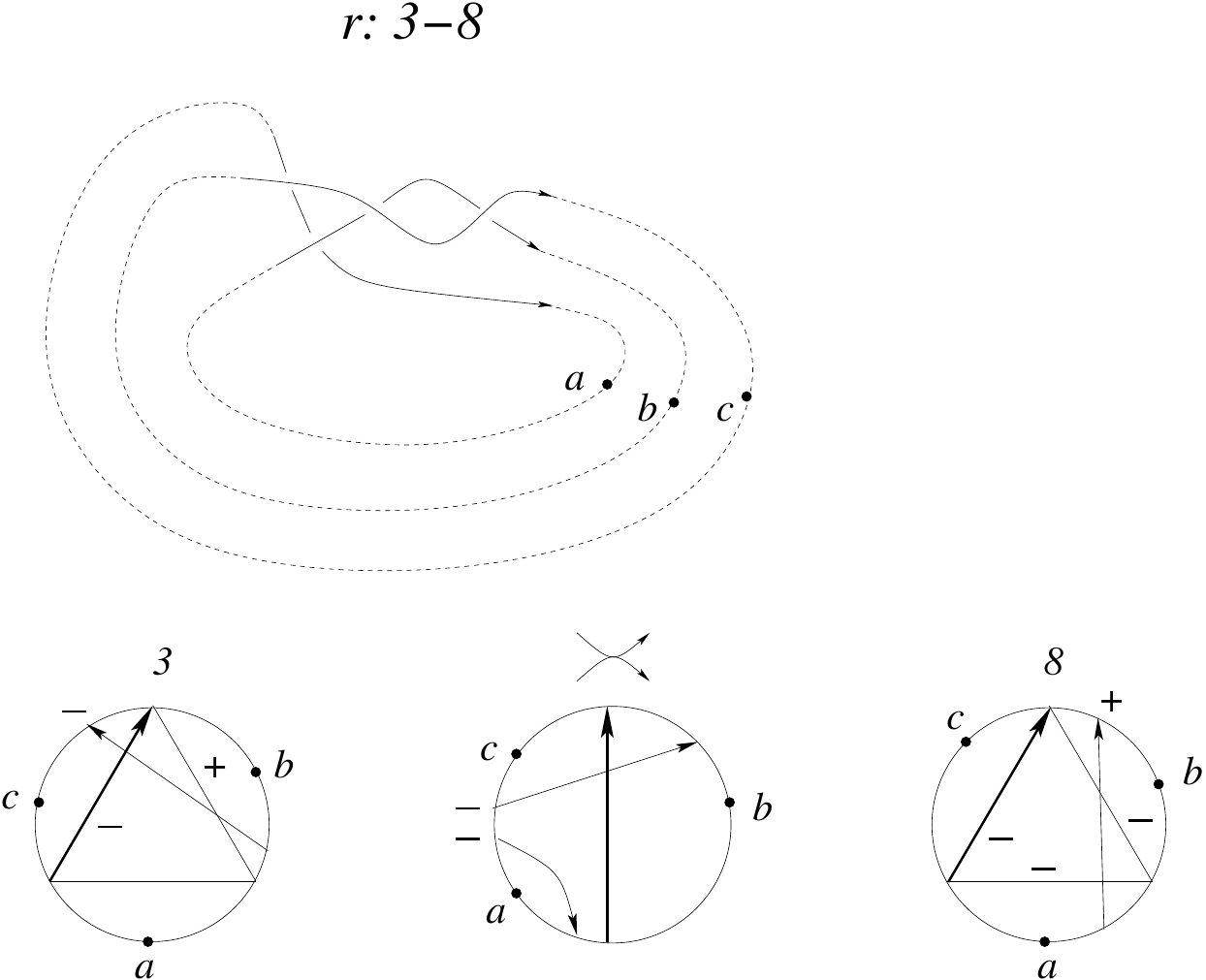}
\caption{\label{r3-8}  $r3-8$}  
\end{figure}

If a loop in $M_n^{semi-reg}$ passes with an ordinary homotopy through a triple crossing where two branches are ordinary tangential then in the R III move exactly one crossing is replaced by another crossing with the same homological marking  but with the opposite sign.
Evidently, we have only to study the case when this crossing is the $n$-crossing $hm$ for the type $r(a,n,a)$ or $r(n,n,n)$. We show the corresponding edges of the graph $\Gamma$ in Fig.~\ref{r1-7} up to Fig.~\ref{r3-8}.  For the numbers of the local types of triple crossings compare Fig.~\ref{loctricross}. The positive triple crossing correspond to the type $1$. As already mentioned, the global type of two triple crossings of an edge is always the same.

\begin{figure}
\centering
\includegraphics{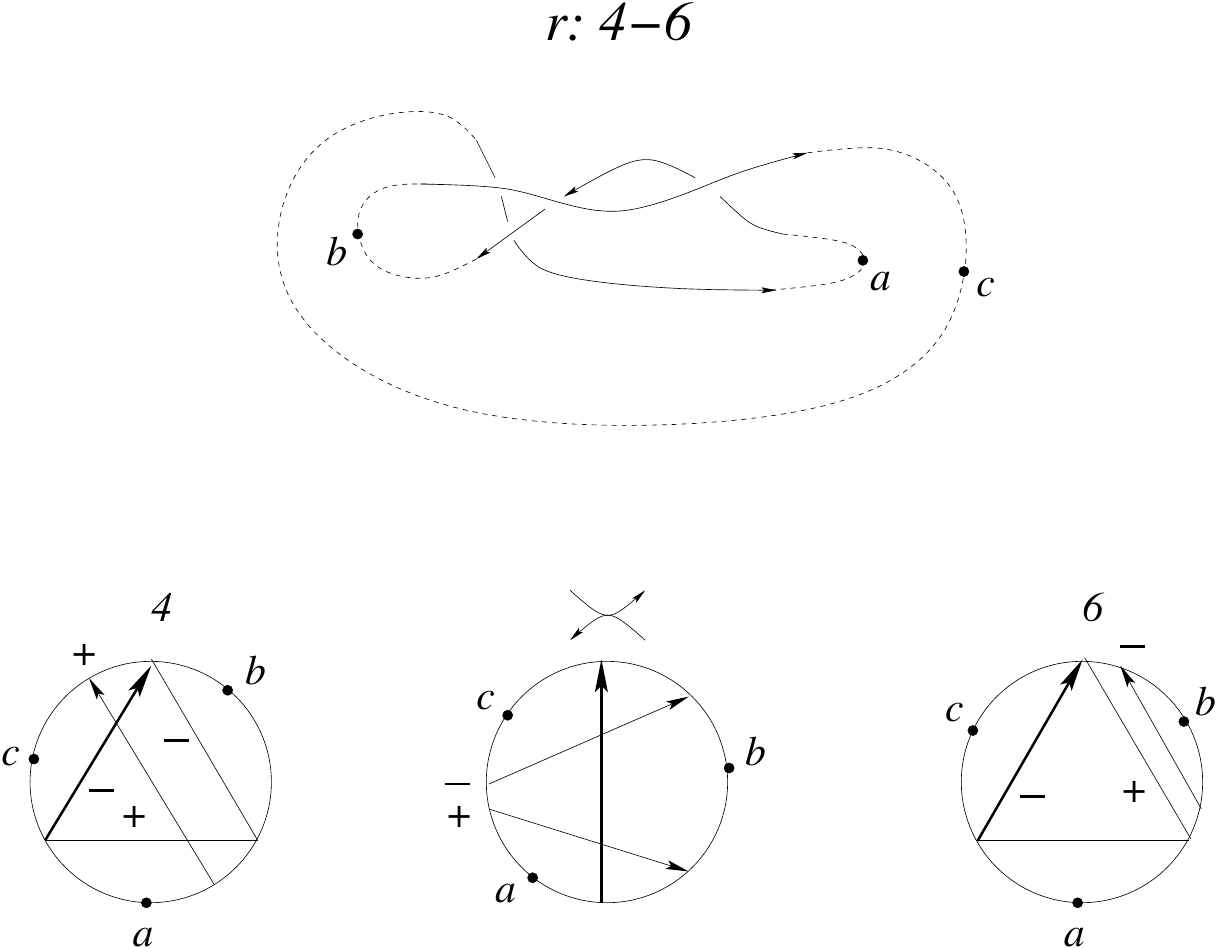}
\caption{\label{r4-6}  $r4-6$}  
\end{figure}

\begin{figure}
\centering
\includegraphics{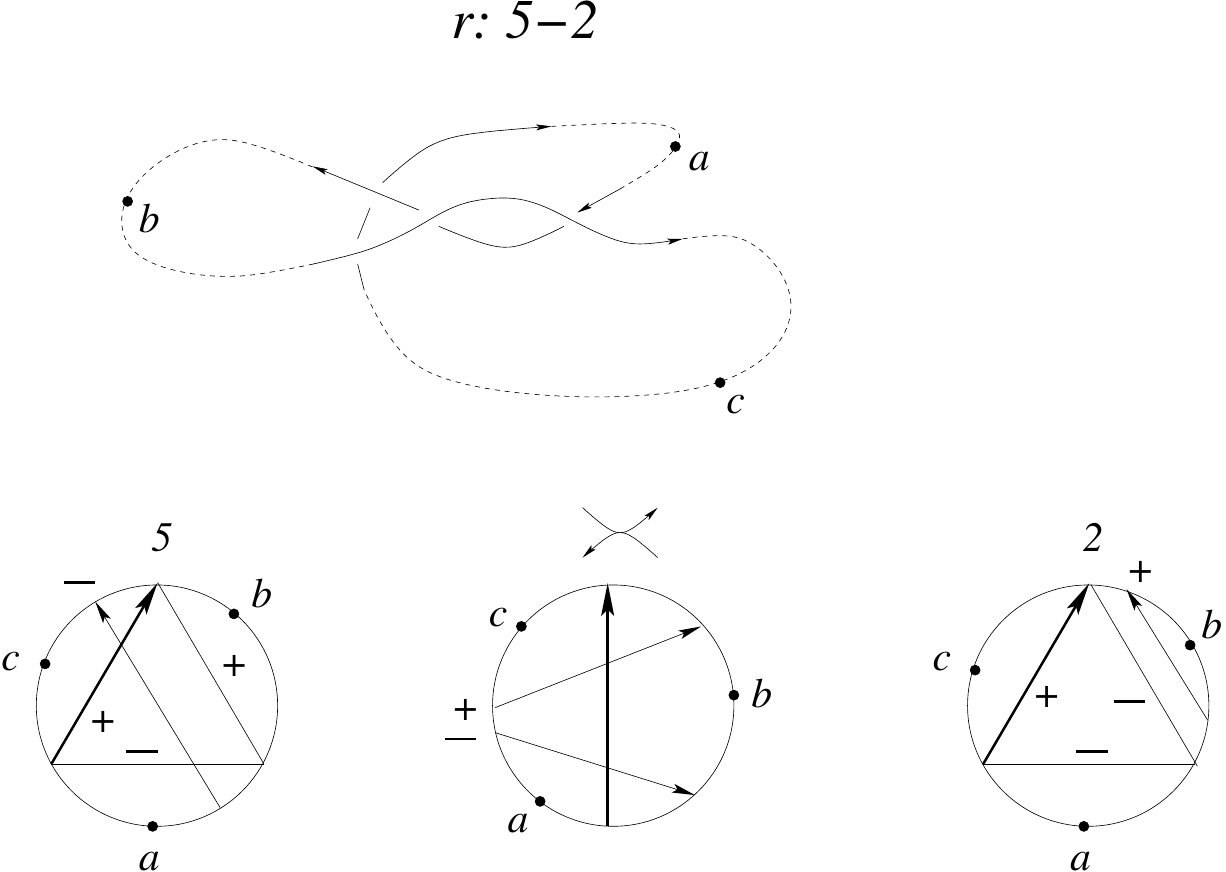}
\caption{\label{r5-2}  $r5-2$}  
\end{figure}

It is very convenient that with our definition the crossing $hm$ enters $W_2(hm)w(hm)$ always with a positive sign. The other new $n$-crossing, lets call it $hm'$, from the R II move does never enter together with $hm$ into a configuration. Consequently, $W_2(hm)$ is the same for the two R III moves in the edge of $\Gamma$.

Let us consider the new crossing $hm'$ from the R II move. The configuration with $hm'$ which has its foot in $dml$ contributes now to $W_2(p)$ by $w(hm')W_2(hm)$ and the linking number $l(p)$ changes by adding $w(hm')$. Going through Fig.~\ref{r1-7} up to Fig.~\ref{r3-8} we see that in each case with the correction term $(w(hm)-1)W_2(hm)w(hm)$ the contributions of the two R III moves in the edge are now the same.

If the crossings from the R II moves do not correspond to the crossing $hm$ then the two crossings $hm$ in the R III moves have the same sign and the correction terms for both R III moves are the same.

For the R III moves of type $r(n,n,n)$ (where we do not consider $W_2(p)$ at all) again $hm$ enters $W_2(hm)w(hm)$ always with a positive sign. Moreover, the new $n$-crossing never enters $l(p)$ because only $(n-a)$-crossings can contribute to $l(p)$ for this type (and $0<a<n$).
$\Box$

\subsection{Proof of Theorem 1}

We have proven that $R_a^{(2)}$ satisfies the equations (a), (b) and (c) and hence it is a 1-cocycle, compare \cite{F2}. The examples in Section 3 show that it represents a non-trivial cohomology class, that it induces a non-trivial pairing  and that it depends non-trivially on the integer parameter $0<a<n$.

\begin{lemma}
$R_a^{(2)}(scan(T \cup nK'))$ is an invariant of $T \cup nK'$ up to semi-regular isotopy.

\end{lemma}

{\em Proof.}

More generally, let $T$ be an arbitrary diagram of an oriented $n$-component string link and let $s$ be a semi-regular isotopy which connects $T$ with a diagram $T'$. We consider the loop $-s \circ $$-scan(T') \circ s \circ scan(T)$ in $M_n^{semi-reg}$, compare Fig.~\ref{scannK} for $scan$. This loop is contractible in $M_n^{semi-reg}$ because $s$ and $scan$ commute, i.e. we can perform them simultaneously. Consequently, $R_a^{(2)}$ vanishes on this loop, because $R_a^{(2)}$ is a 1-cocycle. It suffices to prove now that each contribution of a Reidemeister move $t$ in $s$ cancels out with the contribution of the same move $t$ in $-s$ (the signs of the contributions are of course opposite). The difference for the two Reidemeister moves is in a branch which has moved under $t$. It suffices to study the weights and the linking numbers in the R III moves. The move $t$ corresponds to the stratum $P_2$ respectively $\bar P_2$ in the tetrahedron equation. First we observe, that there are no triple crossings of type $r(n,n,n)$ at all in the scan-arc. Indeed, the $n$-curl does not contain any $n$-crossings (compare  Fig.~\ref{scannK}) and moreover all new crossings from R II moves of the $n$-curl with $T \cup nK'$ have markings $<n$. Hence, we can assume that $t$ is of type $r(a,n,a)$. If $t$ is a positive triple crossing now, then the weights are the same just before the branch moves under $t$ and just after it has moved under $t$. Indeed, this follows from the fact that for the positive global tetrahedron equation the contribution from the stratum $-P_2$ cancels always out with that from the stratum $\bar P_2$ besides that the linking number $l(p)$ could have changed. But it follows from the proof in Section 4.3 that this could only happen if there would be a stratum of type $r(n,n,n)$ in the tetrahedron, which is not the case. 

If we move the branch further away then the invariance follows from the  already proven fact that the values of the 1-cocycles do not change if the loop passes through a stratum of $\Sigma^{(1)} \cap \Sigma^{(1)}$, i.e. two simultaneous Reidemeister moves. We use now again the graph $\Gamma$.
The meridian $m$ which corresponds to an arbitrary  edge of $\Gamma$ is a contractible loop in $M_n^{semi-reg}$, no matter what is the position of the branch which moves under everything. Let's take an edge where one vertex is a triple crossing of local type 1, i.e. a positive triple crossing. Reidemeister II moves do not contribute to $R_a^{(2)}$.  Consequently the contribution of the other vertex of the edge doesn't change neither because the contributions from the two  R III moves together sum up to 0.  Using the fact that the graph $\Gamma$ is connected we obtain the invariance with respect to the position of the moving branch for all  local and global types of Reidemeister moves $t$ of type III. 
$\Box$
\vspace{0,2 cm}

Notice that $R_a^{(2)}$ does not have the scan-property for a branch which moves over everything else because the contributions of the strata $+P_3$ and $-\bar P_3$ in the positive tetrahedron equation do not cancel out together at all. But of course one of the "dual" 1-cocycles will have this property.
\vspace{0,4 cm}

\begin{lemma}

$R_a^{(2)}$ is a 1-cocycle  in $M_n$ if and only if for the underlying long knot $v_2(\tilde K)=0$. 

\end{lemma}

{\em Proof.}

R II and R I moves do not enter $R_a^{(2)}$ and therefore it is sufficient to prove that it vanishes for the single R III move in the meridian of the codimension 2 strata which correspond to a cusp with a transverse branch in the projection, compare \cite{F2}. Consequently we have to study the contribution to $R_a^{(2)}$ of sliding a positive small $n$-curl through a crossing with homological marking $a$ or $n$. 
\vspace{0,2 cm}

\begin{figure}
\centering
\includegraphics{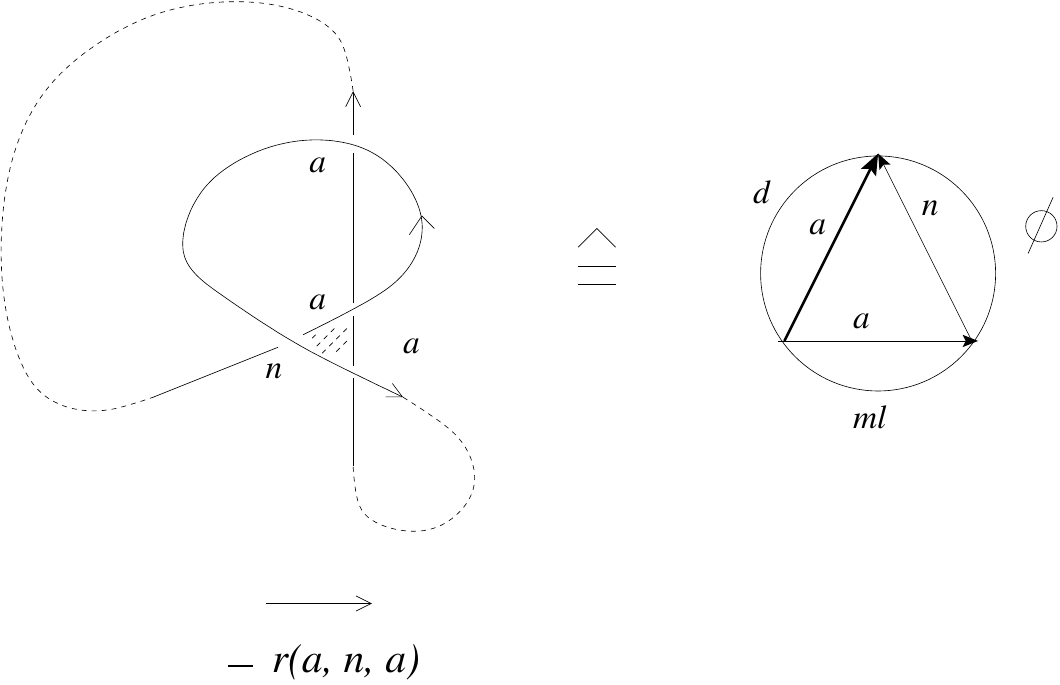}
\caption{\label{nslide} The sliding of an $n$-curl through a crossing of marking $a$}  
\end{figure}

Let us consider the first case. It gives rise to a single R III move $p$, compare Fig.~\ref{nslide}, and $R_a^{(2)}$ is invariant if and only if the contribution of this move vanishes.

We see immediately that $W_2(hm)=0$ because the arc $hm$ is empty (besides perhaps exactly one head or one foot of a crossing $a$, which does not contribute to the weights) and hence the $n$-crossing $hm$ is isolated (i.e. it does not cut any other crossing) and does not contribute. All other $n$-crossings are automatically $f$-crossings because their foots are necessarily in the arc $dml$. It follows that $W_2(p)=v_2(\tilde K)$ for this move, which contributes hence $\pm v_2(\tilde K)$.

Sliding a negative $n$-curl through an $a$-crossing gives the same result because of the Whitney trick and the fact that  $R_a^{(2)}$ is a 1-cocycle.
\vspace{0,2 cm}

Let us now consider the second case. If the new crossing is the crossing $hm$ then the same argument as before applies for $W_2(hm)$. (The only difference is that the crossing which forms with $ml$ a R II move is now an $n$-crossing.) 
Evidently, the $n$-curl can never be the crossing $d$ (because the branch moves over or under it). However, it could be the crossing $ml$, if the branch moves over the small $n$-curl. In that case the region $ml$ is empty and there could well be $0$-crossings  which cut $hm$ and contribute to $W_2(hm)$. However, because the region $ml$ is empty, there can't be any $(n-a)$-crossing with its foot in the region $ml$ and which cuts the crossing $ml$ (there could be just one such $n$-crossing). Consequently, $l(p)=0$ and such a triple crossing $p$ does never contribute to $R_a^{(2)}$ neither (we do not even need here that $v_2(\tilde K)=0$).
$\Box$
\vspace{0,2 cm}

This finishes the proof of Theorem 1.

Institut de Math\'ematiques de Toulouse, UMR 5219

Universit\'e Paul Sabatier

118, route de Narbonne 

31062 Toulouse Cedex 09, France

thomas.fiedler@math.univ-toulouse.fr

\end{document}